\documentclass{fic-l}
\def\currentvolume{56}
\def\currentyear{2009}
\def\copyrightyear{2009}

\usepackage{amssymb}
\usepackage{amsmath}

\usepackage{tikz}

\theoremstyle{definition}

\theoremstyle{remark}

\numberwithin{equation}{subsection}

\def\Aut{\text{Aut}}

\def\b1{\text{\bf 1}}

\def\CA{{\mathcal A}}

\def\CC{{\mathcal C}}
\def\CD{{\mathcal D}}
\def\CF{{\mathcal F}}
\def\CK{{\mathcal K}}
\def\CE{{\mathcal E}}
\def\CH{{\mathcal H}}
\def\CI{{\mathcal I}}

\def\CG{{\mathcal G}}
\def\CL{{\mathcal L}}
\def\CM{{\mathcal M}}
\def\CN{{\mathcal N}}
\def\CO{{\mathcal O}}
\def\CP{{\mathcal P}}
\def\CR{{\mathcal R}}
\def\CS{{\mathcal S}}

\def\CV{{\mathcal V}}

\def\gr{\text{gr}}

\def\Ad{\text{Ad}}

\def\End{\text{End}}
\def\Hom{\text{Hom}}
\def\Sym{\text{Sym}}

\def\deg{\text{deg}}
\def\eff{\text{eff}}

\def\#{\,\check{}}

\def\fD{{\frak D}}

\def\fF{{\frak F}}

\def\Div{{\text{Div}}}

\def\deg{{\text{deg}}}

\def\tr{{\text{tr}}}

\def\id{\text{id}}
\def\tr{\text{tr}}

\def\Coker{\text{Coker}}
\def\Ker{\text{Ker}}

\def\Res{\text{Res}}
\def\Spec{\text{Spec}}

\def\Spf{\text{Spf}}

\def\limleft{\mathop{\vtop{\ialign{##\crcr
  \hfil\rm lim\hfil\crcr
  \noalign{\nointerlineskip}\leftarrowfill\crcr
  \noalign{\nointerlineskip}\crcr}}}}
\def\limright{\mathop{\vtop{\ialign{##\crcr
  \hfil\rm lim\hfil\crcr
  \noalign{\nointerlineskip}\rightarrowfill\crcr
  \noalign{\nointerlineskip}\crcr}}}}

\def\hra{\hookrightarrow}
\def\iso{\buildrel\sim\over\rightarrow}

\def\lra{\longrightarrow}

\begin{document}

\title{$\CE$-Factors for the Period Determinants of Curves}

\author{Alexander Beilinson}

\address{
Department of Mathematics\\
 University of Chicago \\ Chicago, IL 60637, USA}
\email{sasha@math.uchicago.edu}

\dedicatory{To Spencer Bloch}


\maketitle

 \bigskip

\hskip 5cm{\itshape \small The myriad beings of the six worlds -- } 

\hskip 4cm{\itshape \small  gods,  humans, beasts, ghosts, demons, and
devils --} 

\hfill{\itshape \small    are our relatives and
friends.}

\medskip 

\hfill{\itshape \small Tesshu ``Bushido", translated by J.~Stevens.}

\section*{Introduction}

\subsection{} Let $X$ be a smooth compact complex curve, $M$ be a
holonomic $\CD$-module on $X$ (so outside a finite subset $T\subset
X$, our $M$ is a vector bundle with a connection $\nabla$). Denote
by  $dR (M)$ the algebraic  de Rham complex of $M$ placed in degrees
$[-1,0]$; this is a complex of sheaves on the Zariski topology
$X_{\text{Zar}}$. Its analytic counterpart $dR^{an}(M)$ is a complex
of sheaves on the classical topology  $X_{cl}$. Viewed as an object
of the derived category of $\Bbb C$-sheaves, this is a perverse
sheaf, which we denote by $B(M)$; outside $ T$, it is  the local
system $M^\nabla_{X\setminus T}$ of $\nabla$-horizontal sections
(placed in degree $-1$). Set $H^\cdot_{\text{dR}} (X,M) := H^\cdot
(X_{Zar} , dR (M))$, $H^\cdot_{\text{B}} (X,M):= H^\cdot (X_{cl},
B(M))$; these are the {\it de Rham} and {\it Betti} cohomology. We
have the {\it period isomorphism} $\rho : H^\cdot_{\text{dR}}
(X,M)\iso H^\cdot_{\text{B}} (X,M)$.

The cohomology $H^\cdot_{\text{dR}}$ and $H^\cdot_{\text{B}}$ have,
respectively, algebraic and topological nature that can be tasted as
follows.  Let $k,k' \subset \Bbb C$ be  subfields. Then: \newline -
For $(X, M)$  defined over $k$, we have  {\it de Rham $k$-structure} $H^\cdot_{\text{dR}}(X_{k} ,M_{k} )$
  on\linebreak  $H^\cdot_{\text{dR}}(X,M)$;
\newline
- A $k'$-structure on $B(M)$, i.e., a perverse $k'$-sheaf  $B_{k'}$
on $X_{cl}$ together with an isomorphism $B_{k'}
\mathop\otimes\limits_{k'} \Bbb C \iso B(M)$,   yields  {\it Betti
$k'$-structure} $H^\cdot (X_{cl},B_{k'} )$ on $H^\cdot_{\text{B}}
(X,M)$.

If both  $(X_k ,M_k )$ and $B_{k'}$ are at hand, then, computing
$\det\rho$ with respect to rational bases, one gets a number whose class
$[\det\rho ]$ in $ \Bbb C^\times /k^\times k^{\prime\times}$ does not depend on the choice of the bases. In his
farewell seminar at Bures \cite{Del}, Deligne, guided by an analogy
between $[\det\rho ]$  and the constant in the  functional equation
of an $L$-function,  asked if $[\det\rho ]$ can be expressed, in
presence of an extra datum of a rational 1-form $\nu$, as the
product of certain factors of local origin at points of $T$ and
div$(\nu )$. He also suggested the existence of a general geometric
format which would yield the product formula (see 0.3 below).  Our
aim is to establish such a format.

\subsection{}  {\it Remarks.} (i) A natural class of $k'$-structures on $B(M)$
comes as follows. Suppose for simplicity that $M$ equals the
(algebraic) direct image of $M_{X\setminus T}$ by $X\setminus T \hra
X$. Let $\pi: \tilde{X}\to X$ be the real blow-up of $X$ at $T$ (so
$\tilde{X}$ is a real-analytic surface with boundary
$\partial\tilde{X}=\pi^{-1}(T)$, and $\pi$ is an isomorphism over $
X\setminus T$). Then $M^\nabla_{X\setminus T}$ extends uniquely to a
local system $M^\nabla_{\tilde{X}}$ on $\tilde{X}$. Following
Malgrange \cite{M}, consider the constructible subsheaf
$M^\tau_{\tilde{X}} $ of $ M^\nabla_{\tilde{X}}$ of sections of
moderate growth (so $M^\tau_{\tilde{X}} $ coincides with $
M^\nabla_{\tilde{X}}$ off $\partial\tilde{X}$, and
$M^\tau_{\tilde{X}} $ equals $ M^\nabla_{\tilde{X}}$ if and only if
$M$ has regular singularities). By  \cite{M} 3.2, one has a
canonical isomorphism
\begin{equation}  R\pi_*  M^\tau_{\tilde{X}}
\iso B(M). \label{ 0.2.1 }\end{equation} Therefore a $k'$-structure
on $M^\tau_{\tilde{X}}$ yields a $k'$-structure on $B (M)$. Notice
that the former is the same as a $k'$-structure on the local system
$ M^\nabla_{\tilde{X}}$, i.e., on $M^\nabla_{X\setminus T}$, such
that the subsheaf $M^\tau_{\tilde{X}}$ is defined over $k'$.

(ii)  By (0.2.1), one has $H^\cdot_{\text{B}}(X,M)= H^\cdot
(\tilde{X},M^\tau_{\tilde{X}})$.  The dual vector space equals
$H^\cdot (\tilde{X},DM^\tau_{\tilde{X}})$, where $D$ is the Verdier
duality functor, which is the homology group of cycles with
coefficients in $M^\vee_{X\setminus T}$ on $X\setminus T$, having
rapid decay at $T$. So $\rho$, viewed as a pairing
$H^\cdot_{\text{dR}} (X,M)\times H^{\cdot}
(\tilde{X},DM^\tau_{\tilde{X}})\to\Bbb C$, is the  matrix of periods
of $M$-valued forms along the cycles of  rapid decay. See \cite{BE}
for many examples.

(iii)  The setting of 0.1 makes sense for proper $X$ of any
dimension. The passage $B$ to perverse sheaves commutes with direct
image functors for proper morphisms $X\to Y$, so the data $(X_k ,M_k
,B_{k'} )$ are functorial with respect to direct image.

\subsection{}
The next format, which yields the product formula, was suggested in
the last expos\' e of \cite{Del}:

(i) There should exist  $\varepsilon$-factorization formalisms for
$\det H^\cdot_{\text{dR}}$ and $\det H^\cdot_{\text{B}}$. These are
natural rules which assign to every non-zero meromorphic 1-form
$\nu$ on $X$ two collections of  lines $\CE_{\text{dR}}(M)_{(x,\nu
)}$ and   $\CE_{\text{B}}(M)_{(x,\nu )}$   labeled by points $x\in
X$. The lines $\CE_{\text{?}}(M)_{(x,\nu)}$  have $x$-local nature;
if $x\notin T\cup \text{div}(\nu )$, then $\CE_{\text{?}}(M)_{(x,\nu
)}$ is naturally trivialized. Finally, one has  {\it
$\varepsilon$-factorization}, alias {\it product formula},
isomorphisms

\begin{equation}\eta_{\text{?}}:
\mathop\otimes\limits_{x\in T\cup \text{div}(\nu)}
\CE_{\text{?}}(M)_{(x,\nu)} \iso \det H^\cdot_{\text{?}}(X,M).
\label{ 0.3.1}
\end{equation}

(ii)  The de Rham and Betti $\varepsilon$-factorizations should
have, respectively, algebraic and topological origin. Thus, if
$X,M,\nu $ are defined over $k$, then the datum $\{
\CE_{\text{dR}}(M)_{(x,\nu)} \}$  is defined over $k$, and a
$k'$-structure on $B(M)$   yields a $k'$-structure on every
$\CE_{\text{B}}(M)_{(x,\nu)}$. One wants these structures to be
compatible with the trivializations of $\CE_{\text{?}}(M)_{(x,\nu
)}$ off $T\cup \text{div} (\nu )$, and  $\eta_{\text{dR}}$,
$\eta_{\text{B}}$ to be defined over $k$, $k'$.

(iii) There should be natural  {\it $\varepsilon$-period}
isomorphisms $\rho^\varepsilon =\rho^\varepsilon_{(x,\nu )} :
\CE_{\text{dR}}(M)_{(x,\nu ) }\iso \CE_{\text{B}}(M)_{(x,\nu )}$ of
$x$-local origin such that the next diagram commutes:
\begin{equation}
\begin{matrix} \otimes \CE_{\text{dR}}(M)_{(x,\nu )} &
\buildrel{\eta_{\text{dR}}}\over\lra
 &\det H^\cdot_{\text{dR}}(X,M)
 \\
 \otimes\rho^\varepsilon_{(x,\nu )} \downarrow\,\, \,\,\,\,&&\,\,\,\,\,\,\,\,\, \rho\downarrow    \\
\otimes \CE_{\text{B}}(M)_{(x,\nu )} &
\buildrel{\eta_{\text{B}}}\over\lra
 &\det H^\cdot_{\text{B}}(X,M)
\end{matrix} \label{ 0.3.2}\end{equation}

Suppose  $(X,M,\nu )$ is defined over $k$. The points in $T\cup
\text{div} (\nu )$ are algebraic over $k$; let   $\{ O_\alpha \}$ be
their partition by the  Galois orbits. By (ii), the lines
$\mathop\otimes\limits_{x\in O_\alpha} \CE_{\text{dR}}(M)_{(x,\nu
)}$ carry $k$-structure. If $B(M)$ is defined over $k'$, then, by
(ii), $\CE_{\text{B}}(M)_{(x,\nu )}$ carry $k'$-structure.
 Writing
$\mathop\otimes\limits_{x\in O_\alpha}\rho^\varepsilon_{(x,\nu )}$
in $k$-$k'$-bases, we get numbers $[\rho^\varepsilon_{( O_\alpha
,\nu )}]\in \Bbb C^\times /k^\times k^{\prime\times}$. Now (0.3.2)
yields the promised product formula \begin{equation} [\det \rho ]=
\mathop\Pi\limits_\alpha \, [\rho^\varepsilon_{( O_\alpha ,\nu)} ].
\label{ 0.3.3}\end{equation}

 We will show that the above picture is, indeed, true.

\subsection{}
Parts of this  format were established earlier: the de Rham
$\varepsilon$-factorization was constructed already in \cite{Del}
(and  reinvented later in \cite{BBE});  the Betti counterpart was
presented (in the general context of ``animation" of Kashiwara's
index formula) in \cite{B}.\footnote{In \cite{Del} it was suggested
that in case when Re$(\nu )$ is exact, Re$(\nu )=df$, the Betti
$\varepsilon$-factorization comes from the Morse theory of $f$; see
4.6 or \cite{B}  3.8 for a proof.} It remains to construct
$\rho^\varepsilon$. The  point is that $\CE_?$ satisfy several
natural constraints, and compatibility with them determines
$\rho^\varepsilon$ almost uniquely. Notice that we work completely
over $\Bbb C$: the $k$- and $k'$-structures  are irrelevant.

The principal constraints are the global product formula (0.3.1) and
its next local counterpart. For $\nu'$  close to $\nu$, the points
of $T\cup\text{div}(\nu')$ cluster around $T\cup\text{div}(\nu)$.
Now the isomorphism $\mathop\otimes\limits_{x'\in
T\cup\text{div}(\nu')}\CE_{\text{?}}(M)_{(x', \nu')}$ $ \iso
\mathop\otimes\limits_{x\in
T\cup\text{div}(\nu)}\CE_{\text{?}}(M)_{(x,\nu )}$ that comes from
{\it global} identifications (0.3.1) can be written as the tensor
product of natural isomorphisms of {\it  local origin} at points of
$T\cup\text{div}(\nu)$. This  {\it local  factorization structure}
(which is a guise, with an odd twist, of the geometric class field
theory)  is fairly rigid:  $\CE_? (M)$ is determined by a rank 1
local system $\det M_{X\setminus T}$ and a collection of  lines
labeled by elements of $T$.\footnote{In the same manner as the
$\nu$-dependence of the classical  $\varepsilon$-factor of a Galois
module $V$  is controlled, via the class field theory,  by $\det
V$.}

 The rest of constraints for  $M\mapsto \CE_? (M)$ are listed in 5.1. We show that there is an isomorphism
$\rho^\varepsilon :\CE_{\text{dR}}\iso \CE_{\text{B}}$ compatible
with them, which is determined uniquely up to a power of a simple
canonical automorphism of $\CE_? $, i.e., $\rho^\varepsilon$ form
 a $\Bbb Z$-torsor $E_{\text{B}/\text{dR}}$. First we recover 
 $\rho^\varepsilon$ from   $\eta$-compatibility (0.3.2)
 for $(\Bbb P^1,
\{0,\infty\})$, $M$ with regular singularities at $\infty$, and
$(\Bbb P^1, \{0,1,\infty\})$, $M$ of rank 1 with regular
singularities. Having $\rho^\varepsilon$ at hand, one has to prove
that it is compatible with the constraints for all  $(X,T,M)$, of
which  (0.3.2) is central. The core  of the argument is  global: we
use a theorem of Goldman \cite{G} and Pickrell-Xia \cite{PX1},
\cite{PX2}, which asserts that the action of the Teichm\" uller
group on the moduli space of unitary local systems with fixed local
monodromies is ergodic. As in \cite{G}, this implies that, when the genus 
of $X$ and the order of $T$ are fixed, the
possible discrepancy of (0.3.2) depends only on the {\it local}
datum of  monodromies at singularities of $M$. An observation that
this discrepancy  does not change upon quadratic degenerations of
$X$ reduces the proof to a few simple computations.

\subsection{}
One can ask for an explicit formula for $\rho^\varepsilon (M)$. An
analytic approach as in  \cite{PS} or \cite{SW} shows that the de
Rham  $\varepsilon$-factors of a $\CD$-module $M$ can be recovered
from the  $\CD^\infty$-module $M^\infty :=\CD^\infty
\mathop\otimes\limits_{\CD}M$. Thus the ratio between the
$\varepsilon$-factors of $M$ and of another $\CD$-module $M'$ (say,
with regular singularities) with $B(M)=B(M')$,  is certain Fredholm
determinant (a variant of $\tau$-function). If $x$ is a regular
singular point of $M$, then $[\rho^\varepsilon_{(x,\nu)} ]$   can be
written explicitly using the $\Gamma$-function, see 6.3 (which is
similar to the fact that the classical $\varepsilon$-factors of
tamely ramified Galois modules are essentially products of
Gau\ss~sums). An example of the product formula is the Euler
identity $\int\limits_0^1 t^{\alpha -1}(1-t)^{\beta
-1}dt=\frac{\Gamma (\alpha )\Gamma (\beta ) }{\Gamma (\alpha +\beta
)}.$

\subsection{} Plan of the article:
\S 1 presents a general story of factorization lines (i.e., of the
local factorization structure); in \S 2--4 the  algebraic and
analytic de Rham $\varepsilon$-factors, and their Betti counterpart
are defined;   \S 5 treats the  $\varepsilon$-period map; in \S 6
the $\varepsilon$-periods are written explicitly in terms of the
$\Gamma$-function.

A different approach to product formula (0.3.3),  based on Fourier
transform, was developed by Bloch, Deligne, and Esnault  \cite{BDE},
\cite{E}  (some essential ideas go back to  \cite{Del} and
\cite{L}; the case of regular singularities was considered earlier,
and for $X$ of arbitrary dimension, in \cite{A}, \cite{LS},
\cite{ST}, and \cite{T}).\footnote{ \cite{BDE} considers the case of
$M$ of virtual rank 0 and the Betti structure  compatible with the
Stokes structures (hence of type considered in 0.3(i)).}  The two
constructions are fairly complementary; the relation between them
remains to be understood.

{\it Questions \& hopes.} ($o$) For Verdier dual $M$, $M^\vee$ the
lines $\CE_? (M)_{(x,\nu )}$ and $\CE_? (M^\vee )_{(x,-\nu )}$
should be naturally dual,\footnote{For $\CE_{\text{B}}$ this is evident from the construction; for $\CE_{\text{dR}}$
one can hopefully deduce it from (2.10.5) applied to $M\oplus
M^\vee$.}  and $\rho^\varepsilon$ should be compatible with duality.

(i)  The period story should exist for $X$ of any dimension,
  with mere lines replaced by finer objects (the homotopy points
of $K$-theory spectra). For the Betti side, see \cite{B}; for the de
Rham one, see \cite{P}.

The meaning of local factorization structure for $\dim X >1$ is
not clear (as of the more general  notion of factorization sheaves
in the setting of algebraic geometry). Is there an 
agebro-geometric analog of the recent beautiful work of 
Lurie on the classification of TQFT?

(ii) There should be a geometric theory of $\varepsilon$-factors
(cf.~5.1) for \' etale sheaves; for  an \'etale sheaf of virtual
rank 0  on  a curve over a finite  field, the corresponding trace of Frobenius 
function  should be equal to the classical
$\varepsilon$-factors.\footnote{The condition of virtual rank 0 is
essential:   the $\varepsilon$-line for the constant  sheaf of rank
1 has non-trivial $\pm 1$ monodromy on the components of $\nu$'s
with odd order of zero, so the trace of Frobenius function is
non-constant on every such component (as opposed to the classical
$\varepsilon$-factor).} Notice that Laumon's construction \cite{L}
(which is the only currently available method to establish the
product formula for classical $\varepsilon$-factors) has different
arrangement: its input is more restrictive (the forms $\nu$ are
exact), while the output is more precise (the  $\varepsilon$-lines
are realized as determinants of true complexes).

(iii) What would be a motivic version of the story?

(iv)    $\Gamma$-function  appears in  Deninger's vision \cite{Den}
of classical local Archimedean $\varepsilon$-factors. Are the two
stories related on a deeper level?

\subsection{}  I am  grateful to S.~Bloch, V.~Drinfeld, and
H.~Esnault  whose  interest was crucial for this work, to P.~Deligne
for the pleasure to play in a garden he conceived,  to B.~Farb for
the information about the Goldman and Pickrell-Xia theorems,  to
V.~Schechtman and D.~Zagier who urged me to write formulas, to the
referee for the help, and to IHES for a serene sojourn. The
research was partially supported by NSF grant DMS-0401164.

The article is a modest tribute to Spencer Bloch, for all his gifts
and  joy  of books,  of the woods, and of our  relatives and friends
-- the numbers.

\section{Factorization lines}

This section is essentially an exposition of geometric class field
theory (mostly) in its  algebraic de Rham version.

\subsection{} We live  over a fixed ground field $k$ of characteristic 0;
``scheme" means ``separated $k$-scheme of finite type". The category
$\CS ch$ of schemes is viewed  as a site for the \'etale topology
(so ``neighborhood" means ``\'etale neighborhood", etc.),  ``space"
means a sheaf on $\CS ch$; for a space $F$ and a scheme $S$ elements
of $F(S)$ are referred to as $S$-points of $F$.  All Picard groupoids
are assumed to be commutative and essentially small. For a Picard
groupoid $\CL$, we denote by  $\pi_0 (\CL )$, $\pi_1 (\CL )$ the
group of isomorphism classes of its objects and the automorphism
group of any its object; for $L\in \CL$ its class is  $[L]\in\pi_0
(\CL )$.

Let $X$ be a smooth (not necessary proper or connected) curve, $T$
its finite subscheme,  $K$  a line bundle on $X$.\footnote{Starting
from \S2, our $K$ equals $\omega_X$.}  For a test scheme $S$, we
write $X_S := X\times S$, $T_S := T\times S$, $K_S := K\boxtimes
\CO_S$; $\pi : X_S \to S$ is the projection. For a Cartier divisor
$D$ on $X_S$ we denote by $|D|$ the support of $D$ viewed as a {\it
reduced} closed subscheme.

Consider the next spaces:
\newline (a) $\Div (X)$: its $S$-points are relative Cartier divisors
$D$ on $X_S /S$ such that $|D|$ is finite over $S$;
\newline  (b)  $2^T$ is a scheme whose $S$-points $c$ are idempotents
in $\CO (T_S )$. Such  $c$ amounts to an open and closed subscheme
$T_S^c$ of $T_S$ (the support of $c$);
\newline (c) $\fD =\fD (X,T)\subset \Div (X)\times 2^T$ consists of those pairs $(D,c )$ that $D\cap T_S \subset T_S^c$;
\newline (d) $\fD^\diamond = \fD^\diamond (X,T;K)$ is formed by triples $(D,c,\nu_P )$
 where $(D,c)\in \fD$ and $\nu_P$ is a
 trivialization of the restriction  of the line bundle $K(D):=K_S (D)$ to the subscheme $P=P_{D,c}:=T^c_S \cup |D|$.

Denote by $\pi_0 (X)$ the scheme of connected components of
$X$.\footnote{Which is the spectrum of the integral closure of $k$
in the ring of functions on $X$.}  One has  projection $\deg :  \Div
(X) \to \Bbb Z^{\pi_0 (X)}$, hence the projections $\fD^\diamond
\to\fD \to \Bbb Z^{\pi_0 (X)}\times 2^T$. Notice that the component
$\fD_{c=0}$ equals $\Div (X\setminus T)$, and   $\fD_{c=1} $ equals
$\Div (X)$.

{\it Remarks.} (i) Every $S$-point of $\fD$ can be lifted
$S$-locally to $\fD^\diamond$.

 (ii) Every $\nu_P$ as in (d) can be extended $S$-locally to a
trivialization $\nu$ of $K (D)$ on a neighborhood $V\subset X_S$ of
$P$. One can view $\nu_P$ as an equivalence class of $\nu$'s, where
$\nu$  and $\nu'$  are equivalent if the function $\nu/\nu'$ equals
1 on $P$. We often write $(D,c,\nu )$ for $(D,c,\nu_P )$.

(iii) Each space $F$ of the list (a)--(d) is smooth in the next
sense: for every closed embedding $S\hra S'$, a geometric point
$s\in S$, and $\phi \in F(S)$ one can find an neighborhood $U'$ of
$s$ in $S'$ and $\phi'\in F(U')$ such that
$\phi|_{U'_S}=\phi'|_{U'_S}$.

(iv)  The geometric fibers of $\Div (X)$ over $ \Bbb Z^{\pi_0 (X)}$
and  of $\fD  $, $\fD^\diamond$ over $\Bbb Z^{\pi_0 (X)}\times 2^T$,
are connected (i.e., every two geometric points of any fiber are
members of one connected family).

A comment about the fiber
 $\fD^\diamond_{(D,c)}$ of $\fD^\diamond / \fD$ over $(D,c)\in \fD (S)$: Suppose $S$ is smooth, so
 $P_{D,c}$ is a relative Cartier divisor in $X_S /S$.
Denote by $O^\times_{D,c}$ the Weil $P_{D,c}/S$-descent of $\Bbb
G_{m P}$, and by $K(D)^\times_{D,c}$ the Weil $P_{D,c}/S$-descent of
the $\Bbb G_{m P}$-torsor of trivializations of the line bundle
$K(D)|_{P_{D,c}}$. Then $O^\times_{D,c}$ is a smooth group
$S$-scheme, and $K(D)^\times_{D,c}$ is an $O^\times_{D,c}$-torsor;
for any $S$-scheme $S'$ an $S'$-point of $K(D)^\times_{D,c}$ is the
same as a trivialization of $K(D)$ over $(P_{D,c})_{S'}$. The latter
relative divisor contains $P_{D_{S'}, c_{S'}}$ (the corresponding
reduced schemes coincide), so we have a canonical surjective
morphism $K(D)^\times_{D,c}\to \fD^\diamond_{(D,c)}$, hence a
canonical $(D,c,\nu_P )\in \fD^\diamond (K(D)^\times_{D,c})$. The
map $K(D)^\times_{D,c}(S')\to \fD^\diamond_{(D,c)}(S')$ is bijective
if $S'$ is smooth over $S$, but {\it not} in general.

{\it Examples.} (i)  Suppose $S=X\setminus T$, $(D,c)= (\ell \Delta
,0)$ where $\Delta$ is the diagonal divisor, $\ell$ is any integer.
Then $O^\times_{D,c}=\Bbb G_{m S}$, and $K(D)^\times_{D,c}$ is the
$\Bbb G_{m}$-torsor $\CK^{(\ell )}$ of trivializations of the line
bundle $K(D)|_\Delta  = K\otimes \omega_X^{\otimes -\ell}|_S$. For
any $S'/S$ an $S'$-point of $\fD^\diamond_{(\ell\Delta ,0)}$ is the
same as an $S'_{\text{red}}$-point of $\CK^{(\ell )}$, i.e.,
$\fD^\diamond_{(\ell\Delta ,0)}$ is the quotient of $\CK^{(\ell )}$
modulo the action of the formal multiplicative group $\Bbb
G\hat{_{m}}$.

(ii) For a point $b\in T$ let $k_b$ be its residue field,
$T_b\subset T$ be the component of $b$, and $m_b$ its multiplicity.
Consider $(D,c)=(nb, 1_b )\in\fD (S )$, where   $S=\Spec\, k_b$, $n$
is any integer, $1_b$ is the characteristic function of $b\in T(S)$.
Then $P_{D,c}=T_b$, so $O^\times_{D,c}=O^\times_{T_b}$ is an
extension of $\Bbb G_{m S}$ by the unipotent radical. One has
$K(nb)^\times_{T_b}:= K(D)^\times_{D,c}\iso \fD^\diamond_{(D,c)}$.
We set $K^\times_{T_b}:= K(0b)^\times_{T_b}$.

\subsection{} Let $\CV$ be a stack, alias a sheaf of categories,
on $\CS ch$. For a space $F$ we denote by $\CV (F)$  the category of
Cartesian functors $V : F\to \CV$. Explicitly, such $V$ is a rule
that assigns to every test scheme $S$ and $\phi \in F(S)$ an object
$V_\phi \in\CV (S)$ together with a base change compatibility
constraint. If $\CV$ is a Picard stack, alias  a sheaf of  Picard
groupoids, then $\CV (F)$ is naturally a Picard groupoid.

Below we denote by $\bar{F}$  the space with $\bar{F}(S):=
F(S_{\text{red}})$. The stack of {\it $\CV$-crystals}
$\CV_{\text{crys}}$ is defined by formula
 $ \CV_{\text{crys}}(S):= \CV (\bar{S})$.
If $F$ is formally smooth (i.e., satisfies the property from Remark
(iii) in 1.1 for every  nilpotent embedding $S\hra S'$), then
$\bar{F}$ is the quotient of $F$ modulo the evident equivalence
relation; therefore objects of $\CV_{\text{crys}} (F)= \CV
(\bar{F})$ are the same as objects $V\in \CV (F)$ equipped with a
{\it de Rham structure}, i.e., a natural identification $\alpha :
V_\phi \iso V_{\phi'}$ for every $\phi ,\phi' \in F(S)$ such that
$\phi|_{S_{\text{red}}}=\phi'|_{S_{\text{red}}}$, which is
transitive and compatible with base change. E.g.~if $F$ is a smooth
scheme, then a vector bundle crystal on $F$ is the same as a vector
bundle on $F$ equipped with a flat connection.

Key examples:   Let $\CL_k$ be the Picard groupoid of $\Bbb
Z$-graded $k$-lines (with ``super"  commutativity constraint for the
tensor structure). Below we call them simply ``lines" or
``$k$-lines"; the degree of a line $\CG$ is denoted by deg$(\CG )$.
An {\it $\CO$-line on} $S$ (or {\it $\CO_S$-line}) is  an invertible
$\Bbb Z$-graded vector bundle on $S$. These objects form a Picard
groupoid $\CL_\CO (S)$; the usual pull-back functors make $\CL_\CO$
a Picard stack. Below $\CL_\CO$-crystals are referred to as {\it de
Rham lines}; they form a Picard stack $\CL_{\text{dR}}$. Instead of
$\Bbb Z$-graded lines, we can consider $\Bbb Z/2$-graded ones; the
corresponding Picard stacks are denoted by $\CL'_\CO$,
$\CL'_{\text{dR}}$. We mostly consider  $\Bbb Z$-graded setting; all
the results remain valid, with evident modifications, for $\Bbb
Z/2$-graded one.

 {\it Remarks.}
 (i) Let $(X',T')$ be another pair as in 1.1, and $\pi : (X',T')\to (X,T)$
be a finite morphism of pairs, i.e., $\pi : X'\to X$ is a finite
morphism of curves such that $\pi (T')\subset T$. It yields a
morphism of spaces $\pi^* : \fD^\diamond (X,T ;K)\to \fD^\diamond
(X',T' ;\pi^* K)$, $(D,c,\nu )\mapsto (\pi^* D,\pi^* c, \pi^* \nu
)$, hence the pull-back functor $\pi_* : \CV (\fD^\diamond
(X',T';\pi^* K ))\to \CV (\fD^\diamond (X,T;K ))$ denoted by
$\pi_*$; if $\CV$ is a Picard stack, then $\pi_*$ is a morphism of
Picard groupoids. If $X'=X$, $T'\subset T$, we refer to $\pi_*$ as
``restriction to $(X,T)$".

{\it Exercise.} If $T'\subset T$,  $T'_{\text{red}}=T_{\text{red}}$,
then the restriction  $\CL_{\text{dR}} (\fD^\diamond (X,T';K ))\to
\CL_{\text{dR}} (\fD^\diamond (X,T;K ))$  is a fully faithful
embedding.

We denote the union of the Picards groupoids $\CL_{\text{dR}}
(\fD^\diamond (X,T'';K ))$ for all $T''$ with
$T''_{\text{red}}=T_{\text{red}}$ by $\CL_{\text{dR}} (\fD^\diamond
(X,\hat{T};K ))$  (here $\hat{T}$ is the formal completion of $X$ at
$T$).

  (ii) The space $\fD^\diamond (X,T;K)$, hence
    $\CL_? (\fD^\diamond (X,T;K))$,  actually depends only on the
restriction of $K$ to $X\setminus T$.  Indeed, for  any divisor
$D_{(T)}$ supported on $T$, there is a canonical identification
$\fD^\diamond (X,T;K)\iso \fD^\diamond (X,T;K(D_{(T)}))$, $(D,c, \nu
)\mapsto (D- D^c_{(T)},\nu )$, where $D^c_{(T)}$ equals $D_{(T)}$ on
$T_S^c$ and to 0 outside. We keep $K$ to be a line bundle on $X$ for
future notational convenience.

(iii)  If $U$ is any open subset of $X$, then $\fD^\diamond (U,T_U
;K_U )\subset   \fD^\diamond (X,T ;K)$, hence we have the
restriction functor $\CV (\fD^\diamond (X,T;K ))\to \CV
(\fD^\diamond (U,T_U ;K_U ))$.

(iv)  Remark (iv) in 1.1 implies that $\pi_1
(\CL_{\text{dR}}(\fD^\diamond )=\CO^\times (\Bbb Z^{\pi_0 (X)}\times
2^{T})$.

\subsection{} Let $\CS m \subset \CS ch$ be the full subcategory of smooth schemes. For $\CV$, $F$ as in 1.2
we denote by $\CV^{\text{sm}} (F )$ the Picard groupoid of Cartesian
functors $F|_{\CS m}\to \CV |_{\CS m}$. One has a restriction
functor $\CV (F)\to \CV^{\text{sm}} (F)$. If $F$ is smooth in the
sense of Remark (iii) in 1.1, then this is a faithful functor.

{\it Exercise.}
 Suppose we have $\CE ,\CE' \in \CL_{\text{dR}} (\fD^\diamond )$
and a morphism $\phi : \CE \to \CE'$ in  $\CL_{\CO} (\fD^\diamond
)$. Then $\phi$  is a morphism in $\CL_{\text{dR}} (\fD^\diamond )$,
if (and only if) the corresponding morphism in $\CL^{sm}_{\CO}
(\fD^\diamond)$ lies in $\CL^{sm}_{\text{dR}} (\fD^\diamond
)$.\footnote{Hint: Use Remark (iii) in 1.1 for embeddings $S\hra S'$
where $S'$ is smooth.}

\bigbreak\noindent\textbf{Lemma.}\emph{
 $\CL_{\text{dR}} (F )\iso
\CL^{\text{sm}}_{\text{dR}} (F)$.}\bigbreak

{\it Proof.} This follows from the fact that $\CL_{\text{dR}}$ is a
stack with respect to the h-topology, and h-locally every scheme is
smooth.  \hfill$\square$

{\it Remark.} By the lemma and  1.1, one can view $\CE\in
\CL_{\text{dR}} (\fD^\diamond )$ as a rule that assigns to every
smooth $S$ and $(D,c)\in\fD (S)$ a de Rham line $\CE_{(D,c)}:=
\CE_{(D,c,\nu_P )}$ on $K(D)_{D,c}^\times$ in a way compatible with
the base change.

\subsection{} {\it For this subsection, $X$ is proper.}
Let Rat$(X,K)=$ Rat$(X)$ be a space whose $S$-points are rational
sections $\nu$ of the line bundle $K_S $ such that $|\text{div}(\nu
)|$ does not contain a connected component of any geometric fiber of
$X_S /S$. There is a natural morphism Rat$(X)\to
\fD^\diamond_{c=1}$, $\nu \mapsto (-\text{div}(\nu ), 1, \nu )$, so
every $\CE\in\CL_? (\fD^\diamond )$ yields naturally an object of
$\CL_? (\text{Rat}(X))$, which we denote again by $\CE$.

 The next fact is a particular case of \cite{BD}  4.3.13:

\bigbreak\noindent\textbf{Proposition.} \emph{Every function on
Rat$(X)$ is constant.
 All $\CO$- and de Rham lines on Rat$(X)$ are constant.}\bigbreak

{\it Proof.}  Let  $L$ be an auxiliary ample line bundle on $X$; set
$V^{(m)}_1 :=  \Gamma (X, K\otimes L^{\otimes m})$, $V^{(m)}_2 :=
\Gamma (X,  L^{\otimes m})$. Let $U^{(m)}\subset \Bbb P ( V^{(m)}_1
\times V^{(m)}_2 )$ be the open subset of those $\phi = (\phi_1
,\phi_2 )$ that neither $\phi_1 $ nor $\phi_2 $  vanishes on any
connected component of $X$. Consider a map $\theta^{(m)}: U^{(m)}\to
\text{Rat}(X)$, $(\phi_1 ,\phi_2 )\mapsto \phi_1 /\phi_2$. We will
check that for $m$ large the $\theta^{(m)}$-pull-back of any $\CO$-
or de Rham line on Rat$(X)$ is trivial, and every function on $
U^{(m)}$ is constant. This implies the proposition, since  geometric
fibers of $\theta^{(m)}$ are connected and the images of $
\theta^{(m)}$ form a directed system of subspaces  whose inductive
limit equals Rat$(X)$ (i.e., every $\nu\in$ Rat$(X) (S)$ factors
$S$-locally through $\theta^{(m)}$   for sufficiently large $m$).

 Notice that  the complement to $U^{(m)}$ in $ \Bbb P (
V^{(m)}_1 \times V^{(m)}_2 ) $ has codimension $\ge 2$ for $m$
large. Therefore every function and every $\CO$- and de Rham line
extend to $ \Bbb P ( V^{(m)}_1 \times V^{(m)}_2 )$. Thus for $m$
large every function on $U^{(m)}$ is constant and every de Rham line
is trivial.

The case of an $\CO$-line requires an extra argument.  Any $(\psi_1
,\psi_2) \in U^{(n)}$ yields an embedding $\Bbb P (\Gamma
(X,L^{\otimes k}))\hra  \Bbb P ( V^{(m)}_1 \times V^{(m)}_2 )$,
$\gamma \mapsto (\gamma \psi_1 ,\gamma\psi_2 )$; here $m=n+k$. For
$k$ large, the preimage of $U^{(m)}$ in   $\Bbb P (\Gamma
(X,L^{\otimes k}))$ is an open dense subset of codimension $\ge 2$,
 and $\theta^{(m)}$ is constant on it. Thus  for any $\CO$-line $\CL$ on Rat$(X)$ the restriction of the corresponding line on $ \Bbb P (
V^{(m)}_1 \times V^{(m)}_2 )$ to
 $\Bbb P (\Gamma (X,L^{\otimes k}))$ is trivial, hence the  line itself is  trivial, so
  $\theta^{(m)*}\CL $  is trivial, and we are done.
  \hfill$\square$

Therefore for any $\CE$ in $\CL_{\CO} (\fD^\diamond )$ or
$\CL_{\text{dR}} (\fD^\diamond )$ the lines $\CE_\nu $ for all
rational non-zero $\nu$ are canonically identified. We denote this
line simply by $\CE (X)$.

\subsection{}
A finite subset of $\{ (D_\alpha ,c_\alpha ,\nu_\alpha)\} $ of
$\fD^\diamond (S)$, is said to be {\it disjoint} if the subschemes
$P_{D_\alpha ,c_\alpha}$ are pairwise disjoint. Then we have
$\Sigma (D_\alpha ,c_\alpha ,\nu_{\alpha}):=(\Sigma D_\alpha ,\Sigma
c_\alpha ,\Sigma \nu_{\alpha} )$ $\in \fD^\diamond (S)$, where
$\Sigma \nu_{\alpha}$ equals $\nu_\alpha$ on   $P_{D_\alpha
,c_\alpha}$.

For $\CE$ in $\CL_? (\fD^\diamond )$, where   $\CL_?$ is a Picard
stack, a {\it factorization structure} on $\CE$ is a rule which
assigns to every disjoint family as above a {\it factorization
isomorphism} \begin{equation}\otimes_\alpha \, \CE_{(D_\alpha
,c_\alpha ,\nu_{\alpha} )}\iso \CE_{\Sigma (D_\alpha ,c_\alpha
,\nu_{\alpha}) }. \label{ 1.5.1}\end{equation} These isomorphisms
should be compatible with base change and satisfy an evident
transitivity property. One defines factorization structure on
objects of $\CL_? (\fD (X,T) )$,  $\CL_? (\Div(X))$,  $\CL_?
(\Div^{\eff}(X))$, or  $\CL_? (2^T )$ in the similar way.

Objects of $\CL_? (\fD^\diamond )$ equipped with a factorization
structure are called {\it ($K$-twisted) factorization objects} of
$\CL_?$ on $(X,T;K)$; they form a Picard groupoid
$\CL_?^\Phi(X,T;K)$. In particular, we have Picard groupoids
$\CL^\Phi_{\CO} (X,T;K)$, $\CL^\Phi_{\text{dR}} (X,T;K)$ of $\CO$-
and de Rham factorization lines.

\bigbreak\noindent\textbf{Proposition.}\emph{
  Factorization objects have local nature:
$U\mapsto  \CL_?^\Phi (U,T_U; K_U)$ is  a  Picard stack
 on $X_{\text{\'et}}$. }\bigbreak

{\it Proof.}   Let $\pi : U\to X$ be an \'etale map. For $\CE \in
\CL_?^\Phi (X,T;K )$ one defines its pull-back $\pi^* \CE$ as
follows. Take any $(D,c,\nu )\in \fD^\diamond (U, T_U ; K_U )$. It
suffices to define $\CE_{(D,c,\nu )}$  \'etale locally on $S$. Write
$(D,c,\nu )=\Sigma (D_\alpha ,c_\alpha ,\nu_\alpha )$ with connected
$P_\alpha = P_{D_\alpha ,c_\alpha}$. Then there is a uniquely
defined $(D'_\alpha ,c'_\alpha ,\nu'_\alpha ) \in \fD^\diamond (X, T
; K )$ such that $D_\alpha$ is a connected component of the
pull-back of $D'_\alpha$ to $U$ and $\pi$ yields an isomorphism
$P_\alpha \iso P'_\alpha$ which identifies $\nu_{\alpha P_\alpha}$
with $\nu'_{\alpha P'_\alpha}$.

Set $\pi^* \CE_{(D,c,\nu )}:= \otimes\, \CE_{(D'_\alpha ,c'_\alpha
,\nu'_\alpha )}$. Due to factorization structure on $\CE$, this
definition is compatible with base change, and $\pi^* \CE
\in\CL_?^\Phi (\fD^\diamond (U, T_U; K_U ))$ so defined has an
evident factorization structure. Thus $\CL_?^\Phi$ is a presheaf of
Picard groupoids on $X_{\text{\'et}}$. We leave it to the reader to
check the gluing property. \hfill$\square$

NB:  The pull-back functor for open embeddings  is defined
regardless of factorization structure (see Remark (iii) in 1.2).

 {\it Remarks.}
(i) The evident forgetful functor $\CL_{\text{dR}}^\Phi (X,T;K) \to
\CL_{\CO}^\Phi (X,T;K)$ is faithful. By 1.2 (and Remark (iii) in
1.1), for $\CE \in \CL_{\CO}^\Phi (X,T;K)$ a de Rham structure on
$\CE$, i.e., a lifting of $\CE$ to $\CL_{\text{dR}}^\Phi (X,T;K)$,
amounts to a rule which assigns to every scheme $S$ and a pair of
points $(D,c,\nu_P ), (D',c',\nu'_P )\in \fD^\diamond (S)$ which
coincide on $S_{\text{red}}$,
    a natural identification (notice that $c=c'$)
\begin{equation}
\alpha^\varepsilon : \CE_{ (D,c,\nu_P )}\iso \CE_{ (D',c,\nu'_P )}.
\label{ 1.5.2}
\end{equation}
The $\alpha^\varepsilon$ should be
transitive and compatible with base change and  factorization.

(ii)  Remarks  in 1.3 and (i)--(iii) in 1.2  remain valid for
factorization lines.  Thus we have a  Picard groupoid
$\CL^\Phi_{\text{dR}} (X,\hat{T};K) $, etc.

(iii) There is a natural Picard functor  \begin{equation}
\mathop\Pi\limits_{b\in T_{\text{red}}} \CL_? (b) \to \CL^\Phi_?
(X,T;K), \label{ 1.5.3}\end{equation} which assigns to $E= (E_b
)\in\Pi \CL_? (b) $ a factorization object $\CE$ with
$\CE_{(D,c,\nu_P )}= \text{Nm}_{ T^c_{\text{red}S} /S}(E)$; here
$T^c_{\text{red}S}$ is the preimage of $T_{\text{red}}$ by the
projection $p: T^c_S \to T$.

(iv)  By Remark (iv) in 1.2 there is a natural isomorphism
\begin{equation}\CO^\times (T_{\text{red}})\times\CO^\times (\pi_0
(X))\iso \pi_1 (\CL_{\text{dR}}^\Phi (X,T;K)). \label{
1.5.4}\end{equation} Here $(\alpha ,\beta )\in \CO^\times
(T_{\text{red}})\times \CO^\times (\pi_0 (X))$ acts on
$\CE_{(D,c,\nu )}$ as multiplication by the locally constant
function $\text{Nm}_{  T^c_{\text{red}S} /S}(\alpha
)\text{Nm}_{\pi_0 (X)_S /S}(\beta^{deg (D)})$. Notice that the
embedding $\CO^\times (T_{\text{red}})\hra \pi_1
(\CL_{\text{dR}}^\Phi (X,T;K))$ comes from (1.5.3).

\subsection{}
As in 1.1, every $(D,c)\in \fD (S)$, $S$ is smooth, yields a
morphism  $K(D)^\times_{D,c} \to \fD^\diamond$, hence a Picard
functor $ \CL_?^\Phi (\fD^\diamond )\to\CL_? (K(D)_{D,c}^\times )$,
$\CE \mapsto  \CE_{(D,c)}$. In particular, following Examples  in
1.1, for $\ell \in \Bbb Z$ we have $\CE^{(\ell )}:= \CE_{(\ell\Delta
,0)}\in \CL_? (\CK^{(\ell )})$, and for $b\in T$, $n\in\Bbb Z$, we
have $\CE^{(n)}_{T_b}:=  \CE_{(nb,1_b )} \in \CL_?
(K(nb)_{T_b}^\times )$. Set $\CE_{T_b}:=\CE^{(0)}_{T_b}\in   \CL_?
(K_{T_b}^\times )$. Notice that $\CE^{(0)}\in  \CL_? (\CK^{(0 )}) $
is canonically trivialized.

If $\CE\in \CL_{\CO}^\Phi (\fD^\diamond )$, then the $\CO$-lines
$\CE^{(\ell )}$ carry a canonical connection along the fibers of the
projection $\CE^{(\ell )}\to X\setminus T$ (see Example (i) in 1.1).
A de Rham structure on $\CE$ provides a flat connection
$\nabla^\varepsilon$ on $\CE^{(\ell )}$ that extends this relative
connection. Since the degrees of lines are locally constant, the
factorization implies  \begin{equation} \deg (\CE^{(\ell )})=
\ell\,\deg(\CE^{(1)}), \quad \deg ( \CE^{(n+1)}_{T_b})=\deg (
\CE^{(n)}_{T_b})+\deg (\CE^{(1)}). \label{ 1.6.1}\end{equation}

Let $\CL_{\text{dR}} (X,T)\subset \CL_{\text{dR}} (X\setminus T)$ be
the Picard subgroupoid of those de Rham lines whose connection at
every $b\in T_{\text{red}}$ has pole of order less or equal the
multiplicity of $T$ at $b$.

For a $\Bbb G_m$-torsor $\CK$ over a  scheme $Y$ we denote by
$\CL_{\text{dR}} (Y;\CK )$  the Picard groupoid of   de Rham lines
$\CG$ on $\CK$  such that $\CG^{\otimes 2}$ is constant along the
fibers (i.e., comes from a de Rham line on $Y$) and the fiberwise
monodromy of $\CG$ equals  $(-1 )^{\text{deg} \CG }$. Thus  we have
Picard groupoids $\CL_{\text{dR}} (X\setminus T )^{(\ell)}:=
\CL_{\text{dR}} (X\setminus T; \CK^{(\ell)})$, $\ell \in \Bbb Z$;
let $\CL_{\text{dR}} (X,T)^{(\ell )} $ $ \subset \CL_{\text{dR}}
(X\setminus T )^{(\ell)} $ be the Picard subgroupoids of those $\CG$
that $\CG^{\otimes 2}\in\CL_{\text{dR}} (X,T)$.

Choose a trivialization $\nu_T$ of the restriction of $K$ to $T$,
i.e., a collection $\{ \nu_{T_b} \}$ of $k_b$-points in
$K^\times_{T_b}$.
 For a factorization line $\CE$ set $\CE_{\nu_{T_b}}:=\CE_{(0 ,1_b, \nu_{T_b} )}=$ the fiber of  $\CE_{T_b}$ at $\nu_{T_b}$.
The next theorem is the main result of this section. The proof for
$T=\emptyset$  is  in 1.7--1.9; the general case is treated in
1.10--1.11.  In 1.11 one finds its reformulation free from the
auxiliary $\nu_T$.

\bigbreak\noindent\textbf{Theorem.}\emph{ For  $\CE \in
\CL_{\text{dR}}^\Phi (X,T;K)$ one has $\CE^{(1)}\in\CL_{\text{dR}}
(X, T )^{(1)}$, and the functor
\begin{equation}\CL_{\text{dR}}^\Phi (X,T;K)\to \CL_{\text{dR}} (X,
T )^{(1)}\times\! \mathop\Pi\limits_{b\in T_{\text{red}}} \!\!
\CL_{k_b},\quad \CE \mapsto (\CE^{(1)},\{ \CE_{\nu_{T_b} }\}),
\label{ 1.6.2}\end{equation} is an equivalence of Picard groupoids.
}\bigbreak

If $K=\omega_X$, then the $\Bbb G_m$-torsor $\CK^{(1)}$ is
trivialized by a canonical section $\nu_1$ (its value at $x\in
X\setminus T$  is the element  in $\omega (x)/\omega$
 with residue 1). The
functor $\nu_1^* : \CL_{\text{dR}} (X, T )^{(1)}\to
\CL_{\text{dR}}(X, T )$ is evidently an equivalence, so the theorem
can be reformulated as follows (here  $\CE^{(1)}_{X\setminus T}:=
\nu_1^* \CE^{(1)}= \CE_{(\Delta ,0,\nu_1 )}$):

\bigbreak\noindent\textbf{Theorem$'$.} \emph{One has a Picard
groupoid equivalence
\begin{equation}\CL_{\text{dR}}^\Phi (X,T;\omega_X )\iso \CL_{\text{dR}} (X, T
)\times\! \mathop\Pi\limits_{b\in T_{\text{red}}}
\!\!\CL_{k_b},\quad \CE \mapsto ( \CE^{(1)}_{X\setminus T},\{
\CE_{\nu_{T_b} }\}). \label{ 1.6.3}\end{equation} }\bigbreak

{\it Variant.} More generally, we can fix a divisor $\Sigma d_b b$
supported on $T$ and take for $\nu_T$ a trivialization of $K(\Sigma
d_b b)$ on $T$.  The corresponding assertion is equivalent to the
above theorem by Remark (ii) in 1.2.

{\it Example.} If $K=\omega_X$ and $T=T_{\text{red}}$, then a
convenient choice is $d_b \equiv 1$, for $K(b)^\times_b $ is
canonically trivialized by $\nu_1$ as above. We  denote the fiber of
$\CE^{(1)}_{T_b} $  at $\nu_1$ by $\CE^{(1)}_b$.

\subsection{} {\it For the subsections 1.7--1.9  we assume that $T=\emptyset$, so $\fD = \Div(X)$.}

For a Picard stack $\CL_?$ and a commutative monoid space $\text{D}$
denote by $\Hom (\text{D}, \CL_? )$ the Picard groupoid of symmetric
monoidal morphisms $\text{D} \to \CL_?$ (we view $\text{D}$ as a
``discrete" symmetric monoidal stack). Thus an object  of  $\Hom
(\text{D}, \CL_? )$ is  $\CF\in\CL_? (\text{D})$ together with a
{\it multiplication structure} which is a rule that assigns to every
finite collection $\{ D_\alpha \}$ of $S$-points of $\text{D}$ a
{\it multiplication
 isomorphism} $\otimes \,\CF_{D_\alpha }\iso \CF_{\Sigma D_\alpha }$
(where $\Sigma$ is the operation in $\text{D}$); the isomorphisms
should be compatible with base change and satisfy an evident
transitivity property.  If $\text{D}^{\text{gr}}$ is the group
completion of $\text{D}$, then  $\Hom (\text{D}^{\text{gr}}, \CL_?
)\iso\Hom (\text{D}, \CL_? )$.

We are interested in $D$ equal to the monoid space of  effective
divisors $\Div^{\eff}(X) = \sqcup\, \Sym^n (X) \subset \Div (X)$;
one has $\Div^{\eff} (X)^{\text{gr}}=\Div (X)$.  A multiplication
structure on $\CF$ being restricted to disjoint divisors makes a
factorization structure. Pulling $\CF$ back to $\fD^\diamond$ is a
Picard functor \begin{equation}\Hom (\Div(X), \CL_?) \to \CL_?^\Phi
(X; K).  \label{ 1.7.1}\end{equation} Let $ \CL_{\text{dR}}^0
\subset \CL_{\text{dR}}$ be the Picard stack of degree 0 de Rham
lines.

\bigbreak\noindent\textbf{Proposition.} \emph{One has $\Hom (\Div(X),
\CL_{\text{dR}}^{0}) \iso \CL_{\text{dR}}^{0\Phi }(X; K)$.
}\bigbreak

{\it Proof.} (a) Let us show that for any
$\CE\in\CL_{\text{dR}}^{0\Phi }(X; K)$ the de Rham lines $\CE^{(\ell
)}$ on $\CK^{(\ell )}$ come from de Rham lines on $X$.

The claim is $X$-local, so we  trivialize $K$ by section a $\nu_0$
and pick a function $t$ on $X$ with $dt$ invertible. Then
$\CK^{(\ell )}$ is trivialized by section $\nu_0^{(\ell )}:= \nu_0
dt^{\otimes -\ell}$; let $z$ be the corresponding fiberwise
coordinate on $\CK^{(\ell )}$.
 Choose $X$-locally a trivialization $e^{(\ell )}$ of $\CE^{(\ell )}$;
let $\theta^{(\ell )}\in \omega_{\CK^{(\ell )}/X}$ be the
restriction of $\nabla (e^{(\ell )})/e^{(\ell )}$ to the fibers.
Then $\theta^{(\ell )}= \Sigma f^{(\ell )}_k (x) z^{k} d\log z$,
where $f^{(\ell )}_k (x)$,   $k\in\Bbb Z$, are functions on $X$; we
want to show that $f^{(\ell )}_k (x)=0$ for $k\neq 0$,\footnote{The
proof   uses only the factorization $\CO$-line structure on $\CE$.}
and $f^{(\ell )}_0 (x)\in\Bbb Z$.

Let $S\subset X\times X$ be a sufficiently small neighborhood of the
diagonal, $x_1, x_2$ be the coordinate functions on $S$ that
correspond to $t$, $D^{(\ell_1 ,\ell_2)}\in \Div (X)(S)$ be the
divisor $\ell_1 \Delta_1 + \ell_2 \Delta_2$. Then  $ (t-x_1
)^{-\ell_1 }(t-x_2 )^{-\ell_2}\nu_0$ is a trivialization of $K(
D^{(\ell_1 ,\ell_2)})$ near $|D^{(\ell_1 ,\ell_2 )}|$. Denote by
$\nu_0^{(\ell_1 ,\ell_2 )}$ the corresponding section of $K(
D^{(\ell_1,\ell_2)})^\times_{D^{(\ell_1 ,\ell_2 )},0}$ (see 1.1);
set $\CK^{(\ell_1 ,\ell_2 )}:= \Bbb G_m \nu_0^{(\ell_1 ,\ell_2
)}\subset   K( D^{(\ell_1,\ell_2)})^\times_{D^{(\ell_1 ,\ell_2
)},0}$, $\CE^{(\ell_1 ,\ell_2 )}:= \CE_{(D^{(\ell_1 ,\ell_2 )},0
)}|_{\CK^{(\ell_1 ,\ell_2 )}}$.

Outside the diagonal in $S$ one has an embedding $i^{(\ell_1 ,\ell_2
)}:\CK^{(\ell_1 ,\ell_2 )}\hra pr_1^* \CK^{(\ell_1 )}\times pr_2^*
\CK^{(\ell_2)}$ defined by the factorization; explicitly, it
identifies $z \nu_0^{(\ell_1 ,\ell_2 )}(x_1 ,x_2 )$ with $(z(x_1 -
x_2 )^{-\ell_2} \nu_0^{(\ell_1 )}(x_1 ), z (x_2 - x_1 )^{-\ell_1}
\nu_0^{(\ell_2 )}(x_2 ))$. Restricting to $\CK^{(\ell_1 ,\ell_2 )}$
the image of $e^{(\ell_1 )}\boxtimes e^{(\ell_2 )}$ by the
factorization isomorphism (1.5.1), we get a trivialization
$e^{(\ell_1 ,\ell_2 )}$ of $\CE^{(\ell_1 ,\ell_2 )}$ outside the
diagonal in $S$. Let $m(\ell_1 ,\ell_2 )$ be its order of pole at
the diagonal, so $(x_1 - x_2 )^{m(\ell_1 ,\ell_2 )} e^{(\ell_1
,\ell_2 )}$ is a trivialization of $\CE^{(\ell_1 ,\ell_2 )}$ on $S$.
Therefore the restriction $\theta^{(\ell_1 ,\ell_2 )}$ of $\nabla
(e^{(\ell_1 ,\ell_2 )})/e^{(\ell_1 ,\ell_2 )}$ to the fibers is a
regular relative form, which  equals  $i^{(\ell_1 ,\ell_2
)*}(pr_1^* \theta^{(\ell_1 )}+ pr_2^* \theta^{(\ell_2 )})=\Sigma
(f^{(\ell_1 )}_k (x_1 )(x_1 - x_2 )^{-k\ell_2} +f^{(\ell_2 )}_k (x_1
)(x_2 - x_1 )^{-k\ell_1})z^k d\log z$.

Since  $\theta^{(\ell,-\ell)}$ has no pole at the diagonal, the
above  formula implies that $f^{(\ell)}_k =0$ for $k\ell <0$.
Similarly, the formula for $\theta^{(\ell ,2\ell)}$ shows that
$f^{(\ell)}_k =0$ for $k\ell >0$. To see that $f^{( \ell )}_0
\in\Bbb Z$, notice that the above picture for $\ell_1 =\ell_2 =\ell$
is symmetric with respect to the transposition involution $\sigma$
of $X\times X$, hence descends to $S/\sigma \subset \Sym^2 X$. Thus
$m(\ell,\ell )$ is even. One has $\nabla (e^{( \ell)})/e^{( \ell)} =
f^{( \ell )}_0 d\log z + g^{( \ell)}(x) dx$ where $g^{( \ell)}(x) $
is a regular function. Then  $\nabla (e^{(\ell,\ell
)})/e^{(\ell,\ell )} +d\log   (x_1 - x_2 )^{m(\ell,\ell)} $ is a
regular 1-form. It equals $  (f^{(\ell)}_0 (x_1 ) + f^{(\ell)}_0
(x_2 ))(d\log z - d\log (x_1 - x_2 )) + m(\ell,\ell)d\log (x_1 - x_2
) + g^{( \ell)}(x_1 ) dx_1 + g^{( \ell)}(x_2 ) dx_2$. Therefore $
f^{(\ell)}_0 (x)= m(\ell,\ell)/2 \in\Bbb Z$, and we are done.

(b)  The next  properties of de Rham lines will be repeatedly used.
Let $\pi: K\to S$ be a smooth morphism of smooth schemes with dense
image.

\bigbreak\noindent\textbf{Lemma.}\emph{ (i) The functor $\pi^*
:\CL_{\text{dR}}(S)\to\CL_{\text{dR}}(K)$ is faithful. If  the
geometric fibers of $\pi$ are connected (say, $\pi$ is an open
embedding), then $\pi^* $ is  fully faithful. \newline (ii) If, in
addition,  $\pi$ is surjective, then a de Rham line $\CE$ on $K$
comes from $S$ if (and only if) this is true over a neighborhood $U$
of the generic point(s) of $S$. \hfill$\square$ }\bigbreak

(c) As was mentioned, $\Div (X)$ is the group completion of
$\Div^{\eff}(X)=\sqcup \Sym^n (X)$. So we have the projection
$\Div^{\eff}(X)\times\Div^{\eff}(X)\to\Div (X)$, $(D_1 ,D_2 )\mapsto
D_1 - D_2$, which identifies $\Div (X)$ with the quotient of
$\Div^{\eff}(X)\times\Div^{\eff}(X)$ with respect to the diagonal
action. Therefore a line $\CE$ on $\Div (X)$ is the same as a
collection of lines $\CE^{n_1 ,n_2}$ on
 $\Sym^{n_1 ,n_2}(X):= \Sym^{n_1}(X)\times  \Sym^{n_2}(X)$ together
with $\Div^{\eff}(X)$-equivariance structure, which is the datum of
identifications of their pull-backs by  $\Sym^{n_1
,n_2}(X)\leftarrow \Sym^{n_1 ,n_2}(X)\times\Sym^{n_3} (X)\to
\Sym^{n_1 +n_3 ,n_2 + n_3}(X)$, $(D_1 ,D_2 )\leftarrow (D_1 ,D_2
;D_3)\to (D_1 +D_3 ,D_2 +D_3)$ that satisfy a transitivity property.

Let us prove the proposition.  We need to show that any $\CE\in
\CL_{\text{dR}}^{0\Phi }(X; K)$, viewed as a mere  de Rham line on
$\fD^\diamond$, is the pull-back by  $\fD^\diamond \to \Div (X)$ of
a uniquely defined  line in $\CL_{\text{dR}}^0 (\Div (X))$, which we
denote by $\CE$ or $\CE_\Div$, and that  the factorization structure
on $\CE$ comes from a uniquely defined multiplication structure on
$\CE_\Div$.

We use the fact that for any $D\in\Div (X)(S)$, $S$ is smooth, the
projection $K(D)^\times_{D}:= K(D)^\times_{D,0}\to S$ satisfies the
conditions of (i), (ii) of the lemma.

To define $ \CE_\Div$ on $S= \Sym^{n_1 ,n_2}(X)$, we apply (ii) of
the lemma to $\CE$ on $K=K(D_1 -D_2 )^\times_{D_1 + D_2}$. Let $U$
be the complement to the diagonal divisor in $X^{n_1}\times
X^{n_2}$. Over $U$ our $\CE$ equals $(\CE^{(1)})^{\boxtimes
n_1}\boxtimes  (\CE^{(-1)})^{\boxtimes n_2}$, and we are done by
(a). The $\Div^{\eff}(X)$-equivariance structure on the datum of $
\CE^{n_1 ,n_2}_\Div$ is automatic by (i) of the lemma (applied to
$K(D_1 -D_2  )_{D_1 +D_2 +2D_3}$).

The factorization structure on $\CE$ yields one on $\CE_\Div$. Let
us show that it extends uniquely to a multiplication structure.  It
suffices to define the multiplication $\otimes
\CE_{D_\alpha}\to\CE_{\Sigma D_\alpha}$ over each $\Pi
\Sym^{n_{1\alpha} ,n_{2\alpha}}(X)$ in a way compatible with the
$\Div^{\eff}(X)$-equivariance structure. Our multiplication equals
factorization over the open dense subset where all $D_{i\alpha}$ are
disjoint, so we have it everywhere by (i) of the lemma. The
compatibility with $\Div^{\eff}(X)$-equivariance holds over the
similar open subset of $\Pi  (\Sym^{n_{1\alpha} ,n_{2\alpha}}(X)
\times\Sym^{n_{3\alpha}}(X))$, hence everywhere, and we are done.
\hfill$\square$

\bigbreak\noindent\textbf{Corollary.}\emph{ The functor
$\CL_{\text{dR}}^{0\Phi }(X;K)\to \CL^0_{\text{dR}}(X)$,    $\CE
\mapsto \CE^{(1)}$, is an equivalence. }\bigbreak

{\it Proof.} Its composition with the equivalence of the proposition
is a functor $ \Hom (\Div^{\eff} (X), \CL_{\text{dR}}^{0})\to
\CL_{\text{dR}}^{0}(X)$ which assigns to $\CF$ its restriction to
the component $X=\Sym^1 (X)$ of $\Div^{\eff} (X)$. This functor is
clearly invertible: its inverse assigns to $\CP \in
\CL_{\text{dR}}^{0}(X)$ the de Rham line $\Sym (\CP )$ on
$\Div^{\eff} (X)= \sqcup \Sym^n (X)$, $\Sym (\CP )_{\Sym^n (X)}:
=\Sym^n (\CP )$ equipped with an evident multiplication structure.
\hfill$\square$

\subsection{} An example of a de Rham factorization line $\CE$ with deg$(\CE^{(1)})=1$:

Suppose $X=\Bbb A^1$, $T=\emptyset$, $K=\CO_X$. We construct $\CE$
in the setting of Remark in 1.3. For a smooth $S$ and $D\in$ Div$
(S)$ the line $\CO_X (D)$ is naturally trivialized by a section
$\nu_D$, $\nu_{\Sigma n_i x_i} = \Pi (t-x_i )^{-n_i}$. Then $\nu_D$
trivializes the $O_{D,0}^\times$-torsor $K(D)^\times_{D,0}$, so the
canonical character $f \mapsto f(D)$ of $O_{D,0}^\times$ yields an
invertible function $\phi_D \in \CO^\times ( K(D)^\times_{D,0})$,
$\phi_D (\nu ):=(\nu/\nu_D )(D)$. Our $\CE_{(D,0)}$ comes from the
Kummer torsor for $\phi^{-1/2}_D$ placed in degree $\text{deg}(D)$,
i.e., it equals $\CO_{K(D)^\times_{D,0}}[\text{deg}(D)]$ as an
$\CO$-line, the connection is given by the 1-form $\frac{1}{2}d\log
\phi_D$.

Notice that if $D' \in$ Div$ (S)$ is another divisor such that
$|D|\cap |D'|=\emptyset$, then the invertible function $\nu_D$ on
$X_S \setminus |D|$ yields $(D,D'):= \nu_D (D')\in\CO^\times (S)$.
One has \begin{equation}(D,D')=
(-1)^{\text{deg}(D)\text{deg}(D')}(D',D). \label{
1.8.1}\end{equation}

Let us define the factorization structure on $\CE$.  Suppose $D=
\sqcup D_\alpha$, so $K(D)^\times_{D,0} = \Pi K(D_\alpha
)^\times_{D_\alpha ,0}$.  Any linear order on the set of indices
$\alpha$ yields an evident identification of the ``constant"
$\CO$-lines $\otimes \CE_{(D_\alpha ,0)}\iso \CE_{(D,0)}$. The
factorization isomorphism (1.5.1) is  its product with
$\mathop\Pi\limits_{\alpha <\alpha'} (D_\alpha ,D_{\alpha'})$. The
choice of order is irrelevant due to the ``super" commutativity
constraint and (1.8.1). Both transitivity property and horizontality
follow since $\Pi \phi_{D_\alpha}=\phi_D \mathop\Pi\limits_{\alpha
\neq\alpha'} (D_\alpha ,D_{\alpha'}) $.

\subsection{} {\it Proof of the theorem in 1.6 in case $T=\emptyset$.}
Let us check that for  $\CE \in \CL_{\text{dR}}^\Phi (X;K)$ one has
$\CE^{(1)}\in\CL_{\text{dR}}(X)^{(1)}$. The claim is $X$-local, so
we can assume that $K$ is trivialized and there is a function $t$ on
$X$ with $dt$ invertible, i.e., $t : X\to \Bbb A^1$ is \'etale. Let
$\CE' \in  \CL_{\text{dR}}^\Phi (X;K)$ be the pull-back of  the
factorization line from 1.8. Then $\CL_{\text{dR}}^\Phi (X;K)$ is
generated by  $\CL_{\text{dR}}^{0\Phi }(X;K)$ and $\CE'$. Our claim
holds for $\CE\in \CL_{\text{dR}}^{0\Phi }(X;K) $ by 1.7 and it is
evident for $\CE'$; we are done.

Let us show that  $\CL_{\text{dR}}^\Phi (X;K)\to
\CL_{\text{dR}}(X)^{(1)}  $,  $\CE \mapsto \CE^{(1)}$, is an
equivalence. Notice that the preimage of
$\CL^0_{\text{dR}}(X)\subset \CL_{\text{dR}}(X)$ equals
$\CL_{\text{dR}}^{0\Phi }(X;K)$, and, by the corollary in 1.7, $
\CL_{\text{dR}}^{0\Phi }(X;K)\iso \CL^0_{\text{dR}}(X)$. Since
$X$-locally there is $\CE$ with $\deg (\CE^{(1)})=1$ by 1.8 and
  $\deg: \pi_0 (\CL_{\text{dR}}(X))/\pi_0 (\CL^0_{\text{dR}}(X))\iso \Bbb Z$, we are done.  \hfill$\square$

\subsection{}
Suppose now  $T\neq \emptyset$. Pick $b\in T_{\text{red}}$, and
consider the $O^\times_{T_b}$-torsor $K^\times_{T_b}$ (see Example
(ii) in 1.1; we follow the notation of loc.cit.).

Let $(\omega)_b := \omega_X (\infty b)/\omega_X$ be the $k_b$-vector
space of polar parts  of rational 1-forms at $b$,  $(\omega)^{\le
n}_b$ be the subspace of polar parts of order $\le n$.  The  Lie
algebra of $O^\times_{T_b}$ equals $\CO (T_b )$. The space $\Omega^1
(K^\times_{T_b})^{\text{inv}}$ of translation invariant 1-forms on
$K^\times_{T_b}$ is its dual. The residue pairing $(\omega )^{\le
m_b }_b \times \CO(T_b )\to k_b$, $(\psi, f) \mapsto \Res_b (f\psi
)$,  identifies it with $(\omega)_b^{\le m_b }$. So one has
\begin{equation} \Omega^1( K^\times_{T_b})^{\text{inv}} \iso
(\omega)_b^{\le m_b } .\label{ 1.10.1}\end{equation}

Let $U$ be a smooth affine curve over $k_b$, $u\in U$ a closed
point; as above, we set $(\omega )_u := \omega_U (\infty
u)/\omega_U$. Let $ \xi :U^o := U\setminus \{ u\}\to K^\times_{T_b}$
be a $k_b$-morphism, which amounts to a trivialization $\nu^\xi$ of
$K$ on  $T_{b\, U^o}\subset X_{U^o}$. Denote by $(\xi )$ the
composition $(\omega)_b^{\le m_b }\iso \Omega^1( K^\times_{T_b}
)^{\text{inv}}\buildrel{\xi^*}\over\lra \omega (U^\circ )\to\omega
(U^\circ )/\omega (U)= (\omega)_u$.

\bigbreak\noindent\textbf{Lemma.}\emph{ (i) After a possible
localization of $U$ at $u$, one can find $D\in\Div (X)(U)$ and a
trivialization $\nu$ of $K(D)$ on a neighborhood $V\subset X_{U}$ of
$ |D|\cup T_{b U}$ such that $|D|\cap X_u$ is supported at $b$,
$|D|\cap T_{b\, U^o}=\emptyset$,  and $\nu|_{T_{b\, U^o}}=\nu^\xi$.
\newline (ii) Suppose $U$ is a neighborhood of $b$, i.e., we have an \'etale
$\pi : U\to X$, $\pi (u)=b$. Then one can find $(D,\nu )$ as in (i)
with $D$ equal to (the graph of) $\pi$ if and only if $-(\xi )$
equals $\pi^*$, i.e., the composition $ (\omega)_b^{\le m_b }\subset
(\omega)_b \buildrel{\pi^*}\over\hra (\omega)_u$.
 }\bigbreak

{\it Proof.} (i) Let us extend $\nu^\xi$ to a rational section $\nu$
of $K$ on an open subset $V$ of $X_U$, $V\supset T_{b\, U}$, which
is defined at $T_{b \, U^o}$.  Shrinking $U$ and $V$, one can find
$\nu$ with $D:=\text{div} (\nu)$  prime to $X_u$ (if $n$ is the
multiplicity of $X_u$ in $D$, then we replace $\nu$ by $f^{-n} \nu$,
where $f$ is any rational function which equals 1 on $T_{b\, U^o}$
and whose divisor contains $X_u$ with multiplicity 1).\footnote{To
find such $f$ (after possible shrinking of $U$), pick local
coordinate $t$ on $X$ at $b$, and $x$ on $U$ at $u$ (so $t(b)=0=x(u
)$, $dt (b)\neq 0\neq dx (u )$); set $f = x (x-t^{m_b} )^{-1}$.}
After further  localization of $U$ and shrinking of $V$, we get $D$
in $\Div (X)(U)$ and $|D|\cap X_u $ is supported at $b$; we are
done.

(ii) A map $\phi : U^\circ \to K^\times_{T_b}$ extends to $U$ if and
only if for every $\beta \in \Omega^1( K^\times_{T_b} )^{\text{inv}}
$ the form $\phi^* (\beta ) \in \omega (U^\circ )$ is regular at
$b$. Thus either of the properties of $\xi$ in the assertion of (ii)
determines $\nu^\xi$  uniquely up to multiplication by an invertible
function on $T_{b  U}$. It remains to present a trivialization $\nu$
of $K( \pi )$  such that  the corresponding $\xi$ satisfies $-(\xi
)=\pi^*$.

Shrinking $X$, we  trivialize $K$ and  pick a  function $t$ with
$dt$ invertible; set $x:= \pi^* (t) \in\CO (U)$. Our $\nu$ is
$(t-x)^{-1}$. The differential of the corresponding $\xi$ is the
Lie$(O^\times_{T_b} )=O(T_b )$-valued 1-form  $\nu^{-1} d_x \nu =-
(1+ t/x + (t/x)^2 +\ldots )dx/x$. So if $\beta\in \Omega^1
(K^\times_{T_b})^{\text{inv}}$ is identified with $\psi (t) \in
(\omega)_b$ by (1.10.1), then $\xi^* (\beta )=- (\Res_b ( 1+ t/x +
(t/x)^2 +\ldots )\psi (t) )dx/x =-\psi (x)$, q.e.d. \hfill$\square$

\subsection{} A de Rham line $\CF $ on $K^\times_{T_b}$ is said to be
{\it translation invariant} if the de Rham line $\cdot^* \CF \otimes
pr_2^* \CF^{\otimes -1}$ lies in $pr_1^*
\CL_{\text{dR}}(O^\times_{T_b} ) \subset
\CL_{\text{dR}}(O^\times_{T_b} \times K^\times_{T_b})$; here $\cdot
: O^\times_{T_b} \times K^\times_{T_b}\to K^\times_{T_b}$ is the
action map, $pr_i$ are the projections to the factors. Such $\CF$'s
form a Picard subgroupoid
$\CL_{\text{dR}}^{\text{inv}}(K^\times_{T_b})$ of
$\CL_{\text{dR}}(K^\times_{T_b})$.

\bigbreak\noindent\textbf{Lemma.}\emph{ (i) We trivialize
$K^\times_{T_b}$, i.e., identify it with  $O^\times_{T_b}$. The
translation invariance of $\CF$ is equivalent to the next
properties: \newline (a) The de Rham line $ pr_1^* \CF \otimes
pr_2^* \CF \otimes \cdot^* \CF^{\otimes -1}$ is constant; \newline
(b) For any smooth curve $U$, a point $u\in U$, and  two maps $\xi_1
,\xi_2  : U^o := U\setminus \{ u\} \to K^\times_T $, the de Rham
line $ \xi^*_1 \CF \otimes \xi^*_2 \CF \otimes (\xi_1\xi_2 )^*
\CF^{\otimes -1} $ on $U^o$ extends to $U$.
\newline
(c) For some (or every) invertible section $e_\CF$ of $\CF$ on
$K^\times_{T_b}$ one has $\nabla (e_\CF )/e_\CF \in\Omega^1
(K^\times_{T_b})^{\text{inv}}$. (ii) There is a natural isomorphism
$\pi_0 (\CL_{\text{dR}}^{\text{inv}}(K^\times_{T_b}))\iso   \Bbb Z
\times ((\omega)_b^{\le m_b} /\Bbb Z )$ where $\Bbb Z \subset
(\omega)_b^{\le m_b}$ are polar parts of 1-forms with simple pole
and integral residue. }\bigbreak

{\it Proof.} (i) (a) is evidently equivalent to invariance of $\CF$.
Since $K^\times_{T_b}$ is a rational variety, (a) amounts to the
fact that $\cdot^* \CF \otimes pr_1^* \CF^{\otimes -1} \otimes
pr_2^* \CF^{\otimes -1}$ extends to a compactification of
$K^\times_{T_b}$. This can be tested on curves, which is  (b).
Finally (c) is equivalent to the translation invariance since
 every invertible function on $O^\times_{T_b}$ is the product of a character by a constant, and
every line bundle on $K^\times_{T_b}$ is trivial.

(ii) One assigns to $\CF$ the pair $(n,\psi )$ where $n=\deg (\CF )$
and $\psi$ is the class of the image of $\nabla (e_\CF )/e_\CF$ by
(1.10.1). \hfill$\square$

We say that $\CG\in \CL_{\text{dR}} (X\setminus T)^{(1)}$ is {\it
compatible} with $\CF\in
\CL^{\text{inv}}_{\text{dR}}(K^\times_{T_b})$ if  for some
neighborhood $U$ of $b$, a trivialization $\nu$  of $\CK^{(1)}$ on
$U$, and a map $\xi : U^o := U\setminus \{b\} \to K^\times_{T_b}$ as
in (ii) of the lemma in 1.10, the de Rham line $\nu^* \CG
\otimes\xi^* \CF$ on $U^o$ extends to $U$ (the validity of this does
not depend on the choice of $U$, $\nu$, and $\xi$). By loc.cit.,
compatibility is equivalent to the next condition: Pick $U$, $\nu$,
and $e_\CF$ as above; let $e_\CG$ be a non-zero rational section of
$\nu^* \CG$. Then the image of $\nabla (e_\CF )/e_\CF$    by
(1.10.1) in $(\omega)_b^{\le m_b} /\Bbb Z  \subset  (\omega)_b /\Bbb
Z $ equals the class of  $\nabla( e_\CG )/e_\CG$.

Let $\CL^\natural_{\text{dR}} (X,T;K)$ be the Picard subgroupoid of
$\CL_{\text{dR}} (X\setminus T)^{(1)}\times \mathop\Pi\limits_{b\in
T_{\text{red}}} \CL^{\text{inv}}_{\text{dR}}(K^\times_{T_b})$ formed
by those collections $(\CG,\{  \CF_{T_b} \})$ that $\CG$ is
compatible with every $\CF_{T_b}$. Then $\CG$ lies automatically in
$\CL_{\text{dR}}(X,T)^{(1)}$. By (ii) of the lemma, the functor
$\CL^\natural_{\text{dR}} (X,T;K)\to
\CL_{\text{dR}}(X,T)^{(1)}\times \Pi \CL_{k_b}$, $(\CG ,\{ \CF_{T_b}
\})\mapsto (\CG, \{\CF_{\nu_{T_b}}\})$, where $\CF_{\nu_{T_b}}$ is
the fiber of $\CF_{T_b}$ at $\nu_{T_b}$ from 1.6,
  is an equivalence of categories. Thus the theorem in 1.6 follows from the next one:

\bigbreak\noindent\textbf{Theorem.}\emph{ For every $\CE \in
\CL^\Phi_{\text{dR}} (X,T;K)$ one has $(\CE^{(1)},\{ \CE_{T_b}\})\in
\CL^\natural_{\text{dR}} (X,T;K)$, and the functor \begin{equation}
\CL^\Phi_{\text{dR}} (X,T;K)\to \CL^\natural_{\text{dR}}
(X,T;K),\quad \CE \mapsto (\CE^{(1)},\{ \CE_{T_b}\}), \label{
1.11.1}\end{equation} is an equivalence of the Picard
groupoids.}\bigbreak

{\it Proof.}   The assertion is $X$-local, and we have proved it for
$T=\emptyset$. So we can assume that $T_{\text{red}}$ is a single
$k$-point $b$. Thus $T_b =T$ and $\fD$ is the disjoint sum of
$\fD_{c=0} $ equal to $\Div (X\setminus T)$ and $\fD_{c=1} $ equal
to $\Div (X)$. If needed, we can assume that
  $K$ is trivialized and there is an \'etale map $t: X\to\Bbb A^1$.

(a) Let us show that $(\CE^{(1)}, \CE_{T})\in
\CL^\natural_{\text{dR}} (X,T;K)$. Notice that $
\CL^\Phi_{\text{dR}} (X,T;K)$ is generated by $
\CL^{0\Phi}_{\text{dR}} (X,T;K)$, the image of (1.5.3), and  the
pull-back by $t$ of the factorization line on $\Bbb A^1$ from 1.8.
Since the assertion is evident for factorization lines of the latter
two types, it suffices to consider the case  of $\CE\in
\CL^{0\Phi}_{\text{dR}} (X,T;K)$.

We know that $\CE^{(1)} $ comes from a de Rham line on $X\setminus
T$ (see 1.7). Let us check that $\CE_{T}$ is translation invariant
using the criterion of (i)(b) in the lemma. For $U$, $u$, $\xi_i$ as
in loc.cit., let us choose $D_i$, $\nu_i$  as in (i) of the lemma in
1.10; then $D_3 :=D_1 +D_2$, $\nu_3 :=\nu_1 \nu_2$ serves  $\xi_3 :=
\xi_1 \xi_2$. The lines $\CE_{(D_i, 1, \nu_i )}$    on $U$ are equal
to $\CE_{D_i } \otimes \xi^*_i \CE_T $ on $U^o$ by factorization;
here $\CE_{D_i}:= \CE_{(D_i ,0)}$. Since $\CE_{D_3}=\CE_{D_1}\otimes
\CE_{D_2}$ by 1.7, the de Rham line $\xi^*_1 \CE_T \otimes \xi^*_2
\CE_T\otimes (\xi_1\xi_2 )^* \CE_T^{\otimes -1} $ on $U^o$ extends
to $U$, q.e.d.

It remains to check that $\CE^{(1)}$ is compatible with $\CE_T$. Let
$\nu$ is a trivialization of $K(\Delta )$ on an open $V\subset
X\times X$ that contains $(b,b)$, $U:= V\cap (\{b\} \times X)$, $\xi
: U^o \to K^\times_{T}$ the map defined by the restriction of $\nu$
to $T\times U^o$. The de Rham line $\CE_{(\Delta,1,\nu )}$ on $U$
equals $\CE^{(1)}\otimes \xi^* \CE_T$ on $U^o$ by factorization.
Since the compatibility means that the latter line extends to $U$,
we are done.

(b) Consider the projection $\pi : \fD^\diamond_{c=1} \to K^\times_T
\times \Div (X)$, $(D,1,\nu_P) \mapsto ( \nu_P |_{T_S} ,D)$. Let us
show for any $\CE \in  \CL^{0\Phi}_{\text{dR}} (X,T;K)$ its
restriction $\CE_1$ to $\fD^\diamond_{c=1}$ comes from a uniquely
defined de Rham line on $ K^\times_T \times \Div (X)$ which we
denote by $\CE_1$ or $\CE_{\Div 1}$.

We use the fact that for every $D\in\Div (X)(S)$, $S$ is smooth, the
projection $K(D)^\times_{D,1}\to   K^\times_T \times S$ satisfies
the conditions of (i), (ii) of the lemma in  1.7.

As in part (c) of the proof in 1.7, we need to define $\CE_{\Div 1}$
on every  $ K^\times_T \times \Sym^{n_1 ,n_2}(X)$ and provide the
$\Div^{\eff}(X)$-equivariance structure. Consider our $\CE$ on
$K(D_1 - D_2 )^\times_{D_1 + D_2 ,1}$. Over $ K^\times_T \times
\Sym^{n_1 ,n_2}(X\setminus T)$ it equals $\CE_T \boxtimes \CE_{(D_1
- D_2 ,0)}$ by factorization, hence it descends to $ K^\times_T
\times \Sym^{n_1 ,n_2}(X\setminus T)$ by 1.7. By (ii) of the lemma
in 1.7, we have $\CE_{\Div 1}$ over the whole $ K^\times_T \times
\Sym^{n_1 ,n_2}(X)$. The $\Div^{\eff}(X)$-equivariance is automatic
by (i) of the lemma (applied to $K(D_1 - D_2 )^\times_{D_1 + D_2 +
2D_3}$).

(c) Our functor sends $ \CL^{0\Phi}_{\text{dR}} (X,T;K)$ to the
Picard subgroupoid $\CL^{0\natural}_{\text{dR}} (X,T;K)$  of
$\CL^\natural_{\text{dR}} (X,T;K)$ formed by all $(\CG ,\{ \CF_b
\})$ with $\deg (\CG )=\deg (\CF_b )= 0$.  Let us prove that  $
\CL^{0\Phi}_{\text{dR}} (X,T;K)\to \CL^{0\natural}_{\text{dR}}
(X,T;K)$ is an equivalence.

We need to show that every $(\CE^{(1)}, \CE_T )\in
\CL^{0\natural}_{\text{dR}} (X,T;K)$ comes from a uniquely defined
$\CE\in\CL^{0\Phi}_{\text{dR}} (X,T;K)$. By the corollary in 1.7,
$\CE^{(1)}$ defines $\CE_0 := \CE|_{\fD_{c=0}}$, which we can view,
by 1.7, as a de Rham line with multiplication structure on $\Div
(X\setminus T)$. As in (b),   $\CE_1 :=\CE|_{\fD^\diamond_{c=1}}$
comes from $K^\times_T \times \Div (X)$. By factorizaion, its
restriction to $K^\times_T \times \Div (X\setminus T)$ equals $\CE_T
\boxtimes \CE_0$. It remains to show that $\CE_T \boxtimes \CE_0$
extends in a unique way to a de Rham line $\CE_1$ on $K^\times_T
\times \Div (X)$.

As in (c) of the proof in 1.7, we should define $\CE_1$ on every
$K^\times_T \times \Sym^{n_1 ,n_2}(X)$ and provide the $\Div
(X)^{\eff}(X)$-equivariance structure. Our $\CE_1$  is defined on an
open dense subset $U$ of triples $(\xi, D_1, D_2 )$, $D_i \in
\Sym^{n_i}(X\setminus T )$. Let $U' \supset U$ be the open subset of
those $(\xi, D_1, D_2 )$ that $D_1 + D_2$  contains $b$ with
multiplicity at most 1. Then $\CE$ extends to $U'$ due to
compatibility of $\CE^{(1)}$ and $\CE_T$. Since the complement to $
U'$ has codimension $\ge 2$, $\CE$ extends to $K^\times_T \times
\Sym^{n_1 ,n_2}(X)$, and we are done.

As in loc.cit., the  $\Div (X)^{\eff}(X)$-equivariance is
identification of the pull-backs of our line by $K^\times_T \times
\Sym^{n_1 ,n_2}(X)\leftarrow K^\times_T \times \Sym^{n_1
,n_2}(X)\times\Sym^{n_3}(X)\to K^\times_T \times \Sym^{n_1 +n_3 ,n_2
+ n_3}(X)$. The two de Rham lines coincide on the dense open subset
$K^\times_T \times \Sym^{n_1 ,n_2}(X\setminus T)\times
\Sym^{n_3}(X\setminus T)$, so they are canonically identified
everywhere, and we are done. The factorization structure on $\CE$ is
evident.

(d) By (c), the theorem is reduced  to the claim that our functor
yields an isomorphism between the quotients $$\pi_0 ( \CL^{\Phi}_{\text{dR}}
(X,T;K))/\pi_0 (\CL^{0\Phi}_{\text{dR}} (X,T;K))\iso \pi_0
(\CL^{\natural}_{\text{dR}} (X,T;K))/\pi_0
(\CL^{0\natural}_{\text{dR}} (X,T;K)).$$ The degree map identifies
the right group with $\Bbb Z\times\Bbb Z$. Our map is evidently
injective; looking at  the image of (1.5.3) and the pull-back by $t$
of the factorization line on $\Bbb A^1$ from 1.8, we see that it is
surjective, q.e.d.  \hfill$\square$

\subsection{}
{\it A complement.} A connection on a trivialized line bundle
amounts to a 1-form; multiplying the trivialization by $f$, we add
to the form $d\log f$. Here is a similar fact in the factorization
story.

Consider the  group  $ \pi_1 (\CL_{\CO}^\Phi (X,T;K))$ of invertible
functions on $\fD^\diamond$ that satisfy factorization property. One
has evident embeddings $\CO^\times (T_{\text{red}})\hra \pi_1
(\CL_{\text{dR}}^\Phi (X,T;K))$ $\hra  \pi_1 (\CL_{\CO}^\Phi
(X,T;K))$ (see Remarks (iii), (iv) in 1.5). Let
$\CL_{\text{dR}}^\Phi (X,T;K)^{\CO\text{-triv}}$ be the kernel of
the Picard functor $\CL_{\text{dR}}^\Phi (X,T;K) \to \CL_\CO^\Phi
(X,T;K)$. This is a mere abelian group (since the functor is
faithful); its elements are pairs $(\CE ,e)$ where $\CE$ is a
factorization de Rham line, $e$ is a trivialization of $\CE$ as a
factorization $\CO$-line.  Let  $\omega (X,T)$ be the space  of
1-forms on $X\setminus T$ whose order of pole at any $b\in T$ is
less or equal to the multiplicity of $T$ at $b$.

\bigbreak\noindent\textbf{Proposition.}\emph{ There is a natural
commutative diagram
\begin{equation}
\begin{matrix}
  \pi_1 (\CL_{\CO}^\Phi (X,T;K))/ \CO^\times (T_{\text{red}})  & \iso
 &\CO^\times (X\setminus T)
 \\
 \downarrow\,\, \,\,&& \downarrow    \\
 \CL_{\text{dR}}^\Phi (X,T;K)^{\CO\text{-triv}}
& \iso
 &\omega (X,T).
\end{matrix} \label{ 1.12.1}
\end{equation}
}\bigbreak

{\it Proof.} (a)  The connection on $\CE^{(\ell)}$ along the fibers
of $\CK^{(\ell )}/X\setminus T$ is determined solely by the
$\CO$-line structure. So the action of any $h \in \pi_1
(\CL_{\CO}^\Phi (X,T))$  on $\CE^{(\ell )}$ is fiberwise horizontal,
i.e., it is multiplication by a function $h^{(\ell)}\in\CO^\times
(X\setminus T)$. The top horizontal arrow  is $h\mapsto h^{(1)}$.

For the same reason, for $(\CE ,e)\in  \CL_{\text{dR}}^\Phi
(X,T;K)^{\CO\text{-triv}}$ the trivializations $e^{(\ell)}$ of
$\CE^{(\ell)}$ are fiberwise horizontal, i.e., $\nabla
(e^{(\ell)})/e^{(\ell)}\in \omega (X\setminus T)$. By the theorem in
1.6,  $\nabla (e^{(1)})/e^{(1)}\in \omega (X,T)$.  The bottom
horizontal arrow is $(\CE ,e)\mapsto \nabla (e^{(1)})/e^{(1)} $.

The  map $ \pi_1 (\CL_{\CO}^\Phi (X,T;K)) \to \CL_{\text{dR}}^\Phi
(X,T;K)^{\CO\text{-triv}}$, $f\mapsto (\CO_{\fD^\diamond},f1)$, with
kernel  $\pi_1 (\CL_{\text{dR}}^\Phi (X,T;K))$ yields the left
vertical arrow. The right one is the $d\log$ map.

The diagram is evidently commutative. It remains to check that its
horizontal arrows are isomorphisms.

(b) For every $h\in \pi_1 (\CL_{\CO}^\Phi (X,T))$ its restriction to
$\fD^\diamond_{c=0}$ comes from a multiplicative function $h_0$ on
$\Div (X\setminus T)$. Similarly, if $T_{\text{red}}$ is a single
$k$-point $b$, then the restriction of $h$ to $\fD^\diamond_{c=1}$
comes from a function $h_1$ on $K^\times_T \times \Div (X)$ such
that for $\xi \in K^\times_T $, $D\in \Div (X\setminus T)$ one has
$h_1 (\xi ,D)= h_1 (\xi ) h_0 (D)$. This follows by a simple
modification of the argument  from, respectively, part (c) of the
proof in 1.7 and part (b) of the proof in 1.11. The details are left
to the reader.

  Let us show that the map $\pi_1 (\CL_{\CO}^\Phi (X,T))/\CO^\times (T_{\text{red}})
\to \CO^\times (X\setminus T)$ is injective.  Suppose we have $h$
such that $h^{(1)}=1$. Since the group space $\Div (X\setminus T)$
is generated by effective divisors of degree 1, one has $h_0 =1$. It
remains to check that $h$ is locally constant on other components of
$\fD^\diamond$. The assertion is $X$-local, so we can assume that
$T_{\text{red}}$ is a single $k$-point $b$, and we look at
$\fD^\diamond_{c=1}$. By above, it suffices to check that the
restriction $h_T$ of $h_1$ to $K^\times_T$ is constant. We use (i)
of the lemma 1.10; we follow the notation of loc.cit. For $\xi :U^o
\to  K^\times_{T}$, consider $h_{(D,1,\nu )}\in \CO^\times (U)$; by
factorization, its restriction to $U^o$ equals $ \xi^* h_{T}
h_{(D,0,\nu )} = \xi^* h_{T}$. Since $\xi^* h_{T}$ is regular at $u$
for every $\xi$, $h_{T}$ is  constant, q.e.d.

A  similar argument shows that the bottom horizontal arrow in
(1.12.1) is injective. The details are left to the reader.

(c) Let us construct a section $ \CO^\times (X\setminus T)\to \pi_1
(\CL_{\CO}^\Phi (X,T))$, $f\mapsto \tilde{f}$, of the map $h\mapsto
h^{(1)}$. Fix a trivialization $\nu_0$ of $K$ on an open subset
$V_0$ of $X$ that contains $T$. For $(D,c,\nu_P )\in\fD^\diamond
(S)$ let us define $ \tilde{f}_{(D,c,\nu_P )}\in\CO^\times (S)$.
Pick $\nu$, $V$ corresponding to  $ (D, c,\nu_P )$ as in Remark (ii)
in 1.1;  we can assume that $V\cap T_S = T^c_S$.  Localizing $S$, we
can decompose $D$ in a disjoint sum of $D'$ and $D''$ such that $D'
\subset V\setminus T_S$ and $D''\subset V_{0 S}$. Set
\begin{equation}\tilde{f}_{(D,c,\nu_P )} := f(D')\{ f,\nu_0 /\nu \}_{|D''|\cup
T^c_S}. \label{ 1.12.2}\end{equation} Here $\{ f,\nu_0 /\nu
\}_{|D''|\cup T^c_S}\in\CO^\times (S)$ is the Contou-Carr\`ere
symbol at $|D''|\cup T^c_S$ (see \cite{CC} or \cite{BBE} 3.3). One
readily checks  that (1.12.2) does not depend on the auxiliary
choices of $\nu$ and the decomposition $D=D'+D''$; its compatibility
with the factorization is evident. So $\tilde{f}\in \pi_1
(\CL_{\CO}^\Phi (X,T))$; clearly $\tilde{f}^{(1)}=f$.

{\it Remark.} If $\nu'_0$ is another trivialization of $K$ near $T$,
$f\mapsto \tilde{f}'$ the corresponding section, then
$\tilde{f}/\tilde{f}'$ is an element of $ \CO^\times
(T_{\text{red}})\subset \pi_1 (\CL_{\CO}^\Phi (X,T))  $ whose value
at $b\in T_{\text{red}}$ equals   $(\nu'_0/\nu_0 )(b)^{n_b}$, where
div$(f)=\Sigma n_b b$.

(d)  To finish the proof, let us construct explicitly a section of
the bottom horizontal arrow in (1.12.1). For $\phi \in \omega
(X,T)$, we construct the corresponding $\CE^\phi =(\CE^\phi ,e )\in
\CL_{\text{dR}}^\Phi (X,T;K)^{\CO\text{-triv}}$ using Remark (i) in
1.5. Since $\CE^\phi$ is trivialized as an $\CO$-line,
$\alpha^\varepsilon$ of (1.5.2) is multiplication by a function
$\alpha^\phi =\alpha^\phi_{ (D',c,\nu'_P )/ (D,c,\nu_P
)}\in\CO^\times (S)$. To determine it, we extend $\nu_P$, $\nu'_P$
to $\nu$, $\nu'$ as in Remark (ii) in 1.1 such that $\nu$ equals
$\nu'$ on $V_{\text{red}}$. Then $\nu/\nu' \in \CO^\times
(V\setminus P)$  equals 1 on  $V_{\text{red}}$, so we have a
function $\log (\nu/\nu' )\in \CO (V\setminus P)$ that vanishes on
$V_{\text{red}}$. The residue $\Res_{P/S} (\log (\nu/\nu' )\phi )\in
\CO (S)$ vanishes on $S_{\text{red}}$, and we set
\begin{equation}\alpha^\phi := \exp (\Res_{P/S} (\log (\nu/\nu'
)\phi )). \label{ 1.12.3}\end{equation} Our $\alpha^\phi$ does not
depend on the auxiliary choice of $\nu$ and $\nu'$: Indeed, $\nu$,
$\nu'$ can be changed  to $f\nu$,  $f'\nu'$ with $f,f'\in\CO^\times
(V)$ that coincide on $V_{\text{red}}$ and equal 1 on $ T^c_S$ (see
Remark (ii) in 1.1); then $\log (f/f')$ is a regular function on $V$
that vanishes on $T^c_S$, so $\Res_{P/S} (\log (f/f' )\phi )=0$, and
we are done. The transitivity of $\alpha^\phi$ and  compatibility
with base change and factorization are evident; we have defined
$\CE^\phi$.

{\it Remark.} Suppose we have $(D,c,\nu_P )\in\fD^\diamond (S)$
where $S$ is smooth. The de Rham structure on $\CE^\phi_{(D,c,\nu_P
)}$ amounts to a flat connection $\nabla^\phi$ on our line bundle,
which is the same as a closed 1-form $\theta^\phi = \nabla^\phi (e
)/e$ on $S$. Choose $\nu$ as in Remark (ii) in 1.1; then (1.12.3)
implies that \begin{equation}\theta^\phi = \Res_{P/S}( (d_S
(\nu)/\nu )\otimes \phi ). \label{ 1.12.4}\end{equation} Here $d_S$
means derivation along the fibers of the projection $V\subset X_S
\to X$, so $d_S (\nu )$ is a section of $\Omega^1_S \boxtimes K$
over $V\setminus P$, and $d_S (\nu )/\nu$ is a section of the
pull-back of $\Omega^1_S$ to $V\setminus P$. Of course, due to
  the lemma in 1.3, one can use (1.12.4) as an alternative definition of $\CE^\phi$.

{\it Example.} Consider the $\CO$-trivialized de Rham line
$(\CE^\phi_{T_b}, e)$ on the $O^\times_{T_b}$-torsor
$K_{T_b}^\times$ (see 1.6).  Its 1-form $\theta^\phi =\nabla (e)/e$
is translation invariant and corresponds to the functional $f
\mapsto \Res_b (f\phi )$ on the Lie algebra $\CO (T_b )$ of
$O^\times_{T_b}$ (cf.~1.11).

It remains to check that the bottom horizontal arrow in (1.12.1)
sends
 $(\CE^\phi ,e)$ to $\phi$, i.e., that $\nabla (e^{(1)})/e^{(1)}=\phi$.
Pick a local trivialization of $K$ and a local  function $t$ on
$X\setminus T$ with non-vanishing $dt$; let $x$ be the corresponding
local function on $S=X\setminus T$. Then $\nu = (t-x )^{-1}$ is a
trivialization of $K_S (\Delta )$ near the diagonal, so we have the
de Rham line $\CE^\phi_{(\Delta, 0,\nu )}$ on $S$. Then $d_S (\nu
)/\nu = (t-x)^{-1}dx$, so, by (1.12.4), one has $\nabla
(e^{(1)})/e^{(1)}
 =\theta^\phi = (\Res_{t=x}((t-x)^{-1} \phi ))dx=\phi$, q.e.d.  \hfill$\square$

\bigbreak\noindent\textbf{Corollary.}\emph{ For  $\CE, \CE' \in
\CL_{\text{dR}}^\Phi (X,T;K)$,  a morphism $ \CE\to \CE'$  in
$\CL_{\CO}^\Phi (X,T;K)$ is horizontal, i.e., is a morphism in
$\CL_{\text{dR}}^\Phi (X,T;K)$, if (and only if) the corresponding
morphism $\CE^{(1)} \to \CE^{\prime (1)}$  of  $\CO$-lines on
$\CK^{(1)}$ is horizontal.  \hfill$\square$ }\bigbreak

{\it Remark.}  If $K=\omega_X$, then in the corollary one can
replace $\phi^{(1)}$  by the  morphism $\phi^{(1)}_{X\setminus T}:
\CE^{(1)}_{X\setminus T}\to \CE^{\prime (1)}_{X\setminus T}$ of
$\CO$-lines on $X\setminus T$ (see 1.6).

\subsection{} The next lemma will be used in 2.12. Assume that
$X\setminus T$ is affine and $K=\omega_X$. The Lie algebra $\Theta
(X\setminus T)$ of vector fields on  $X\setminus T$ acts naturally
on $\fD^\diamond (X,\hat{T};\omega  ):= \limleft \fD^\diamond
(X,nT;\omega  )$. Therefore we have the notion of $\Theta
(X\setminus T)$-action on any $\CO$-line $\CE$ on $\fD^\diamond
(X,\hat{T};\omega  )$. If $\CE$ carries a factorization structure,
then one can ask our action to be compatible with it. Ditto for  a
de Rham structure.

Suppose that $\CE$ is a de Rham line. The flat connection yields
then a $\Theta (X\setminus T)$-action on $\CE$, which we refer to as
the {\it standard} action.  It is evidently compatible
 with the de Rham structure.

\bigbreak\noindent\textbf{Lemma.}\emph{ Any $\Theta (X\setminus
T)$-action $\tau$ on $\CE$  compatible with the de Rham structure
equals the standard one. }\bigbreak

{\it Proof.} Let 
$\tau_0$ be the standard action. Then $\theta \mapsto \tau (\theta
)-\tau_0 (\theta )$ is a Lie algebra homomorphism from $\Theta
(X\setminus T)$ to the Lie algebra of de Rham line endomorphisms of
$\CE$. The latter Lie algebra is commutative; the former one is
perfect. Thus our homomorphism is 0, i.e., $\tau =\tau_0$.
\hfill$\square$

\subsection{}
The whole story makes sense in the relative setting. The input is a
smooth (not necessary proper) $Q$-family of curves $q: X \to Q$
(where $Q$ is a scheme), a relative divisor $T \subset X$  such that
$T_{\text{red}}$ is finite and \'etale over $Q_{\text{red}}$, and a
line bundle $K$ on $X$. It yields a space $\fD^\diamond
=\fD^\diamond (X /Q,T ;K)$ over $Q$. If $\CL_?$ is any sheaf of
Picard groupoids on the category $\CS ch_{/Q}$ of $Q$-schemes
(equipped with the \'etale topology), then we have the Picard
groupoid $\CL_? (\fD )$ and $\CL_? (\fD^\diamond )$ defined as in
1.2, and the Picard groupoid of factorizarion objects $\CL^\Phi_? (X
/Q,T ;K)$ as in 1.5. The $\CL_?$'s we need are  $\CL_\CO$,
$\CL_{\text{dR}/Q}$, and $\CL_{\text{dR}}$, where $\CL_\CO$,
$\CL_{\text{dR}}$ are as in 1.2, and $\CL_{\text{dR}/Q}(S)$ is
formed by $\CO$-lines equipped with an action of the universal
relative formal groupoid on $S/Q$ (if $S/Q$ is smooth, then this is
the same as a flat relative connection, cf.~1.1). All the results
above immediately generalize to the relative setting. Thus, as in
1.4, for proper $q$ every $\CE$ in $\CL_\CO (\fD^\diamond )$ or in
$\CL_{\text{dR}/Q} (\fD^\diamond )$ yields an $\CO$-line $\CE (X
/Q)$ on $Q$; if $\CE$ lies in $\CL_{\text{dR}} (\fD^\diamond )$,
then  $\CE (X /Q)\in \CL_{\text{dR}} (Q)$. The theorem in 1.6
remains valid both in the setting of $\CL_{\text{dR}/Q}$ and
$\CL_{\text{dR}}$, etc.

\subsection{} Suppose $k=\Bbb C$, and $X$ is any complex smooth curve.
All the above definitions and results render  into the analytic
setting without problems. In fact, the story simplifies since
$\CL^\Phi_{\text{dR}}(X,T_{\text{red}})\iso \CL^\Phi_{\text{dR}}(X,T
)$.\footnote{Since for any $\Bbb A^1$-torsor $K/S$ the pull-back
functor $\CL_{\text{dR}}(S)\to\CL_{\text{dR}}(K)$ is an equivalence.
} It is equivalent to the Betti version of the story with de Rham
lines replaced by local systems of  $\Bbb C$-lines. And the Betti
version works if we replace this $\Bbb C$  by any commutative ring
$R$ of coefficients.

{\it Remark.} The fact that $\CL^\Phi_{\text{dR}}(X,T; K)=
\CL^\Phi_{\text{dR}}(X,T_{\text{red}};K)$ permits to consider in
case $K=\omega_X$ a {\it canonical} equivalence (1.6.3) as in
Example in 1.6.

Every de Rham factorization line $\CE$  in the analytic setting
carries a canonical automorphism $\mu$ which acts on $\CE_{(D,c,\nu
)}$ as multiplication by the (counterclockwise) monodromy of the de
Rham line $\CE_{(D,c,z\nu )}$, $z\in\Bbb C^\times$,  around $z=0$.
Notice that $\mu$ is multiplicative, i.e., we have a homomorphism
$\mu : \pi_0 (\CL^\Phi_{\text{dR}}(X,T ))\to \pi_1
(\CL^\Phi_{\text{dR}}(X,T ))$. Same is true for the Betti
factorization lines.

\section{The de Rham $\varepsilon$-lines: algebraic theory}

This section recasts the story of \cite{Del} and \cite{BBE} in the
factorization line format.

\subsection{}
We follow the setting and notation of 1.1, so $X$ is a smooth curve
over $k$, $T\subset X$ is a  finite subscheme.  {\it From now on we
assume  that $K$ from 1.1 equals $\omega= \omega_X$, so
$\fD^\diamond =\fD^\diamond (X,T;\omega)$.}

Let $M$ be a (left) holonomic $\CD$-module on $(X,T)$, i.e., a
holonomic module on $X$ which is smooth on  $X\setminus T$. We say
that $T$ is {\it compatible} with $M$ if $\det M_{X\setminus
T}\in\CL_{\text{dR}} (X,T)$ (see 1.6).

\bigbreak\noindent\textbf{Theorem-construction.}\emph{ $M$ defines
naturally a de Rham factorization line  $\CE_{\text{dR}}(M)$ $\in
\CL^\Phi_{\text{dR}}(X,\hat{T})$ with
$\CE_{\text{dR}}(M)^{(1)}_{X\setminus T}=(\det M_{X\setminus
T})^{\otimes -1}$. It is functorial with respect to isomorphisms of
$M$'s, has local origin, and lies in $\CL^\Phi_{\text{dR}}(X,T)$ if
$T$ is compatible with $M$. For proper $X$  there is a canonical
isomorphism of $k$-lines $\eta_{\text{dR}}:
\CE_{\text{dR}}(M)(X)\iso \det H^\cdot_{\text{dR}}(X,M)$. The
construction is compatible with base change of  $k$,
 filtrations on $M$, and direct images for finite morphisms of $X$'s.
}\bigbreak

We construct $\CE_{\text{dR}}(M)$ as a  factorization $\CO$-line in
2.5, and define a de Rham structure on it  in 2.10. The
identification $\CE_{\text{dR}}(M)^{(1)}_{X\setminus T}\iso (\det
M_{X\setminus T})^{\otimes -1}$ is established in (2.6.1); we check
that it is horizontal in 2.11. $\eta_{\text{dR}}$ is defined in 2.7,
the compatibilities are discussed in 2.8. The compatibility of $T$
with $M$ becomes relevant only in 2.10. Let us begin with necessary
preliminaries.

\subsection{}
{\it $\CL$-groupoids and $\CL$-torsors: a dictionary.} Let $\CL$ be
a Picard groupoid with the product operation $\otimes$.    Below
``{\it $\CL$-groupoid"} means ``enriched category over $\CL$".  Thus
this is a collection of objects $J$ and a rule which assigns to
every  $j, j'\in J$ an object $\lambda (j/j') \in\CL$, and to every
$j,j',j''\in J$ a {\it composition} isomorphism $\lambda
(j/j')\otimes\lambda (j'/j'')\iso \lambda (j/j'')$ which satisfies
associativity property. Then $J$ is automatically a mere groupoid
with $\Hom (j',j):= \Hom (1_\CL ,\lambda (j/j'))$.

Let $J_1$, $J_2$ be  $\CL$-groupoids.  Their {\it tensor product }
$J_1 \otimes J_2$ is an $\CL$-groupoid whose objects $j_1 \otimes
j_2$ correspond to  pairs $j_1\in J_1$, $j_2\in J_2$,  $\lambda
(j_1\otimes j_2 /j'_1\otimes j'_2 ):=  \lambda (j_1
/j'_1)\otimes\lambda (j_2/j'_2)$, and the composition $\lambda
(j_1\otimes j_2 /j'_1\otimes j'_2 )\otimes \lambda (j'_1\otimes j'_2
/j''_1\otimes j''_2 )\to \lambda (j_1\otimes j_2 /j''_1\otimes j''_2
)$ equal to $ (\lambda (j_1 /j'_1)\otimes\lambda (j_2/j'_2)) \otimes
(\lambda (j'_1 /j''_1)\otimes\lambda (j'_2/j''_2) ) \to (\lambda
(j_1 /j'_1)\otimes\lambda (j'_1 /j''_1))\otimes (\lambda (j_2/j'_2))
\otimes\lambda (j'_2/j''_2) ) \to  \lambda (j_1
/j''_1)\otimes\lambda (j_2/j''_2)$ where the first arrow is the
commutativity constraint, the second one is the tensor product of
the composition maps for $J_1$, $J_2$.

An {\it $\CL$-morphism}  $\phi : J_1 \to J_2$ is an $\CL$-enriched
functor, i.e., rule which assigns to every $j\in J_1$ an object
$\phi (j)\in J_2$, and to every $j,j'\in J_1$ an identification
$\phi : \lambda (j/j')\iso \lambda (\phi (j)/\phi (j'))$ compatible
with the composition. Such a $\phi$ yields a morphism of mere
groupoids $J_1 \to J_2$. All $\CL$-morphisms form naturally an
$\CL$-groupoid $\Hom_\CL (J_1 ,J_2 )$. Precisely, there is an
$\CL$-groupoid structure on $\Hom_\CL (J_1 ,J_2 )$ together with an
$\CL$-morphism  $\Hom_\CL (J_1 ,J_2 )\otimes J_1 \to J_2$ that lifts
the action map $\Hom_\CL (J_1 ,J_2 )\times J_1 \to J_2$ of mere
groupoids, and such pair is unique (up to a unique 2-isomorphism).
The composition $\Hom_\CL (J_2 ,J_3 )\times \Hom_\CL (J_1 ,J_2 )\to
\Hom_\CL (J_1 ,J_3 )$ lifts naturally to a morphism of
$\CL$-groupoids $\Hom_\CL (J_2 ,J_3 )\otimes\Hom_\CL (J_1 ,J_2 )\to
\Hom_\CL (J_1 ,J_3 )$, etc.

For an $\CL$-groupoid $J$ its {\it inverse} $J^{\otimes -1}$ is an
$\CL$-groupoid whose objects are in bijection $j\leftrightarrow
j^{\otimes -1}$ with elements of $J$, and $\lambda (j^{\otimes
-1}/j^{\prime \otimes -1})= \lambda (j'/j)$.

There are two equivalent ways to define the notion of {\it
$\CL$-torsor}: (a) This is a mere groupoid $F$ equipped with a
$\CL$-action, i.e., a functor $\otimes : \CL \times F \to F$
together with an associativity constraint, such that for some (hence
every) object $f\in F$ the functor $\CL \to F$, $\ell \mapsto \ell
\otimes f$, is an equivalence of groupoids; (b) This is an
$\CL$-groupoid such that the image of $\lambda$ meets every
isomorphism class in $\CL$.  To pass from (a) to (b), we lift the
groupoid structure on $F$ to $\CL$-groupoid one with  $\lambda
(f/f'):= f\otimes f^{\prime \otimes -1}$ (the latter  is an object
of $\CL$ together with an isomorphism $(f\otimes f^{\prime \otimes
-1})\otimes f' \iso f$; the pair is defined uniquely up to a unique
isomorphism).

For a non-empty  $\CL$-groupoid $J$    and an $\CL$-torsor $F$  both
$\CL$-groupoids $F \otimes J$ and  Hom$_\CL ( J,F )$ are
$\CL$-torsors.  In particular, we have the product and ratio of
$\CL$-torsors (with $\CL$ being a unit). Notice that there is a
natural equivalence $F_1 \otimes F_2^{\otimes -1} \iso$   Hom$_\CL (
F_2 ,F_1 )$ which assigns to
 $f_1 \otimes f_2^{\otimes -1}$ the $\CL$-morphism $f'_2 \mapsto \lambda (f'_2 /f_2 )\otimes f_1$.

{\it Remarks.} (i)   For any non-empty $\CL$-groupoid $J$ the
$\CL$-morphism $J\to \CL \otimes J$, $j\mapsto 1_{\CL}\otimes j$, is
a universal $\CL$-morphism to an $\CL$-torsor.

(ii) Every $\CL$-morphism  between $\CL$-torsors  is an equivalence.
Thus for non-empty $J_i$'s every $\CL$-morphism $J_1 \to J_2$ yields
an equivalence of $\CL$-torsors Hom$_\CL ( J_2 ,\CL )\iso$ Hom$_\CL
( J_1 ,\CL )$ and $\CL \otimes J_1 \iso \CL\otimes J_2$; in
particular,  the $\CL$-torsor Hom$_\CL ( J,\CL )$ does not change if
we replace $J$ by any its non-empty subset. E.g.~ every $j\in J$
yields an identification  of $\CL$-torsors Hom$_\CL ( J,\CL )\iso
\CL$, $\lambda\mapsto \lambda (j)$, and $\CL \iso \CL \otimes J$,
$\ell \mapsto \ell\otimes J$.

(iii) If $F_i = $ Hom$_\CL (J_i ,\CL )$ where $J_i$ are any
$\CL$-groupoids, then $F_1 \otimes F_2^{\otimes -1}$ identifies
naturally with the $\CL$-torsor whose objects are maps $\mu : J_1
\times J_2 \to \CL$, $(j_1 ,j_2) \mapsto \mu (j_1 /j_2 )$,  together
with identifications $\lambda ( j'_1 /j_1 )\otimes \mu (j_1 /j_2
)\otimes\lambda (j_2 /j'_2 )\iso \mu (j'_1 /j'_2 )$ which are
associative with respect to the composition of $\lambda$ on both
$J_i$. Here $f_1 \otimes f_2^{\otimes -1}$ corresponds to $\mu (j_1
/j_2 ):=f_1 (j_1 )\otimes f_2 (j_2 )^{\otimes -1}$.

Suppose a group $G$ acts on a non-empty $\CL$-groupoid $J$. Then it
acts on the $\CL$-torsor $\CL \otimes J$ by transport of structure.
We can view this action as a monoidal functor $g\mapsto \lambda_g$
from $G$ (considered as a discrete monoidal category) into the
monoidal category $\End_{\CL}(\CL \otimes J)$, which is naturally
equivalent to $\CL$. Explicitly,  $\lambda_g :=\lambda (g/\id_J
)\iso \lambda (g(j)/j)$, $j\in J$;   the product isomorphism
$\lambda_{g_1}\otimes\lambda_{g_2}\iso \lambda_{g_1 g_2 }$ is the
composition $\lambda (g_1 (g_2 (j))/g_2 (j))\otimes \lambda (g_2
(j)/j))\iso \lambda ((g_1 g_2 )(j)/j)$. A monoidal functor $G\to
\CL$ is sometimes called {\it (central) $\CL$-extension} $G^\flat$
of $G$.\footnote{ If $\CL$ is the Picard groupoid of $A$-torsors,
$A$ is an abelian group, then $G^\flat$ amounts to a central
extension of $G$ by $A$, which is the reason for the terminology.}

The group $G$ acts on $G^\flat$ in adjoint way. Namely, for $h\in G$
the isomorphism Ad$_h : \lambda_g \iso \lambda_{hgh^{-1}}$ is the
composition $\lambda_g \iso \lambda_{g}\otimes \lambda_h
\otimes\lambda_{h^{-1} }\iso \lambda_{h}\otimes \lambda_{g}
\otimes\lambda_{h^{-1} }\iso \lambda_{hgh^{-1}}$, the first arrow is
the tensoring with the inverse to the composition map
$1_{\CL}\buildrel\sim\over\leftarrow\lambda_h \otimes\lambda_{h^{-1}
}$, the second one is the commutativity constraint. Equivalently,
Ad$_h$ is determined by the condition that the composition $\lambda
(g(j)/j) \iso \lambda_g \buildrel{\text{Ad}_h}\over\lra
\lambda_{ghg^{-1}} \iso \lambda (hgh^{-1}(h(j))/h(j))=\lambda (h
(g(j))/h(j))$ coincides with the action of $h$.

For commuting $g,h\in G$  we denote by $\{ g ,h \}^\flat \in\pi_1
(\CL ):=\Aut_\CL (1_\CL )$ their commutant in $G^\flat$, i.e., the
action of Ad$_g$ on $\lambda_h$ or the ratio of  $\lambda_g
\otimes\lambda_{h}\to \lambda_{gh}$ and $\lambda_g
\otimes\lambda_{h}\to \lambda_{h} \otimes\lambda_{g}\to
\lambda_{hg}=\lambda_{gh}$ where the first and the last arrows are
the product, the middle one is the commutativity constraint
(cf.~\cite{BBE} A5).

\subsection{} {\it A digression on lattices and relative determinants}
 (see e.g.~\cite{Dr} \S 5).

Let $S$ be a scheme, $P$ be a relative effective divisor in $X_S /S$
finite over $S$. Let $E$ be any quasi-coherent $\CO_{X_S}$-module
such that for some neighborhood $V$ of $P$ the restriction of $E$ to
$V\setminus P$  is coherent and locally free.

A {\it $P$-lattice in} $E$ is an $\CO_{X_S}$-submodule $L$  of $E$
which is  locally free (hence coherent) on a neighborhood of $P$,
and equals  $E$ on $X_{S}\setminus P$. Denote by $\Lambda_P (E)$ the
set of $P$-lattices in $E$. We assume that it is non-empty. Then
$\Lambda_P (E)$  is directed by the inclusion ordering. Since $P$ is
finite over $S$, for every $P$-lattices $L \supset L'$ the
$\CO_S$-module $  \pi_{*} (L /L' )$ is locally free of finite
rank.\footnote{Let us show that  $  \pi_{*} (L /L' )$ is
$\CO_S$-flat. We can assume that $X$ is affine, $L$, $L'$ are
locally free. Since $  \pi_{*} (L /L' )   =  \pi_{*} (L )/ \pi_{*}
(L')$ and   $\pi_* L$, $\pi_* L'$ are $\CO_S$-flat,  it suffices to
check that for any geometric point $s$ of $S$ the map $\pi_* ( L')_s
\to \pi_* (L)_s $ is injective. This is clear, since $ \pi_*
(L^{(\prime )})_s =  \Gamma (X_s ,L^{(\prime )}_s )$ and  $L'_s \to
L_s$ is injective (being an isomorphism at the generic points of
$X_s$).}

$\Lambda_P (E)$ carries a natural $\CL_\CO (S)$-groupoid structure.
Namely, for $P$-lattices $L$, $L'$  one has their  {\it relative
determinant} $\lambda_P (L / L' ) \in\CL_\CO (S)$; for  $L$, $L'$,
$L''$ there is a canonical composition isomorphism
\begin{equation}{\lambda}_P (L / L' )\otimes {\lambda}_P (L' / L'' )
\iso {\lambda}_P (L / L'' ) \label{ 2.3.1}\end{equation} which
satisfies associativity property. This datum is uniquely determined
by a demand that for $L \supset L'$ one has $\lambda_P (L/L'):=\det
\pi_{*} (L / L' )$, and for $L\supset L' \supset L''$ the
composition is the standard isomorphism defined by the short exact
sequence $0\to \pi_{*}(L' /L'' )\to   \pi_{*}(L /L'' )\to\pi_{*}(L
/L' )\to 0$.

We denote by $\CD et_{P/S}(E)$ the  $\CL_\CO (S)$-torsor
Hom$_{\CL_\CO (S)}( \Lambda_P (E), \CL_\CO (S)) $. Its objects,
referred to as {\it determinant theories on $E $ at $P$}, are rules
$\lambda$ that assigns to every $L\in \Lambda_P (E) $ a line
$\lambda (L)\in\CL_\CO (S)$ together with identifications $\lambda_P
(L/L') \otimes\lambda (L' ) \iso \lambda (L)$ compatible with
(2.3.1). Here one can replace $\Lambda_P (E) $ by any its non-empty
subset. $\CD et_{P/S}(E)$ is compatible with base change change.

For  quasi-coherent $\CO_{X}$-modules  $E_1$, $E_2$ as above set
$$\CD et_{P/S} (E_1 /E_2 ) :=\CD et_{P/S} (E_1 )\otimes\CD et_{P/S}
(E_2 )^{\otimes -1}.$$ By Remark (iii) in 2.2, objects of this
$\CL_\CO (S)$-torsor, referred to as {\it relative determinant
theories on $E_1 / E_2 $ at $P$}, can be viewed as rules $\mu$ which
assign to every $P$-lattices $L_1$ in $E_1$ and $L_2$ in $E_{2 } $
a line  $\mu (L_1 /L_2 ) \in \CL_\CO (S)$ together with natural
identifications
\begin{equation}
{\lambda}_P (L'_1 /L_1 ) \otimes  \mu (L_1 /L_2) \otimes {\lambda}_P (L_2 /L'_2) \iso  \mu (L'_1/ L'_2 ) \label{ 2.3.2}
\end{equation}
which are associative with respect to composition  (2.3.1) on the
two sides. We can also restrict ourselves to $L_i$ in any non-empty
subsets of $\Lambda_P (E_i )$.

{\it Remarks.} (i)  Let $E (\infty P):=\limright E(nP)$ be the
localization of $E$ with respect to  $P$. An evident morphism of
$\CL_\CO (S)$-groupoids $\Lambda_P (E)\to \Lambda_P (E(\infty P))$
yields an equivalence $\CD et_{P/S}(E)\iso \CD et_{P/S}(E(\infty
P))$; same for relative determinant theories. Thus $\CD et_{P/S}(E)$
depends only on the restriction of $E$ to the complement of $ P$.
Notice that $\Lambda_P (E(\infty P))$ is directed by both inclusion
ordering and the opposite one.

(ii) Let us order the set of pairs of $P$-lattices $(L_1 ,L_2 )$ by
the product of either of the inclusion orderings. Let $I$ be any of
its  directed subsets. Suppose we have a rule $\lambda$ that assigns
to every $(L_1 ,L_2 )\in I$ a line $\lambda (L_1 /L_2 )\in\CL_\CO
(S)$ together with natural identifications (2.3.2) defined for $(L_1
,L_2 )\ge (L'_1 ,L'_2)$ and associative for    $(L_1,L_2 )\ge (L'_1
,L'_2 )\ge (L''_1 ,L''_2 )$. Then $\lambda$ extends uniquely to a
relative determinant theory.

(iii) The group $\Aut (E_{V\setminus P})$ acts on  $\Lambda_P
(E(\infty P))$ as on an $\CL_\CO (S)$-groupoid. As in 2.2, this
defines an $\CL_\CO (S)$-extension $\Aut (E_{V\setminus P})^\flat$
of $\Aut (E_{V\setminus P})$. Thus for commuting $g',g\in \Aut
(E_{V\setminus P})$ we have $\{ g',g\}^\flat = \{ g',g\}^\flat_P
\in\CO^\times (S)$. Example: if $g'$ is multiplication by a function
$f\in \CO^\times (V\setminus P)$, then $\{ g',g\}^\flat_P$ equals
the Contou-Carr\`ere symbol $\{ \det g,f\}_P$ at $P$ (see
e.g.~\cite{BBE} 3.3).\footnote{A short proof: Both expressions are
compatible with base change. Since the datum of $E$, $V$, $P$, $g$,
$f$ can be extended $S$-locally to a smooth base, we can assume that
$S$ is smooth. Then it suffices to check the equality at the generic
point of $S$, where the Contou-Carr\`ere symbol  is the usual tame
symbol. The rest is a standard computation.} In particular, if
$f\in\CO^\times (V)$, then $\{ g',g\}^\flat_P=f(-\text{div}(\det
g))$; here  $\text{div} (\det g) $ is the part of the divisor
supported on $V$, i.e., at $P$.

$P$-lattices have local nature with respect to $P$. In particular,
if $P$ is the disjoint  sum of components $P_\alpha$, then a
$P$-lattice $L$  amounts to a collection of $P_\alpha$-lattices
$L_\alpha$, and there is an evident canonical factorization
isomorphism
 \begin{equation} \otimes\, {\lambda}_{P_\alpha} (L_\alpha /L'_\alpha )
   \iso {\lambda}_P (L/L') \label{ 2.3.3   }\end{equation} compatible with composition (2.3.1).
   Therefore $\CD et_{P/S}(E)$ is the $\CL_\CO (S)$-torsor product of $\CD et_{P_\alpha /S}(E)$. So
every collection of determinant theories $\lambda_\alpha$ on $E$ at
$P_\alpha$ yields a determinant theory $\otimes \lambda_\alpha$ on
$M$ at $P$, $(\otimes\lambda_\alpha )(L):=\otimes \lambda_\alpha
(L_\alpha )$. Same for relative determinant theories.

\subsection{} We return to the story of 1.1.  Suppose we have $(D,c)\in\fD (S)$.
 Let us apply the format of 2.3 to $E_1 =M_S =M\boxtimes\CO_S$, $E_2 = \omega M_S
 :=(\omega\otimes M)_S$, and $P= P_{D,c}$.
  The connection  $\nabla = \nabla_M : M\to \omega M$ yields  a relative
  determinant theory $\mu^\nabla_P =\mu (M)^\nabla_P \in \CD et_{P/S}
  (M/  \omega M):=  \CD et_{P/S}
  (M_S /  \omega M_S )$ defined as follows.

  Let $L$, $L_\omega$ be $P$-lattices
  in $ M_S$, $ \omega  M_S$ such that $\nabla (L)\subset L_\omega$. Let
   $\CC (L,L_\omega )= \CC (L,L_\omega )_{M,P} $ be the complex $M_S /L
\buildrel{\nabla}\over\to \omega  M_S /L_\omega $ in degrees $-1$
and $0$, i.e., it is the quotient $dR(M)\boxtimes\CO_{S}   /\CC one
( L  \buildrel{\nabla}\over\to L_\omega )$. Then
   $\pi_* \CC (L,L_\omega )$ is
    a complex of quasi-coherent $\CO_S$-modules.

\bigbreak\noindent\textbf{Lemma.}\emph{  $\pi_* \CC (L,L_\omega )$
has $\CO_S$-coherent cohomology. }\bigbreak

{\it Proof.} We can assume that $T^c_S =T_S$, since  the  assertion
is $S$-local. Its validity  does not depend on the choice of $L$,
$L_\omega$. Take them to be ``constant" $T$-lattices  equal to
  $M$, $\omega M$ on $X\setminus T$, and we are reduced to
  the case of $S=\Spec \, k$, $P=T$.

  Let $\bar{X}$ be the smooth projective curve
  that contains $X$, $T^{\infty}:= \bar{X}\setminus X$. Let us
  extend $M$ to a holonomic $\CD$-module on  $\bar{X}$, which we denote also by
  $M$; let $L$, $L_{\omega}$ be $T\cup T^{\infty}$-lattices as above. Since
  $\CC (L,L_{\omega})_{M,T\cup T^{\infty}}=
  \CC (L,L_{\omega})_{M,T}\oplus \CC (L,L_{\omega})_{M, T^{\infty}} $,
  it suffices to check that $\pi_* \CC (L,L_{\omega})_{M,T\cup T^{\infty}}$ has finite-dimensional cohomology. Therefore
we are reduced to the case of projective $X$.

Now the lemma follows since $\CC (L,L_{\omega})=dR(M) /\CC one
(L\buildrel{\nabla}\over\to L_{\omega})$, and the cohomology of $X$
with coefficients in $dR(M)$, $L$, $L_{\omega}$ are finite
dimensional. \hfill$\square$

Pairs $(L,L_\omega )$ as above form a directed set $I_\CC$ as in
Remark (ii) in 2.3. Set
  \begin{equation}
\mu^\nabla_P (L/L_\omega ):=\det\pi_* \CC (L,L_\omega ). \label{
2.4.1}\end{equation}
If $(L',L'_\omega )\in I_\CC$ is  such that
$L\supset L'$, $L_\omega \supset L'_\omega$, then the short exact
sequence $0\to  \CC one ( L/L' \buildrel{\nabla}\over\to L_\omega
/L'_\omega )\to \CC (L',L'_\omega )\to\CC (L,L_\omega )\to 0$
together with the identification $\det \pi_*
 \CC one ( L/L' \buildrel{\nabla}\over\to L_\omega /L'_\omega )= \det \pi_* (L_\omega /L'_\omega )\otimes \det \pi_* (L/L')^{\otimes -1}$
yields an isomorphism
\begin{equation}{\lambda}_P  (L'/L)\otimes
\mu^\nabla_P (L/L_\omega )\otimes{\lambda}_P  (L_\omega /L'_\omega
)\iso \mu^\nabla_P (L'/L'_\omega) \label{ 2.4.2}
\end{equation}
 which
evidently satisfies associativity property. By Remark (ii) in 2.3,
we have defined $\mu_P^\nabla =\mu (M/\omega M)_P^\nabla \in \CD
et_{P/S} (M/\omega M)$.

{\it Remarks.} (i) Sometimes it is convenient to consider a smaller
directed set formed by pairs $(L, L_\omega )$ that equal $M$,
$\omega M$ outside $T^c_S$, or  its subset of those $(L, L_\omega )$
that are locally constant with respect to $S$.

(ii) Denote by $j_T$ the embedding $X\setminus T \hra X$. Then $j_{
T*}M:=j_{ T*}M_{X\setminus T}$ is holonomic $\CD$-module as well; by
Remark (i) from 2.3, one has $\CD et_{P/S}(M/\omega M)=\CD et_{P/S}
(  j_{T*}M /\omega j_{T*}M)$. Suppose $T^c_S =T_S$. Then $L$,
$L_{\omega}$ are $P$-lattices in $j_{ T*}M$, $\omega j_{T*} M$, and
$\CC one (\CC (L,L_{\omega})_{M,P},\CC (L,L_{\omega})_{j_{T*}M,P}) $
$\iso dR(\CC one (M\to j_{T*}M))\otimes\CO_{S}$. Thus $\mu(M/\omega
M)_P^\nabla =\mu (j_{T*}M/\omega j_{T*}M)_P^\nabla\otimes \det
R\Gamma_{\!\text{dR} \, T}(X,M )$.

The above construction has local nature with respect to $P$. If $P$
is disjoint  sum of components $P_\alpha$, and $L$, $L_\omega$
correspond to collections of $P_\alpha$-lattices $L_\alpha$,
$L_{\omega \alpha}$, then $\CC (L,L_\omega )= \oplus\, \CC (L_\alpha
,L_{\omega\alpha})$.  Passing to the determinants, we see that
\begin{equation}\mu(M/\omega M)^\nabla_P = \otimes \mu (M/\omega
M)^\nabla_{P_\alpha}. \label{ 2.4.3}\end{equation}

\subsection{}  Now we can construct the promised  $\CE_{\text{dR}}(M)$ as a factorization $\CO$-line.

For $(D,c )\in \fD (S)$, any trivialization $\nu$ of $\omega (D)$ on
a neighborhood $V$ of $P=P_{D,c}$ yields naturally  $\mu^\nu_P =\mu
(M/\omega M)^\nu_P \in \CD et_{P/S}(M/\omega M)$. Namely,  the
multiplication by $\nu$ isomorphism $M_{V\setminus P} \iso \omega
M_{V\setminus P}$ identifies the $\CL_\CO (S)$-groupoids $\Lambda_P
(M(\infty P))\iso  \Lambda_P (\omega M(\infty P))$;  passing to $\CD
et_{P/S}$, we get   $\mu^\nu_P$. Explicitly, for every $P$-lattices
$L$, $L_\omega$ one has a canonical identification
\begin{equation}\mu^\nu_P (L/L_\omega )\iso \lambda_P (\nu
L/L_\omega )= \lambda_P (\omega L(D)/L_\omega ), \label{
2.5.1}\end{equation} and identifications (2.3.2) are (2.3.1)
combined with the isomorphism $$\nu_{L/L'} : \lambda_P (L/L')\iso
\lambda_P ( \nu L/ \nu L')$$ that comes from multiplication by $\nu$.

The r.h.s.~of (2.5.1) does not depend on $\nu$. Thus every $L$
provides an identification $e_L$ between $\mu_P^\nu$ for all
trivializations $\nu$ of $\omega (D)$ near $P$; it is characterized
by property that (2.5.1) transforms $e_L (L/L_\omega )$ into the
identity map for $\lambda_P (\omega L(D)/L_\omega )$.

{\it Exercise.} Let $\nu_1 ,\nu_2$ be any trivializations of $\omega
(D)$ on $V$; set  $f:=\nu_2 /\nu_1\in\CO^\times (V)$. Consider the
identifications $e_L ,e_{L'}: \mu^{\nu_1}_P \iso \mu^{\nu_2}_P$.
Show that \begin{equation}e_{L' }= f (\text{div}(L/L'))e_L .\label{
2.5.2}\end{equation} Here div$(L/L'):=$ div$(\phi_{L'}/\phi_{L})$,
where $\phi_{L}$, $\phi_{L'}$ are  local trivializations of the line
bundles $\det (L)$, $\det (L')$ at $P$; this is a relative Cartier
divisor supported  at $P$.

If $\nu_1 |_P =\nu_2 |_P$, i.e., $f|_P=1$, then we define a {\it
canonical} identification  \begin{equation}e: \mu^{\nu_1}_P \iso
\mu^{\nu_2}_P \label{ 2.5.3}\end{equation} as follows. The function
$f(D)\in \CO^\times (S)$ equals 1 on $S_{\text{red}}$ since $f|_P
=1$. Let $f(D)^{\frac{1}{2}}$ be the branch of the root that equals
1 on $S_{\text{red}}$. Pick a lattice $L_0$ which equals $M$ off
$T^c_S$ and is  $S$-locally constant. Then\footnote{The reason for
the normalization will become clear in 2.10.}  $e :=
f(D)^{\frac{\text{rk}(M)}{2}}e_{L_0} $. By (2.5.2), the auxiliary
choice of $L_0$ is irrelevant: indeed, if $L_0$, $L'_0$ that satisfy
our condition, then $f($div$(L_0/L'_0))=1$, for the divisor
div$(L_0/L'_0)$ is $S$-locally constant and supported on $T^c_S$,
and $f|_{T^c_S}=1$. The identifications $e$ are evidently
transitive.

\medskip

{\it Remark.} Suppose $L$ is an arbitrary lattice. Then  $$e= f(\text{div} (L/L_0
))f(D)^{\text{rk}(M)/2}e_L$$ where $L_0$ is any lattice as above (see (2.5.2)
and (2.5.3)).

Our $e$ provides a canonical identification of $\CO_S$-lines
\begin{equation}
r:= \id_{\mu^\nabla_P}\otimes e^{\otimes -1} : \mu^\nabla_P \otimes
(\mu^{\nu_1}_P )^{\otimes -1} \iso \mu^\nabla_P \otimes
(\mu^{\nu_2}_P )^{\otimes -1}. \label{ 2.5.4}
\end{equation}
Therefore for $(D,c,\nu_P )\in\fD^\diamond (S)$ the $\CO_S$-line
$\mu^\nabla_P \otimes (\mu^\nu_P )^{\otimes -1}$  does not depend on
the choice of $\nu$  such that $\nu|_P =\nu_P$. Set
\begin{equation}\CE_{\text{dR}}(M)_{(D,c,\nu_P )}:=\mu^\nabla_P
\otimes (\mu^\nu_P )^{\otimes -1} . \label{ 2.5.5}\end{equation} The
construction is compatible with base change, so $\CE_{\text{dR}}(M)$
is  an $\CO$-line on $\fD^\diamond$. By (2.4.3) and similar property
of $\mu^\nu_P$, it carries a factorization structure. So we have
defined $\CE_{\text{dR}}(M)\in\CL_\CO^\Phi  (X,T )$.

{\it Example.} If $M$ is supported at $T$, then
$\CE_{\text{dR}}(M)_{(D,c,\nu )}=\det
R\Gamma_{\!\text{dR}\,T^c}(X,M)$.

{\it Summary.} Suppose we have  $(D,c,\nu_P )\in\fD^\diamond (S)$.
By (2.5.1),  every $L$ and $\nu$ such that $\nu|_P =\nu_P$ yields an
isomorphism \begin{equation}nr_{L,\nu }:
\CE_{\text{dR}}(M)_{(D,c,\nu_P )}\iso \mu_P^\nabla (L/\omega L(D)).
\label{ 2.5.6}\end{equation} If $f|_P =1$, then, by Remark,
\begin{equation}
r_{L,f\nu} =f(\text{div}(L_0 /L))
f(D)^{-\frac{\text{rk}(M)}{2}}a_{L,\nu }. \label{
2.5.7}
\end{equation}
In particular, $r_{L,f\nu} =r_{L,\nu }$ for
reduced $S$.

\subsection{}\label{2.6}
\bigbreak\noindent\textbf{Lemma.}\emph{ (i) On $X\setminus T$ there
is a canonical isomorphism
\begin{equation}\CE_{\text{dR}}(M)^{(1)}_{X\setminus T} \iso (\det
M_{X\setminus T})^{\otimes -1}.  \label{ 2.6.1}\end{equation} \\
(ii) Suppose $T' \subset T$ is such that $M$ is smooth at
$T\setminus T'$. Let $\CE_{\text{dR}}(M)' $ be the restriction to
$(X,T)$ of the $\varepsilon$-line of $M$  on $(X,T')$. Then there is
a canonical identification $\CE_{\text{dR}}(M)\iso
\CE_{\text{dR}}(M)'$.}\bigbreak

{\it Proof.} (i) Recall that $\CE^{(1)}_{X\setminus T}
=\CE_{(D,0,\nu )}$ where $S=X\setminus T$, $D=\Delta =P$ is the
diagonal divisor, and $\nu$ is (the principal part of) a form with
logarithmic singularity at $\Delta$ with residue 1. Take
$L=M_{X\setminus T}$. Then $\CC (L,\omega L)=0$, hence
$\CE_{\text{dR}}(M)^{(1)}_{X\setminus T}
\buildrel{r_{L,\nu}}\over\lra \mu^\nabla_P (L/\omega L(D))\iso
\lambda_P (\omega L/\omega L(D))=  \det \pi_* (\omega L(D)/\omega
L)^{\otimes -1}$. Now (2.6.1) is this isomorphism followed by the
residue map $\pi_* (\omega L(D)/\omega L ) \iso L=M_{X\setminus T}$.

(ii) Evident.   \hfill$\square$

\subsection{}\label{2.7}

 \bigbreak\noindent\textbf{Proposition.}\emph{ For
proper $X$ there is a canonical isomorphism
\begin{equation}\eta_{\text{dR}}: \CE_{\text{dR}}(M)(X)\iso
\det H^\cdot_{\text{dR}}(X,M). \label{ 2.7.1}\end{equation}
}\bigbreak

{\it Proof.}   By 1.4, (2.7.1) amounts to
 a natural identification $\eta_{\text{dR}}:\CE_{\text{dR}}(M)_\nu  \iso
 \det H^\cdot_{\text{dR}}(X,M)\otimes\CO_S$ which is  defined
 for every $S$-family of rational 1-forms $\nu$ on $X$ and is compatible with base change.

Set $P=P_{\text{div}(\nu),1}=T\cup |$div$(\nu )|$. Since $X$ is
proper, for a $P$-lattice $L$ in $M(\infty P )$ the complex of
$\CO_S$-modules $R\pi_* (L)$ is perfect; set $\lambda (L):=\det
R\pi_* (L)$. Then $\lambda $ is a determinant theory on $M$ at $P$
in an evident manner. Replacing $M$ by $\omega M$, we get
$\lambda_\omega \in \CD et_{P/S}(\omega M)$, hence $\lambda
\otimes\lambda_\omega^{\otimes -1}\in\CD et_{P/S}(M/\omega M)$. One
has an isomorphism \begin{equation}\mu^\nu_P \iso \lambda
\otimes\lambda_\omega^{\otimes -1}, \label{ 2.7.2}\end{equation}
namely, $ \mu^\nu_P (L/L_\omega ):= \lambda_P (\nu L/L_\omega
)=\lambda_\omega (\nu L)\otimes\lambda_\omega (L_\omega )^{\otimes
-1} \iso \lambda (L)\otimes\lambda_\omega (L_\omega )^{\otimes -1}=
(\lambda \otimes\lambda_\omega^{\otimes -1} )(L/L_\omega )$ where
$\iso$ comes from  isomorphism  $\nu^{-1} :  \nu L \iso  L$.

For $P$-lattices $L$ in $ M$, $L_\omega $ in $ \omega M$ with
$\nabla ( L )\subset L_\omega$ (see 2.4), set $dR (L,L_\omega ):=
\CC one (L\buildrel{\nabla}\over\to L_\omega )\subset dR (M)$. Since
$dR (M)/ dR (L, L_\omega )= \CC (L,L_\omega )$, we see that $dR (M)$
carries a 3-step filtration with successive quotients $L_\omega$,
$L[1]$, $\CC (L,L_\omega )$. Applying $\det R\pi_*$, we get an
isomorphism

\begin{equation}
\det R\pi_* \CC (L,L_\omega )\otimes\lambda (L)^{\otimes -1}\otimes
\lambda ( L_\omega ) \iso \det H^\cdot_{\text{dR}}(X,M)\otimes\CO_S
. \label{ 2.7.3}
\end{equation}
To get $ \eta_{\text{dR}}$, we combine it with (2.7.2) (and
(2.5.5)). The construction does not   depend on the auxiliary choice
 of $L$, $L_\omega$.  \hfill$\square$

\medskip

{\it Example.} Suppose  $M$ is a $\CD$-module on $\Bbb P^1$ with
regular singularities at 0 and $\infty$, where it is, respectively,
the $*$- and the !-extension. Since $R\Gamma_{\!\text{dR}}(\Bbb
P^1,M)=0$, our $\eta_{\text{dR}}$ is an isomorphism \begin{equation}
\eta_{\text{dR}} : \CE_{\text{dR}}(M)_{(0,t^{-1}dt )}\otimes
\CE_{\text{dR}}(M)_{(\infty ,t^{-1}dt )}\iso k . \label{
2.7.4}\end{equation} To compute it explicitly, pick a
$t\partial_t$-invariant vector subspace $V$ of $\Gamma (\Bbb P^1
\setminus \{ 0,\infty \},M)$, which freely generates $M_{\Bbb P^1
\setminus \{ 0,\infty \}}$ as an $\CO$-module and such that the only
possible integral eigenvalue of $t\partial_t$ on $V$ is 0. Set $L:=
\CO_{\Bbb P^1}(-(\infty ))\otimes V$, $L_\omega := t^{-1}dt \otimes
L=\omega_{\Bbb P^1}((0)+(\infty) )\otimes L $.  The condition on $V$
implies that there are $\CO$-linear embeddings $i : L\hra M$,
$i_\omega : L_\omega \hra \omega M$, which extend the evident
isomorphisms on $\Bbb P^1 \setminus \{ 0,\infty \}$,  such that
$\nabla (L)\subset L_\omega$ and $i_\omega (\phi \ell )=\phi i (\ell
)$ for any $\phi \in \omega _{\Bbb P^1}$, $\ell \in L$. Such $i$,
$i_\omega$ are unique. The complex $\CC (L,L_\omega )= \CC
(L,L_\omega )_0 \oplus \CC (L,L_\omega )_\infty$ is acyclic, so
(2.5.6) provides trivializations $\iota_0$ of
$\CE_{\text{dR}}(M)_{(0,t^{-1}dt )}$ and $\iota_\infty$ of
$\CE_{\text{dR}}(M)_{(\infty ,t^{-1}dt )}$.

\bigbreak\noindent\textbf{Lemma.}\emph{ One has $
\eta_{\text{dR}}(\iota_0 \otimes\iota_\infty )=1$.}\bigbreak

{\it Proof.}  The determinant of the complex $R\Gamma (\Bbb P^1 ,dR
(L,L_\omega ))$ has two natural trivializations $\alpha_1$,
$\alpha_2$: the first one comes since the complex is acyclic, the
second one  from the identification $\det R\Gamma (\Bbb P^1 ,dR (L,
L_\omega ))= \det R\Gamma (\Bbb P^1 , L_\omega )\otimes \det R\Gamma
(\Bbb P^1 , L)^{\otimes -1} $ and the multiplication by $t^{-1}dt$
isomorphism $L\iso L_\omega$. Now (2.7.3) identifies $ \iota_0
\otimes\iota_\infty \otimes \alpha_1$ with 1, and $ \iota_0
\otimes\iota_\infty \otimes \alpha_2$ with $\eta_{\text{dR}} (
\iota_0 \otimes\iota_\infty )$. Since
 $R\Gamma (\Bbb P^1 ,L)=R\Gamma (\Bbb P^1 ,L_\omega )=0$, one has $\alpha_1 =\alpha_2$; we are done. \hfill$\square$

\subsection{} The next constraints follow directly from the construction:

(i) For  a finite filtration $M_\cdot $ on $M$,  there is a
canonical isomorphism \begin{equation} \CE_{\text{dR}} ( M )\iso
\otimes\, \CE_{\text{dR}} (\gr_i M )  \label{ 2.8.1}\end{equation}
which satisfies transitivity property with respect to refinement of
the filtration.

{\it Remark.}  If $M=\oplus M_\alpha$, then every linear ordering of
the indices yields a filtration on $M$, hence an isomorphism
$\CE_{\text{dR}} ( M )\iso \otimes\, \CE_{\text{dR}} ( M_\alpha ) $.
This isomorphism does not depend on the ordering. Thus
$\CE_{\text{dR}}$ is a symmetric monoidal functor.

(ii) Let $\pi : (X',T') \to (X,T)$ be a finite morphism of pairs
(see Remark (i) in 1.2) which is \'etale over $X\setminus T$. As in
loc.~cit., we have a morphism of Picard groupoids $\pi_* :
\CL^\Phi_{\CO} (X',T')\to \CL^\Phi_{\CO} (X,T)$. We also have the
$\CD$-module direct image functor $\pi_*$ which  is exact. If $M'$
is a holonomic $\CD$-module on $(X',T')$, then $\pi_* M'$ is
holonomic $\CD$-module on $(X,T)$. Notice that $\pi_* M'$ coincides
with the ``naive" direct image $\pi_\cdot M'$ outside $T$, and  $dR
(\pi_* M')$ is canonically quasi-isomorphic to $\pi_\cdot dR (M')$
as a dg module over the de Rham dg algebra of $X$. Therefore one has
a canonical isomorphism \begin{equation}\CE_{\text{dR}} ( \pi_*
M')\iso \pi_* \CE_{\text{dR}} (  M')  \label{ 2.8.2}\end{equation}
 compatible with the composition of $\pi$'s and with (2.8.1).

{\it Exercise.}  Consider the standard isomorphism
$H^\cdot_{\text{dR}}(X,\pi_* M')\iso  H^\cdot_{\text{dR}}(X', M') $.
If $X$ is proper, then (2.8.2) yields an isomorphism
$\CE_{\text{dR}} ( \pi_* M')(X) \iso \CE_{\text{dR}} (  M')(X')$.
Show that $\eta_{\text{dR}}$'s identify the second isomorphism with
the determinant of the first one.

\subsection{}
{\it Another digression on lattices and relative determinants.} For
a Clifford algebra explanation of the  next constructions,  see
\cite{BBE} 2.14--2.17. In this subsection and the next one we use
$\Bbb Z/2$-graded  lines instead of $\Bbb Z$-graded ones; as in
1.2,  the corresponding Picard groupoids are marked by $'$.

Let $S$, $P$ be as in  2.3. Let $E$, $E^\circ$ be $\CO_{X}$-modules
as in 2.3, $V$ a neighborhood of $P$ such that both $E$, $E^\circ$
are locally free over $V\setminus P$,
  $\psi : E_{V\setminus P}\times E^\circ_{V\setminus P} \to
\omega_{V\setminus P/S}$ be a non-degenerate pairing. For a
$P$-lattice $L$ in $E(\infty P)$ its {\it $\psi$-dual} $L^\psi$ is
the $P$-lattice in $E^\circ (\infty P)$ such that $\psi$ is a
non-degenerate $\omega_{V/S}$-valued pairing  between $L_V$ and
$L_V^\psi$. The map
 $\tau_\psi : \Lambda_P (E(\infty P))\to \Lambda_P (E^\circ (\infty P))$,
 $L\mapsto \tau_\psi (L):= L^\psi$, is an order-reversing bijection. It lifts to an isomorphism
 of $\CL'_\CO (S)$-groupoids:

\bigbreak\noindent\textbf{Lemma.}\emph{  For every $L, L' \in
\Lambda_P (E(\infty P))$ there is a canonical isomorphism
\begin{equation}
\tau_\psi : {\lambda}_P (L/L')\iso    {\lambda}_P (L^{ \psi}/
L^{\prime \psi} ) \label{ 2.9.1}
\end{equation}
of $\Bbb Z/2$-graded lines compatible with the
composition.}\bigbreak

{\it Proof.} It suffices to define (2.9.1) for $L'\supset L$ and
check the compatibility with composition for $L''\supset L'\supset
L$.

The  pairing $\ell,\ell^\circ \mapsto (\ell ,\ell^\circ )_\psi
:=\Res_{P /S}\,\psi ( \ell,\ell^\circ )$  yields a  duality between
the vector bundles $\pi_* (L'/L )$ and $\pi_* ( L^{ \psi}/ L^{\prime
\psi} )$, hence a duality between the determinant lines $\lambda_P
(L'/L )\otimes \lambda_P ( L^{ \psi}/ L^{\prime \psi} )\iso \CO_S$,
$(\wedge \ell_i )\otimes (\wedge \ell_j^\circ )\mapsto
(-1)^{\frac{n(n-1)}{2}}\det (\ell_i ,\ell_j^\circ )_\psi$ where $n=
\text{rk}\, \pi_* (L'/L)$. Then (2.9.1) is characterized by the
property that $\id_{\lambda_P (L'/L)}\otimes\tau_\psi$ identifies
the latter pairing with the composition $\lambda_P (L'/L )\otimes
\lambda_P (L/L')\iso \CO_S$. The compatibility of $\tau_\psi$ with
composition is left to the reader. \hfill$\square$

{\it Remark.} Here is a sketch of a Clifford algebra interpretation
of $\tau_\psi$. Consider the Clifford $\CO_S$-algebra generated by
$\pi_* E_{V\setminus P}\oplus  \pi_* E^\circ_{V\setminus P}$
equipped with the hyperbolic form $\Res_{P/S} \psi$. Let $N$ be any
``invertible" continuous Clifford module. For any $P$-lattice $L$
the $\CO_S$-submodlule $N^{L\oplus L^\psi}$ of vectors  killed by
$L\oplus L^\psi$ lies in $\CL'_\CO (S)$, and $\lambda_P (L/L') =
N^{L\oplus L^\psi} \otimes (N^{L'\oplus L^{\prime\psi}})^{\otimes
-1}= \lambda_P (L^\psi /L^{\prime\psi}) $; the composition is
$\tau_\psi$.

Passing to $\CD et'_{P/S}$, our $\tau_\psi$ yields  $\mu^\psi \in\CD
et'_{P/S} (E/E^\circ )$. Thus for $L\in\Lambda_P (E(\infty P))$,
$L_\omega \in\Lambda_P (E^\circ (\infty P))$ one has $\mu^\psi
(L/L_\omega )= \lambda_P (L^\psi /L_\omega )$, and identifications
(2.3.2) are (2.3.1) combined with  $\tau_\psi$.

{\it Exercise.} Suppose we have another non-degenerate pairing $
E_{V\setminus P}\times E^\circ_{V\setminus P} \to \omega_{V\setminus
P/S}$; one can write it as $\psi_g (\cdot ,\cdot )=\psi (g\cdot
,\cdot )=\psi (\cdot ,{}^\psi \! g \cdot )$ where  $g\in\Aut
(E_{V\setminus P})$ and ${}^\psi \!g\in\Aut (E_{V\setminus P})$ is
the $\psi$-adjoint to $g$. Then $L^{\psi_g} = {}^\psi\!
g^{-1}(L^\psi )= (g(L))^\psi$ and
\begin{equation}\tau_{\psi_g}={}^\psi\!g^{-1} \tau_\psi = \tau_\psi
g . \label{ 2.9.2}\end{equation}  E.g.~for $f\in \CO^\times (V)$ one
has  $L^{f\psi} = L^{\psi}$ and $\tau_{f\psi }= f
(\text{div}(L/L'))\tau_\psi : {\lambda}_P (L/L')\iso    {\lambda}_P
(L^{ \psi}/ L^{\prime \psi})$;  here div$(L/L')$ was defined in 2.5.

As in 2.2 and Remark (iii) in 2.3, we denote by  $g\mapsto
\lambda_g$ the $\CL'_{\CO}(S)$-extension
$\Aut_{\CL'_{\CO}(S)}(\Lambda_P (E(\infty P)))^{\flat\prime}$ of the
group $\Aut_{\CL'_{\CO}(S)}(\Lambda_P (E(\infty P)))$; same for $E$
replaced by $E^\circ$.     Passing to $\CD et'_{P/S}$, (2.9.2)
yields then an  isomorphism
\begin{equation}
\mu^{\psi_g}\iso \lambda_{{}^\psi \! g^{-1}} \otimes\mu^\psi \iso
\mu^\psi \otimes \lambda_{g} . \label{ 2.9.3}
\end{equation}

Suppose now that $E^\circ = E$ and $\psi $ is symmetric. Then
$\tau_\psi$  is an involution of the $\CL_\CO^{\prime}(S)$-groupoid
$\Lambda_P (E(\infty P))$. Therefore, since $\mu^\psi
=\lambda_{\tau_\psi} \in \CD et'_{P/S} (E/E )=\CL'_\CO (S)$, the
composition yields a canonical identification \begin{equation}a_\psi
:  \mu^\psi \otimes\mu^\psi \iso\CO_S. \label{ 2.9.4}\end{equation}
Explicitly, the isomorphism $\mu^\psi \iso \lambda_P (\tau_\psi
(L)/L)$ identifies $a_\psi$ with the pairing $\lambda_P (\tau_\psi
(L)/L)\otimes\lambda_P (\tau_\psi (L)/L)$ $\to\CO_S$, $l_1 \otimes l_2
\mapsto \tau_\psi (l_1)l_2 = l_1 \tau_\psi (l_2 )$.

\subsection{} Let us construct on $\CE =\CE_{\text{dR}}(M)\in\CL_\CO^\Phi (X,T )$   a de Rham structure
such that the identification of (2.6.1) is horizontal.  By the
corollary in 1.12, it is uniquely defined by this property, and the
constraints from 2.8 are automatically compatible with the de Rham
structure. As in 2.1, we assume that $T$ is compatible with $M$.

As in Remark (i) in 1.2, we need to present for every scheme $S$ and
a pair of points $(D,c,\nu_P ), (D',c',\nu'_P )$ $\in \fD^\diamond (S)$
which coincide on $S_{\text{red}}$, a natural identification (notice
that $c=c'$) \begin{equation}\alpha^\varepsilon : \CE_{ (D,c,\nu_P
)}\iso \CE_{ (D',c,\nu'_P )}. \label{ 2.10.1}\end{equation} The
$\alpha^\varepsilon$ should be transitive and compatible with base
change and  factorization.

Set $P:=T^c_S \cup |D|\cup |D' |=P_{D,c}\cup P_{D',c}$. Localizing
$S$, we find an open neighborhood $V$ of $P$, $V\cap T_S =T^c_S$,
together with a datum $\nu,\nu', \kappa $, where: \newline (a)
$\nu$, $\nu'$ are trivializations of $\omega_{V/S}(D)$,
$\omega_{V/S}(D')$ which coincide on $V_{\text{red}}$ and such that
$\nu|_{P_{D,c}}=\nu_P$, $\nu|_{P_{D',c}}=\nu'_P$ (cf.~Remark (ii) in
1.1);  \newline (b) $\kappa :  M_{V\setminus P} \times
M_{V\setminus P}\to \CO_{V\setminus P}$ is   a non-degenerate
symmetric bilinear form.
\newline We construct $\alpha^\varepsilon$ explicitly using this datum.

The notation from 2.9 are in use. One can view $\kappa$ as a
non-degenerarate pairing $ M_{V\setminus P} \times \omega
M_{V\setminus P}\to \omega_{V\setminus P/S}$, which yields
$\mu_P^\kappa = \mu^{\kappa} \in \CD et'_{P/S}(M/\omega M) $. We get
\begin{equation}\mu_P^{\nabla/\kappa}:=\mu^\nabla_P \otimes
(\mu_P^{\kappa})^{\otimes -1} ,\quad \mu_P^{\kappa /\nu}:=
\mu^\kappa_P \otimes (\mu^\nu_P )^{\otimes -1}\in \CL'_{\CO}(S);
\label{ 2.10.2}\end{equation} recall that $\mu_P^\nu \in \CD
et'_{P/S}(M/\omega M) $ corresponds to the multiplication by $\nu$
identification of $M(\infty P )$ and $\omega M(\infty P)$.  Let us
rewrite (2.5.5) as an identification \begin{equation} \CE_{(D,c,\nu
)}\iso
 \mu^\nabla_P \otimes ( \mu_P^\kappa )^{\otimes -1} \otimes\mu_P^\kappa\otimes(\mu_P^\nu )^{\otimes -1}= \mu_P^{\nabla /\kappa} \otimes
 \mu_P^{\kappa /\nu}. \label{ 2.10.3}\end{equation} There is
a similar identification for $\CE_{ (D' ,c,\nu' )}$.

Notice that $ \mu_P^{\kappa /\nu}$ is the line that corresponds to
the symmetric pairing $\kappa/\nu = \nu^{-1}\kappa    : \omega
M_{V\setminus P} \times \omega M_{V\setminus P}\to
\omega_{V\setminus P/S} $. So, by (2.9.4), one has a canonical
trivialization $a_{\kappa/\nu}: \mu^{\kappa/\nu}\otimes
\mu^{\kappa/\nu}\iso \CO_S$.

Let $\beta  :  \CE_{ (D,c,\nu)}\iso \CE_{ (D',c,\nu' )}$ be an
isomorphism obtained by means of (2.10.3) from the tensor product of
$\beta_1 := \id_{\mu^{\nabla/\kappa}_{P} } $ and an identification
$\beta_2 :  \mu_P^{\kappa/\nu}\iso \mu_P^{\kappa/\nu'}$ such that
$\beta_2^{\otimes 2}= a_{\kappa/\nu'}^{-1}a_{\kappa/\nu}$ and
$\beta_2$ equals identity on $S_{\text{red}}$. We set
$\alpha^\varepsilon := \gamma\beta$, where
\begin{equation}\gamma := \exp \Res_{P/S}(\log (\nu'/\nu )\phi_\kappa
)\in\CO^\times (S).\label{ 2.10.4}\end{equation} Here $\nu'/\nu$ is
an invertible function on $V\setminus P$ that equals 1 on
$V_{\text{red}}$, so $\log (\nu'/\nu)$ is a nilpotent function on
$V\setminus P$, and $\phi_\kappa := \frac{1}{2} \nabla_M (\det
\kappa^{-1})/ \det \kappa^{-1}$ $ \in \Gamma (V\setminus P ,
\omega_{V/S} )$ where $\det \kappa^{-1} $ is the trivialization of
$\det M_{V\setminus P}^{\otimes 2}$ defined by $\kappa$.

\bigbreak\noindent\textbf{Proposition.} \emph{$\alpha^\varepsilon$
does not depend on the auxiliary choice of $\nu,\nu'$, and $\kappa
$.}\bigbreak

{\it Proof.} (a)  Let us show that $\alpha^\varepsilon$ does not
depend on  $\kappa$ for fixed $\nu$, $\nu'$. Suppose we have two
forms $\kappa$ and $\kappa_g$, so
 $\kappa_g (\cdot,\cdot )= \kappa (g\cdot,\cdot )$ for a $\kappa$-self-adjoint
$g\in \text{Aut} (M_{V\setminus P})$. Consider the corresponding
$\beta$, $\gamma$ and $\beta_g$, $\gamma_g$. Then $\beta_g /\beta $
and $ \gamma /\gamma_g$ are functions on $S$ that equal 1 on
$S_{\text{red}}$; we want to check that they are equal.

We have an $\omega$-valued symmetric bilinear form $\psi
:=\kappa/{\nu}$ on $\omega M_{V\setminus P}$. Our $g$, viewed as an
automorphism of $\omega M_{V\setminus P}$, is $\psi$-self-adjoint
and $\kappa_{g}/{\nu} =\psi_g$. Since Ad$_{\tau_\psi} (g)= g^{-1}$,
we have the isomorphism $\Ad_{\tau_\psi} :\lambda_{g^{-1}}\iso
\lambda_g$ (see 2.2). Consider the identification $\mu^{\psi_g}\iso
\lambda_{g^{ -1}}\otimes \mu^{\psi}$ of (2.9.3). The next lemma
follows directly from the definition of $a_\psi $ and $a_{\psi_g}$
as the composition in the $\flat\prime$-extension (see (2.9.4)):

\bigbreak\noindent\textbf{Lemma.}\emph{  $a_{\psi_g}$  equals the
composition $\mu^{\psi_g}\otimes\mu^{\psi_g}\iso \lambda_{g^{
-1}}\otimes \mu^{\psi}\otimes \lambda_{g^{ -1}}\otimes
\mu^{\psi}\iso \lambda_{g^{ -1}}\otimes \lambda_{g^{ -1}} \otimes
\mu^\psi \otimes \mu^\psi \iso  \CO_S$. Here the second arrow  is
the commutativity constraint, the third one is tensor product of  a
map $b_\psi : \lambda_{g^{ -1}}\otimes \lambda_{g^{ -1}}\to \CO_S$,
$\ell_1\otimes \ell_2 \mapsto \ell_1\Ad_{\tau_\psi} (\ell_2 )$,
 and $a_\psi$. \hfill$\square$ }\bigbreak

There is a similar assertion for $\psi$ replaced by $\psi'
:=\kappa/{\nu'}$. Combining them, we see that $(\beta_g /\beta )^2 =
b_\psi /b_{\psi'}$. Since $\tau_{\psi'}= h\tau_\psi $ where $h$ is
the multiplication by $\nu'/\nu$ automorphism (see (2.9.2)), one has
$ \Ad_{\tau_{\psi'}}  =\{ h ,g\}_P^\flat \Ad_{\tau_\psi} :
\lambda_{g^{-1}}\iso \lambda_g$ (see 2.2). Therefore   $(\beta_g
/\beta )^2 = (\{ h ,g\}^{\flat})^{ -1}$; by Remark (iii) in 2.3,
this equals the  Contou-Carr\`ere symbol  $ \{ \nu' /\nu , \det g
\}_P$.

Now   $\phi_{\kappa}-\phi_{\kappa_g} =\frac{1}{2}d\log  (\det g)$,
hence $\gamma /\gamma_g =\exp \Res_{P/S} (\frac{1}{2}\log (
\nu'/\nu)d\log (\det g))$ $=  \{ (\nu'/\nu)^{\frac{1}{2}},\det g \}_P
$, and we are done.

(b) It remains to show that $\alpha^\varepsilon$ does not depend on
the choice of the lifts $\nu$, $\nu'$ of $\nu_P$, $\nu'_P$.  One can
change $\nu$, $\nu'$ to $f\nu$, $f'\nu'$ where $f,f'\in\CO^\times
(V)$ are such that $f$ equals $f'$ on $V_{\text{red}}$,  $f$ equals
1 on $P_{D,c}$, $f'$ equals 1 on $P_{D',c}$.

 By (a), in the computation we are free to use $\kappa$ of our choice.
We work  $X$-locally, so we can assume that $c=1$ and pick $\kappa$
to be a non-degenerate symmetric form on $M_{X\setminus T}$. Then
$\phi_\kappa \in\omega (X,T)$ due to the compatibility of  $T$ with
$M$, see 2.1. The function $\log (f'/f )$ is regular on $V$ and
vanishes on $T^c_S$. Therefore $\log (f'/f )\phi_\kappa$ is regular
at $P$, hence its residue vanishes. Thus
  $\gamma (f'\nu',f\nu )=   \gamma (\nu',\nu )$. It remains to show that $\beta (f'\nu',f\nu )=   \beta (\nu',\nu )$.

Let $L_0$ be any $T$-lattice in $M$,  $L_0^\kappa$ be the
$\kappa$-orthogonal $T$-lattice. Since $\omega L^\kappa_0 =
\tau_{\kappa/\nu} (\omega L_0 (D))=\tau_{\kappa/{f\nu}} (\omega L_0
(D))$, one has  $\mu^{\kappa/\nu} \iso \lambda_P (\omega L^\kappa _0
/\omega L_0 (D))\buildrel{\sim}\over\leftarrow \mu^{\kappa/{f\nu}
}$. Let $e_{L_0}^\kappa : \mu^{\kappa/\nu} \iso \mu^{\kappa/{f\nu}
}$ be the composition, and $r_{L_0} : \CE_{(D,c,\nu)}\iso
\CE_{(D,c,f\nu)}$ be the tensor product of  $\id_{\mu_P^{\nabla
/\kappa }}$ and $e_{L_0}^\kappa$ (we use identifications (2.10.3)
for $\nu$ and $f\nu$). We see that $r_{L_0}$ coincides with the
isomorphism $\id_{\mu^\nabla_P} \otimes e_{L_0}^{\otimes -1}$ where
$e_{L_0} : \mu^\nu_P \iso \mu^{f\nu}_P$ was defined in 2.5. Thus $
f(D)^{-\frac{\text{rk}(M)}{2}} r_{L_0}:    \CE_{(D,c,\nu)}\iso
\CE_{(D,c,f\nu)}$ is the canonical isomorphism $r$ from (2.5.4).

We want to check that $r$ and the similar isomorphism for $f'$,
$\nu'$ identify\linebreak $\beta (f'\nu',f\nu )$ with  $\beta (\nu',\nu )$.
Indeed, $e_{L_0}^\kappa$ identifies  $a_{ \kappa/{f\nu}}$ with $
f(D)^{\text{rk}(M)}a_{ \kappa/{\nu}}$ (see Exercise in 2.9), so
$r_{L_0}$   identify $\beta (f'\nu',f\nu ) $ with $
f(D)^{\frac{\text{rk}(M)}{2}} f'(D')^{-\frac{\text{rk}(M)}{2}}
\beta (\nu' ,\nu) $. By above, this implies the assertion for $r$.
\hfill$\square$

The isomorphisms $\alpha^\varepsilon$  are  transitive (since such
are  $\alpha^\varepsilon$ with fixed  $\kappa$) and evidently
compatible with base change and factorization, so we have defined a
de Rham structure on $\CE$. The horizontality of (2.6.1) will be
checked in Example (i) of 2.11.

{\it Remark.}  Let us fix a non-degenerate symmetric bilinear form
$\kappa$ on $M_{X\setminus T}$ (like in part (b) of the proof of the
lemma). Then the above construction can be reformulated as follows.
There is a canonical isomorphism \begin{equation}\CE_{\text{dR}} (M)
\iso \CE_1 \otimes \CE_2 \otimes\CE_3  ,\label{
2.10.5}\end{equation} where $\CE_i $ are the next de Rham
factorization lines:
\newline - $\CE_{1(D,c,\nu )}:= \mu_P^{\nabla /\kappa}$
in (2.10.2); it depends only on $c$, i.e., $\CE_1$ comes from
$2^T$.\newline - $\CE_{2(D,c,\nu )}:= \mu_P^{\kappa /\nu}$ in
(2.10.2), so $\CE_2 \otimes\CE_2$ is canonically trivialized (as a
de Rham factorization line) by $a_{\kappa /\nu}$. Notice that
$\CE_2$ does not depend on the connection $\nabla_M$.
\newline -
$\CE_3 :=\CE^{-\phi_\kappa}\in \CL_{\text{dR}}^\Phi
(X,T)^{\CO\text{-triv}}$ (see  1.12), i.e., it is a de Rham
factorization line equipped with an $\CO$-trivialization $e$ with
$\nabla (e^{(1)})/e^{(1)}=-\phi_\kappa $, where $\phi_\kappa :=
\frac{1}{2}\nabla ( \det \kappa^{-1})/\det \kappa^{-1} \in \omega
(X,T)$ (see (1.12.3)). Isomorphism (2.10.5) is (2.10.3)$\otimes e$.

\subsection{}
As in 1.3, the de Rham structure on $\CE$ can be viewed as a datum
of integrable connections $\nabla^\varepsilon$ on
 $\CE_{ (D,c,\nu_P )}$ for $S$ smooth.  The next explicit construction of $\nabla^\varepsilon$  is a paraphrase of the above:

Our problem is $X$-local, so we can fix $\kappa$ as in the above
remark and
 a $T$-lattice $L$ in $M$; we can assume that $c=1$.
 Choose any $\nu$ as in Remark (ii) in 1.1.
We have identification $r_{L,\nu}: \CE_{(D,c,\nu_P )}\iso
\mu^\nabla_P (L/L(D))= \mu^\nabla_{P} (L  / \omega L^\kappa )\otimes
{\lambda}_{P} (\omega L^\kappa / L(D))$ of (2.5.6). Let $\nabla^1
=\nabla^1_{L,\kappa}$ be the ``constant" connection  on $
\mu^\nabla_{P} (L/\omega L^\kappa )$, and $\nabla^2
=\nabla^2_{L,\nu,\kappa}$ be the connection   on ${\lambda}_{P}
(\omega L^\kappa /  L(D))$ for which the pairing $a_{\kappa/\nu}:
{\lambda}_{P} (\omega L^\kappa /  L(D))\otimes {\lambda}_{P} (\omega
L^\kappa /  L(D))\to \CO_S$, $\ell_1 \otimes\ell_2 \mapsto
\tau_{\kappa/\nu}(\ell_1 )\ell_2$, of (2.9.4)  is horizontal. We get
a connection $\nabla_{\nu,\kappa} := \nabla^1 \otimes\nabla^2$  on
$\CE_{ (D,c,\nu )}$. Then
\begin{equation}
\nabla^\varepsilon=
\nabla_{\nu,\kappa} -  \theta_{\nu,\kappa} \label{
2.11.1}
\end{equation}
where $\theta_{\nu,\kappa} :=  \Res_{P/S} (d_S
\nu /\nu \otimes \phi_\kappa )\in \Omega^1_S$. Here $d_S$ is the
derivation along the fibers of $ X_S / X$, so  $d_S \nu /\nu$ is a
section of $\pi_V^* \Omega^1_S$  on $V\setminus P$.

{\it Examples.} (i) Let us compute the connection on $\CE^{(\ell )}$
(see 1.6). So let $S$ be a copy of $X\setminus T$, $P=\Delta$,
$D=\ell \Delta$. Let $t$ be a local coordinate on $X\setminus T$,
$x$ be the corresponding coordinate on $S$, $z$ be the coordinate on
$\Bbb G_m$, $\nu := z(t-x)^{-\ell}dt  $. Our $\CE^{(\ell )}$ is the
de Rham line  $ \CE_{(D , 0, \nu )}$ on $S\times \Bbb G_m$.

We take $L=M_{X\setminus T}$, so $L^\kappa =L$ and $\mu_P^\nabla
(L/\omega L^\kappa )$ is trivialized. Therefore $\CE^{(\ell
)}=\lambda_P ( \omega M /\omega M (D))$. The choice of $t$
trivializes $\omega$ and all the vector bundles $\pi_* \CO_{X\times
S}(m\Delta )/ \CO_{X\times S}(n\Delta )$, which provides an
identification $\CE^{(\ell )}\iso (\det M)^{\otimes -\ell}$. The
pairing $a_{\kappa/\nu}$ is equal (up to sign) to $(z^{-n} \det
\kappa )^{-\ell}$, $n:= \text{rk}(M)$, so
 $ \nabla_{\nu,\kappa}= \nabla_{(\det M)^{\otimes -\ell}}
+\ell\phi_\kappa +\frac{\ell n}{2}z^{-1}dz$. One has  $d_{S\times
\Bbb G_m} \nu/\nu =z^{-1}dz + \ell (t-x)^{-1}dx$. Hence
$\theta_{\nu,\kappa}= \ell\phi_\kappa$, and
\begin{equation}
\nabla^\varepsilon = \nabla_{(\det M)^{\otimes
-\ell}} +\frac{\ell n}{2}z^{-1}dz . \label{ 2.11.2}
\end{equation}

In case $\ell =1$ and $z\equiv 1$, (2.11.2) says that  (2.6.1)  is
horizontal with respect to $\nabla^\varepsilon$ and the connection
on $(\det M_{X\setminus T})^{\otimes -1}$.

(ii)  Suppose $M$ has regular singularities   at $b\in T$. Let  $t$
be a parameter at $b$. Consider $P=b$, $D=\ell b$,  and a family of
1-forms $\nu_z :=zt^{-\ell }dt$, $z\in \Bbb G_m$. Let us compute the
connection $\nabla^\varepsilon$ on the line bundle $\CE_{(b,\nu_z )}
:= \CE_{ (D,1_b ,\nu_z )}$ on $\Bbb G_m$.

Let $L$ be any  $t\partial_t$-invariant $P$-lattice  in  $M(\infty
b)$; denote by $r$ the trace of  $t\partial_t $ acting on $L/tL$,
$n:=\text{rk}(M)$. Let $\nabla_0 $ be a connection on $\CE_{(b,\nu_z
)}$ such that $r_{L,\nu_z}$ of (2.5.6) identifies it with the
``constant" connection on $z$-independent line $\mu^\nabla_b
(L/\omega L(D))$. Then  \begin{equation}\nabla^\varepsilon =\nabla_0
+(\frac{\ell n}{2} -r)z^{-1}dz .\label{ 2.11.3}\end{equation}

Indeed, consider the above construction with $\nu =\nu_z$ and
$\kappa = t^\ell  \kappa_0$ where $\kappa_0$ is a non-degenerate
symmetric bilinear form on $L$. Then $L^\kappa = L(D)$, so
$\nabla_{\nu,\kappa} =\nabla_0$. The trivialization  $\det
\kappa^{-1} $ of $\det M_{X\setminus \{ b\}}^{\otimes 2}$ has pole
of order $\ell n$ at $b$, thus the form $\phi_\kappa $ has
logarithmic singularity at $b$ with residue $r-\ell n/2$. Since $d_z
(\nu_z )/\nu_z =z^{-1}dz$, one has $\theta_{\nu,\kappa} = (r-\ell
n/2  )z^{-1}dz$, and we are done by (2.11.1).

\subsection{} Let $q:X \to Q$, $T$ be as in 1.14. Let $M$ be a
coherent $\CD_{X /Q}$-module which is $\CO_{Q}$-flat and is a vector
bundle on $X \setminus T$. We call such $M$ a {\it flat} $Q$-family
of holonomic $\CD$-modules on $(X /Q,T )$. The notion of
compatibility of $T$ and $M$ is defined as in 5.1.

Let $dR_{X/Q}(M )=\CC one (\nabla )$ be the relative de Rham
complex.
 If for some (hence every) $T$-lattices $L$ in $M$, $L_{\omega}$ in $\omega M :=\omega_{X/Q}\otimes M$  with
 $\nabla (L)\subset L_{\omega}$  the   complex
$q_{*} (dR_{X/Q}(M )/\CC one (L\buildrel{\nabla}\over\to
L_{\omega})) $ has  $\CO_{Q}$-coherent cohomology, then we call $M$
a {\it nice} $Q$-family of $\CD$-modules.

{\it Exercises.} Suppose $Q=\Spec\, \Bbb C [s]$, $X =\Spec\, \Bbb C
[t,s]$, $T$ is the divisor  $t=0$. \newline (i)  Show that  $M$
generated by a section $m$ subject to the relation $t\partial_t m =s
m$ is  nice, and that the $\CD_{X /Q}$-module $j_{T *}M =M[t^{-1}]$
is {\it not } coherent (cf.~\cite{BG}). \newline (ii) Show that $M$
generated by section $m$ subject to the relation $t^n \partial_t m
=s m$, $n>1$, is nice over the subset $s\neq 0$.

 By a straightforward relative version of the constructions of this section,
 every nice family compatible with  $T$ gives rise to a relative factorization
 line $\CE_{\text{dR}} (M )\in \CL^\Phi_{\text{dR}/Q}(X /Q,T)$ (see 1.14). The construction
 is compatible with base change. For proper $X /Q$, the $\CO_Q$-complex
 $Rq_{\text{dR}*} (M ):= Rq_* dR_{X /Q} (M )$ is perfect, and
 we  have an isomorphism of $\CO$-lines $\eta_{\text{dR}}:\CE_{\text{dR}}(M ) (X /Q)
 \iso \det Rq_{\text{dR}*} (M )$.

 Suppose  that $Q$ is smooth and
the relative connection on $M$ is extended to a flat connection (so
our nice family is isomonodromic).

\bigbreak\noindent\textbf{Proposition.}\emph{ The relative
connection on $\CE_{\text{dR}} (M)$ extends naturally to a flat
absolute connection which has local origin and is compatible with
base change and constraints from 2.8. Thus $\CE_{\text{dR}} (M)\in
\CL^\Phi_{\text{dR}}(X /Q,T )$. For proper $X /Q$,
$\eta_{\text{dR}}$ is horizontal (for the Gau\ss-Manin connection on
the target).}\bigbreak

{\it Proof.} Let $L=L((X ,T)/Q)$ be the Lie algebra  of
infinitesimal symmetries of $(X ,T )/Q$; its elements are pairs
$(\theta_{X},\theta_Q )$ where $\theta_{X}$, $\theta_Q$ are vector
fields on $X$, $Q$ such that $\theta_X$ preserves $T$ and $dq
(\theta_{X })=\theta_Q$. Our $L$ acts on $\fD^\diamond (X /Q, T
;\omega )$ and on $\CE :=\CE_{\text{dR}} (M )$ by transport of
structure. This action is compatible with constraints from 2.8 and
the relative de Rham structure.

{\it Variant.}  Let $T'\subset T$ be a component of $T$; set $L' :=
L((X \setminus T' ,T \setminus T')/Q)\supset L$.  Then $L'$ acts
naturally on  $\fD^\diamond (X /Q, \hat{T} ;\omega )$ (see 1.13
where we considered the ``vertical" part of this action). If $M
=j_{T'  *}M$, then this action lifts naturally to  $\CE_{\text{dR}}
(M )$ (as follows directly from the construction of $\CE_{\text{dR}}
(M )$). The $L'$-action extends the $L$-action (pulled back to
$\fD^\diamond (X /Q, \hat{T} ;\omega )$) and satisfies similar
compatibilities.

\bigbreak\noindent\textbf{Lemma.}\emph{ The Lie ideals $L_0 \subset
L$, $L'_0 \subset L'$ act on  $\CE$  via
$\nabla^\varepsilon$.}\bigbreak

{\it Proof of Lemma.} It suffices to check this $Q$-pointwise, so,
due to compatibility with the base change, we can assume that $Q$ is
a point. By the compatibility with the first constraint in 2.8, it
suffices to consider the cases when $M$ is supported at $T$ and
$M=j_{T*}M$. In the first situation the lemma is evident. If
$M=j_{T*}M$, then it suffices to consider the case of $L'_0$ for
$T'=T$, i.e., $L'_0 =\Theta (X\setminus T)$. The $L'_0$-action is
compatible with the de Rham structure, so we are done by the lemma
in 1.13. \hfill$\square$

We want to define the connection in a manner compatible with the
localization of $X$, so it suffices to do it in case when $X$ and
$Q$ are affine. Then $L/L_0$ is the Lie algebra of vector fields on
$Q$. Therefore, by the lemma, $\nabla^\varepsilon$ extends in a
unique manner to an absolute flat connection such that $L$ acts via
this connection.

{\it Remark.} In the situation with $T'$ the Lie algebra $L'$ acts
on $\CE$ via the connection as well (by the same lemma).

All the properties stated in the proposition, except the last one,
are evident from the construction. Let us show that
$\eta_{\text{dR}}$ is horizontal. We work $Q$-locally, so we can
assume that $Q$ is affine and $X \setminus T$ admits a section $s$.
We can enlarge $T$ to $T^+ := T \sqcup T'$, $T' :=s(Q)$. By the
first constraint from 2.8,  the assertion for $M$ reduces to that
for  $s_{*} s^* M$ and $j_{T' *} M =M (\infty T')$. The first case
is evident. In the second case the Gau\ss-Manin connection comes
from the action of the Lie algebra $L'$ (for $T^+$ and $T'$), and we
are done by the remark. \hfill$\square$

\subsection{}
{\it Compatibility of $\eta_{\text{dR}}$ with quadratic
degenerations of $X$.} We will show that $\eta_{\text{dR}}$ remains
constant (in some sense) when $X$ degenerates quadratically and $M$
stays constant outside the node. Notice that the family is {\it not}
isomonodromic (so 2.12 is not applicable). We will need the result
in \S 5; the reader can presently skip the subsection.  Consider the
next data (a), (b):

(a) A smooth proper curve $Y$, a finite subscheme $T\subset Y$, two
points $ b_+ , b_- \in (Y\setminus T) (k)$, a rational 1-form $\nu$
on $Y$  invertible off $T\cup\{ b_\pm \}$ and having poles of order
1 at $b_\pm$ with $\Res_{b_\pm}\nu =\pm 1$. Let $t_\pm$ be formal
coordinates at $b_\pm$ such that $d\log t_\pm = \pm \nu$.

(b) A $\CD$-module $N$ on $(Y,  T\cup b_\pm )$ which has regular
singularities at $b_\pm$, is the
  $*$-extension at $b_+$ and the $!$-extension at $b_-$.
 We also have a $t_{\pm}\partial_{t_\pm}$-invariant
$ b_\pm $-lattice $L$ in $M$, and an identification of the
$b_\pm$-fibers $\alpha : L_{ b_+} \iso L_{ b_-}$. Let $A_\pm \in
\End (L_{ b_\pm})$ be the action of $\pm t_\pm \partial_{t_\pm}$ on
the fibers; we ask that $\alpha A_+ =A_- \alpha$,  and that the
eigenvalues of $A_+$ (or $A_-$) and their pairwise differences
cannot be non-zero integers. Then the restriction of $L$ to the
formal neighborhoods $Y\hat{_\pm} =\text{Spf} \, k[[t_\pm ]]$ of
$b_\pm$ can be identified in a unique way with $L_{b_\pm} [[ t_\pm
]]$ so that $L_{b_\pm}\subset \Gamma ( Y\hat{_\pm},L)   $ is $t_\pm
\partial_{t_\pm}$-invariant.

Datum (a) yields a proper family of curves $X$ over $Q=\Spf\, k
[[q]]=\limright Q_n$, $Q_n = \Spec\, k[q]/q^n$, which has quadratic
degeneration at $q=0$. The 0-fiber $X_0$ is $Y$ with $b_\pm$ glued
to a single point $b_0 \in X_0 (k)$; let  $j_{b_0}$ be the embedding
$Y\setminus \{ b_+ ,b_- \}=X_0 \setminus \{ b_0 \} \hra X_0$.
Outside $b_0$ our $X$ is trivialized, i.e., $\CO_{X\setminus \{
b_0\}}= \CO_{X\setminus \{ b_0\}}[[q]]$. The formal completion of
the local ring at $b_0$ equals $k[[t_+ ,t_- ]]$ with $q=t_+ t_-$,
and the glueing comes from the embedding  $k[[t_+ ,t_- ]]\hra k((t_+
))[[q]] \times k((t_- ))[[q]]$, $t_+ \mapsto (t_+ ,q/t_- )$, $t_-
\mapsto (q/t_+ ,t_- )$. Set $\CR:= k((t_+ ))[[q]] \times k((t_-
))[[q]]/k[[t_+ ,t_- ]]$. We have a short exact sequence ($\CR$ is
viewed as a skyscraper at $b_0$) \begin{equation}0\to \CO_X \to
j_{b_0 * }\CO_{X\setminus \{ b_0 \}}= j_{b_0 * } \CO_{X_0 \setminus
\{ b_0 \}}[[q]] \to \CR\to 0, \label{ 2.13.1}\end{equation} where
the right projection assigns to $f=\Sigma f_n q^n$ the image  of
$(f_+ ,f_- )$ in $\CR$, $f_\pm =\Sigma f_n (t_\pm )q^n \in k((t_\pm
))[[q]]$ are the  expansions of $f$ at $b_\pm$.

Our family of curves has standard nodal degeneration, so we have the
dualizing line bundle $\omega_{X/Q}$. Our $\nu$ defines a rational
section $\nu_Q$ of  $\omega_{X/Q}$, which is ``constant" on
$X\setminus \{ b_0\}$ with respect to the above trivialization, and
is invertible near $b_0$.

Below for an $\CO_X$-module $F$ a relative connection on $F$ means a
morphism $\nabla : F\to \omega_{X /Q}\otimes F$ such that $\nabla
(f\phi)=d(f)\otimes \phi+ f\nabla (\phi)$, where $d$ is the
canonical differentiation $d: \CO_X \to\omega_{X/Q}$. Set $dR_{X/Q}
(F):=\CC one (\nabla)$, $Rq_{\text{dR}*}F:=Rq_* dR_{X/Q} (F)$.

Datum (b) yields an $\CO$-module $M$ on $X$ equipped with a relative
connection $\nabla$. Our $M$ is locally free over $X\setminus T_Q$.
Outside $b_0$ it is constant with respect to the above
trivialization: one has $M|_{X \setminus \{ b_0 \}}=N|_{Y\setminus
\{ b_+ ,b_- \}}[[q]]= L|_{Y\setminus \{ b_+ ,b_- \}}[[q]]  $. The
restriction $M_0$ of $M$  to $X_0$ equals $L$ with fibers
$L_{b_\pm}$ identified by $\alpha$. On the formal neighborhood of
$b_0$ our $M$ equals $M_{b_0}[[t_+ ,t_- ]]$, and the glueing comes
from the trivializations of $L$ on $Y\hat{_\pm}$ (see (b)) and
  the gluing of functions. Therefore we have a short exact sequence
\begin{equation}
0\to M\to (j_{b_0 *} L |_{Y\setminus \{ b_+ ,b_- \}} ) [[q]] \to M_{b_0}\otimes \CR
\to 0, \label{ 2.13.2}
\end{equation}
where the right projection
assigns to  $\ell =\Sigma \ell_n q^n\in(j_{b_0 *} L |_{Y\setminus \{
b_+ ,b_- \}} )$ the image in of $(\ell_+ ,\ell_- )$ $M_{b_0}\otimes
\CR$, $\ell_\pm=\Sigma \ell_n (t_\pm )q^n  \in L_{b_\pm}((t_\pm
))[[q]]=M_{b_0}((t_\pm ))[[q]]$ are the expansions of $\ell$ at
$b_\pm$ with respect to the formal trivializations of $L$ on
$Y\hat{_\pm}$. On $X \setminus \{ b_0 \}$ the relative connection
$\nabla$ comes from the $\CD$-module structure on $N$; on the formal
neighborhood of $b_0$ this is the relative connection on
$M_{b_0}[[t_+ ,t_- ]]$ with potential $At_+^{-1} dt_+ = -At_-^{-1}
dt_-$.

{\it Remarks.}  (2.13.2) is an exact sequence of $\CO_X$-modules
equipped with relative connections. The projection
$(t_+^{-1}k[t_+^{-1}] \oplus k[t_-^{-1}])[[q]] \to \CR$ is an
isomorphism. The relative connection on $M_{b_0}\otimes \CR$ is
$\nabla (m\otimes t_\pm^a q^b)= \nu (A(m)\pm a m) t_\pm^a q^b$.

Set $D=-\text{div}(\nu )-b_+ -b_-$.  Then div$(\nu_Q )= -D_Q$ does
not intersect $b_0$, so we have  an $\CO_Q$-line
$\CE_{\text{dR}}(M)_{\nu_Q}:= \CE_{\text{dR}}(M)_{(D_Q , 1_{T},\nu_Q
)}$. An immediate modification of the construction in 2.7 (to be
spelled out in the proof of the proposition below) yields an
isomorphism \begin{equation}\eta_{\text{dR}}:
\CE_{\text{dR}}(M)_{\nu_Q} \iso \det Rq_{\text{dR}*} M. \label{
2.13.3}\end{equation} Our aim is to compute it explicitly. Notice
that since our family $(X,T_Q,M, \nu_Q )$ is trivialized outside
$b_0$, one has a canonical identification
\begin{equation}\CE_{\text{dR}}(M)_{\nu_Q} \iso
\CE_{\text{dR}}(N)_{(D, 1_{T},\nu )}[[q]]. \label{
2.13.4}\end{equation}

(i) Consider an embedding $\CC one (A) [[q]] \hra dR_{X/Q}
(M_{b_0}\otimes \CR)$ whose components are $M_{b_0}[[q]] \hra
M_{b_0}\otimes \CR$, $m\otimes f(q)\mapsto m \otimes (0,f(q))$, and
$M_{b_0}[[q]] \hra \omega_{X/Q}\otimes M_{b_0}\otimes \CR$,
$m\otimes f(q)\mapsto \nu\otimes m \otimes (0,f(q))$. By the
condition on $A$, this is a quasi-isomorphism. Thus (2.13.2) yields
an isomorphism \begin{equation} \det Rq_{\text{dR}*}M \iso \det
R\Gamma_{\!\text{dR}}(Y,j_{b_\pm *}N) \otimes \det\CC one (A)[[q]].
\label{ 2.13.5 }\end{equation}

Let $L^-$ be a $b_\pm$-lattice in $N$ that equals $L$ outside $b_-$
and $t_- L$ at $b_-$. Set $C_! :=\CC one (\nabla : L^- \to \omega_Y
(\log b_+ )  L)$, $C_* :=\CC one (\nabla : L \to \omega_Y (\log
b_\pm )  L)$. Recall that $N$ is the !-extension at $b_-$ and the
$*$-extension at $b_+$, so the condition on $A$  assures that the
embeddings $C_! \hra dR (N)$ and $C_*\hra dR (j_{b_\pm *} N)$ are
quasi-isomorphisms. Since $C_* /C_! $ equals $ \CC one (A)$ (viewed
as a skyscraper at $b_-$), we see that $\CC one (dR (N)\to
dR(j_{b_\pm *} N))\iso \CC one (A)$, hence \begin{equation}\det
R\Gamma_{\!\text{dR}}(Y,N)\iso \det R\Gamma_{\!\text{dR}}(Y, j_{b_\pm
*}N) \otimes \det\CC one (A). \label{ 2.13.6}\end{equation} Combining it with
(2.13.5), we get an isomorphism \begin{equation}\det
Rq_{\text{dR}*}M \iso \det R\Gamma_{\!\text{dR}}(Y, N)[[q]]. \label{
2.13.7}\end{equation}

(ii) Consider a $\CD$-module on $\Bbb P^1 \setminus \{ 0,\infty \}$
which equals $\CO_{\Bbb P^1 \setminus \{ 0,\infty \} }\otimes
M_{b_0}$ as an $\CO$-module, $\nabla (f\otimes m)= df\otimes m +f
\otimes A(m)$. Let $\bar{N}$ be the !-extension to $\infty$ and the
$*$-extension to 0 of $N$. Consider the embeddings $t_+ :
Y\hat{_+}\hra \Bbb P^1$, $t_-^{-1}:   Y\hat{_-}\hra \Bbb P^1$ which
identify $Y\hat{_\pm}$ with the formal neighborhoods of 0 and
$\infty$.   The trivializations of $L$ on $Y\hat{_\pm}$ from (b)
identify the pull-back of $\bar{N}$ with $N|_{Y\hat{_\pm}}$. Since
the pull-back of $t^{-1}dt$ equals $\nu|_{Y\hat{_\pm}}$, we get the
identifications $\CE_{\text{dR}} (N)_{(b_+ ,\nu )}\iso
\CE_{\text{dR}}(\bar{N})_{(0,t^{-1}dt)}$, $\CE_{\text{dR}} (N)_{(b_-
,\nu )}\iso\CE_{\text{dR}}(\bar{N})_{(\infty ,t^{-1}dt)}$. Combined
with (2.7.4) (for $M$ in loc.~cit.~equal to $\bar{N}$), they produce
an isomorphism \begin{equation} \CE_{\text{dR}} (N)_{(b_+ ,\nu
)}\otimes \CE_{\text{dR}} (N)_{(b_- ,\nu )}\iso k. \label{
2.13.8}\end{equation} Since $\CE_{\text{dR}} (N)_\nu =
\CE_{\text{dR}} (N)_{(D,1_T,\nu )}\otimes \CE_{\text{dR}} (N)_{(b_+
,\nu )}\otimes \CE_{\text{dR}} (N)_{(b_- ,\nu )}$, we rewrite it as
\begin{equation}\CE_{\text{dR}} (N)_{(D,1_T,\nu )}\iso
\CE_{\text{dR}} (N)_\nu .\label{ 2.13.9}\end{equation}

\bigbreak\noindent\textbf{Proposition.}\emph{ The diagram
\begin{equation}
\begin{matrix}
 \CE_{\text{dR}}(M )_{\nu_Q} & \buildrel{\eta_{\text{dR}}}\over\lra
 &\det Rq_{\text{dR}*}M
  \\   \downarrow\,\, \,\,\,\,&&\,\,\,\,\,\,\,\,\, \downarrow    \\
 \CE_{\text{dR}}(N)_{\nu }[[q]] & \buildrel{\eta_{\text{dR}}}\over\lra
 &\det R\Gamma_{\!\text{dR}}(Y,N)[[q]],
\end{matrix} \label{ 2.13.10}
\end{equation}
where the vertical arrows are (2.13.9)$\circ$(2.13.4) and (2.13.7),
commutes. }\bigbreak

{\it Proof.} We check the assertion modulo $q^{n+1}$. Thus we
restrict our picture to $Q_n := \Spec\, R_n$, $R_n:=k[q]/q^{n+1}$;
we get a $Q_n$-curve $X_n$, $M_n =M\otimes R_n$, etc.

Let $F$ be a $\{ b_+ ,b_-\}$-lattice in $j_{b_\pm *}N$ such that
$\nabla (F)\subset \nu F=\omega_Y (\log b_\pm )\otimes F$; set $dR
(F):= \CC one (\nabla : F\to \nu F)$, $R\Gamma_{\!\text{dR}}(Y,F):=
R\Gamma (Y,dR(F))$. Then there is a canonical isomorphism
\begin{equation}\eta_{\text{dR}}(F) : \CE_{\text{dR}}(N)_{(D,1_T
,\nu )}\iso \det R\Gamma_{\!\text{dR}} (Y,F) \label{
2.13.11}\end{equation} defined as in the proposition in 2.7.
Precisely, pick any $ T\cup |D|$-lattices $E$, $E_\omega$  in $N$,
$\omega N$ such that $\nabla (E)\subset E_\omega$.
 Denote by $FE$ the $T\cup \{ b_+ ,b_-\}\cup |D|$-lattice  in $j_{b_\pm *}N$
that equals $F$ off $ T\cup |D|$ and $E$ off $b_\pm$; similarly,
$FE_\omega $ equals $E_\omega$ off $b_\pm$ and $\nu F$ off $T\cup
|D|$.  Now follow the construction from the proposition in 2.7, with
$L$, $L_\omega$, $dR (M)$ from loc.~cit.~replaced by $FE$,
$FE_\omega$ and $dR (F)$. Namely, $dR(F)$ carries a 3-step
filtration with successive quotients $FE_\omega $, $FE[1]$, $\CC
(E,E_\omega )$, and $\eta_{\text{dR}}(F)$ is the composition
$\CE_{\text{dR}}(N)_{(D,1_T ,\nu )}\iso \det \Gamma (Y, \CC
(E,E_\omega ))\otimes \lambda (FE_\omega /\nu FE)\iso \det \Gamma
(Y, \CC (E,E_\omega ))\otimes\det R\Gamma (Y, FE_\omega )\otimes
(\det R\Gamma (Y, \nu FE))^{\otimes -1}\iso \det \Gamma (Y, \CC
(E,E_\omega ))\otimes\det R\Gamma (Y, FE_\omega )\otimes \det
R\Gamma (Y,  FE[1])) \iso  \det R\Gamma (X,dR (F))$.

For example, for $F=L^-$ from (i) above one has $dR (L^- )= C_!$, so
$R\Gamma_{\!\text{dR}}(Y,L^- ))$ $=R\Gamma_{\!\text{dR}}(Y,N)$, and we get
$\eta_{\text{dR}} (L^- ): \CE_{\text{dR}}(N)_{(D,1_T ,\nu )}\iso
\det R\Gamma_{\!\text{dR}}(Y,N)$. We can also view $L^-$ as a lattice
in $N$, and compute $\eta_{\text{dR}}: \CE_{\text{dR}}(N)_\nu \iso
R\Gamma_{\!\text{dR}}(Y,N)$ using it (as in the proposition in 2.7).
Now the lemma in 2.7 implies that $\eta_{\text{dR}} (L^- )$ equals
the composition $\CE_{\text{dR}}(N)_{(D,1_T ,\nu
)}\buildrel{(2.13.9)}\over\lra \CE_{\text{dR}}(N)_\nu
\buildrel{\eta_{\text{dR}}}\over\lra R\Gamma_{\!\text{dR}}(Y,N)$.

\medskip

{\it Exercise.} If $F' \subset F$ is a sublattice with $\nabla
(F')\subset \nu F'$, then $dR (F)/dR (F')=\CC one (\nabla :F/F' \to
\nu F/\nu F')$. Thus $dR(F)$ carries a 3-step filtration with
successive quotients $dR (F')$, $\nu F/\nu F'$, $F/F'$, hence $\det
R\Gamma_{\!\text{dR}}(Y,F)\iso \det  R\Gamma_{\!\text{dR}}(Y,F')\otimes
\det \Gamma (Y, \nu F/\nu F' )\otimes \det \Gamma (Y,  F/ F'
)^{\otimes -1}$. The multiplication by $\nu$ isomorphism\linebreak $F/F' \iso
\nu F/\nu F'$ cancels the last two factors, i.e., we have $\det
R\Gamma_{\!\text{dR}}(Y,F)\iso\Gamma_{\!\text{dR}}(Y,F')$. Show that
this isomorphism equals
$\eta_{\text{dR}}(F')\eta_{\text{dR}}(F)^{-1}$.

One can repeat the above story with $Y$ replaced by $X_n$, $j_{b_\pm
*}N$ by $j_{b_0 *}N\otimes R_n$, and $E$, $E_\omega$ by $E\otimes
R_n$, $E_\omega \otimes R_n$. For a $b_0$-lattice $G$ in $j_{b_0
*}N\otimes R_n$ (i.e., an $\CO_{X_n}$-submodule, which is $
R_n$-flat and equals $j_{b_0 *}N\otimes R_n$ outside $b_0$) such
that $\nabla (G)\subset \nu G$, we get an isomorphism
\begin{equation}\eta_{\text{dR}}(G) : \CE_{\text{dR}}(N)_{(D,1_T
,\nu )}\otimes R_n\iso \det Rq_{\text{dR}*} G. \label{
2.13.12}\end{equation} For $G=M_n$, this is (2.13.3) combined with
(2.13.4). If $G$ is a ``constant" lattice, $G=F\otimes R_n$, then
$Rq_{\text{dR}*} G=R\Gamma_{\!\text{dR}}(Y,F) \otimes R_n$ and
  $\eta_{\text{dR}}(G)=\eta_{\text{dR}}(F)\otimes \id_{ R_n}$.

By above,  the proposition means that $\eta_{\text{dR}}(L^-
)\eta_{\text{dR}}(M_n )^{-1}:  \det Rq_{\text{dR}*} M  \iso \det
R\Gamma_{\!\text{dR}}(Y,N)\otimes R_n$ coincides with (2.13.7). Set
$L^{(n)} :=\CI^{n+1}L$, where  $\CI$ is the ideal of $\{ b_+ ,b_-
\}$ in $\CO_Y$. Then $ L^{(n)} $ lies in both $L^- $ and $M_n$. Set
$P_N := L^- /L^{(n)}$ and  $P_M := M_n /L^{(n)}\otimes R_n$; let
$B_N $, $B_M$ be their endomorphisms $\nu^{-1}\nabla$. One has
evident isomorphisms $\CC one (B_M)\iso
 \Gamma (X_n ,dR(M_n)/dR (L^{(n)} )\otimes R_n )$,
$\CC one (B_N )\iso \Gamma (Y, dR (L^- )/dR (  L^{(n)} ))$. Thus $
\det Rq_{\text{dR}*} M = \det R\Gamma_{\!\text{dR}}(Y, L^{(n)})\otimes
\det \CC one (B_M )$, $\det R\Gamma_{\!\text{dR}}(Y,N) \iso \det
R\Gamma_{\!\text{dR}}(Y, L^{(n)})\otimes \det \CC one (B_N )$, so
both $\eta_{\text{dR}}(L^- )\eta_{\text{dR}}(M)^{-1}$ and (2.13.7)
can be rewritten as
 isomorphisms $\det \CC one (B_M ) \iso \det  \CC one (B_N )\otimes R_n$.

Both $\det \CC one (B_M )$ and $ \det  \CC one (B_N )$ are naturally
trivialized (since\linebreak $\det \CC one (B_M )$ $=\det (P_M)\otimes \det
(P_M[1])$, etc.). By Exercise, $\eta_{\text{dR}}(L^-
)\eta_{\text{dR}}(M_n)^{-1}$ identifies  these trivializations.
Isomorphism (2.13.7) comes due to the fact that $\CC one (B_M )$ and
$\CC one (B_N ) \otimes R_n$ are naturally quasi-isomorphic: we have
the evident embeddings $i_M : M_a \otimes R_n \hra M_n /L^{(n)}$,
$i_N : M_a \hra L^- /L^{(n)}$ such that $B_M i_M =i_M A$, $B_N i_N
=i_N A$, which yield quasi-isomorphisms $\CC one (A)\otimes R_n \hra
\CC one (B_M )$, $\CC one (A)\hra \CC one (B_N )$. Therefore the
ratio of $\eta_{\text{dR}}(L^- )\eta_{\text{dR}}(M_n)^{-1}$ and
(2.13.7)  equals  the ratio of the determinants of $B_M$ and $B_N$
acting on the quotients $(M_n /L^{(n)} ) / M_a \otimes R_n$, $(L^-
/L^{(n)} )/M_a$. The first quotient is the direct sum of
components $M_a \otimes (t_+^i ,q^i t_-^{-i})\otimes R_n$ and $M_a
\otimes (q^i t_+^{-i} , t_-^i )\otimes R_n$ for $1\le i\le n$, the
second one is the direct sum of  $M_a \otimes (t_+^i ,0)$ and $ M_a
\otimes (0,t_-^i )$ for $1\le i\le n$. Both $B_M$ and $B_N$ act on
them as $A+i\,\id $ and $A-i\,\id $, so the two determinants are
equal. We are done.  \hfill$\square$

\subsection{} Suppose $k=\Bbb C$.
The definitions and constructions of this section render immediately
into the complex-analytic setting of 1.15.  Thus every triple
$(X,T,M)$, where $X$ is a smooth (not necessarily compact) complex
curve,  $T$ its finite subset, $M$ a holonomic  $\CD$-module on
$(X,T)$, yields a factorization line
$\CE_{\text{dR}}(M)\in\CL^\Phi_{\text{dR}}(X,T)$ in the
complex-analytic setting of 1.15. If $X$ and $M$ came from an
algebraic setting, then $\CE_{\text{dR}}(M)$ is an analytic
factorization line produced by the algebraic one (defined
previously).  If an algebraic family of $\CD$-modules is nice (see
2.12), then the corresponding analytic family is  nice.

We work in the analytic setting. Let $q: X \to Q$, $i: T \hra X$ be
as in 1.14; we assume that $T$ is \'etale over $Q$ (see 1.15). Let
$M$ be a flat family of $\CD$-modules on $(X /Q,T )$  which admits
locally a $T$-lattice, see 2.12. Consider the sheaf-theoretic
restriction $F:= i^\cdot dR_{X /Q}(M )$ of the relative de Rham
complex to $T$. Since $q|_{T}^\cdot \CO_Q = \CO_T$, this is a
complex of $\CO_{T}$-modules.

\bigbreak\noindent\textbf{Lemma.}\emph{ $M$ is nice if and only if
$F$ has $\CO_{T}$-coherent cohomology. }\bigbreak

{\it Proof.}  The assertion is $Q$-local, so we can assume that $T$
is a disjoint sum of several copies of $Q$.  Since $q_{*} (dR_{X
/Q}(M )/\CC one (L\buildrel{\nabla}\over\to L_{\omega}))$ is the
direct sum of pieces corresponding to the components of $T$, we are
reduced to the situation when  $X$ equals $U\times Q$, where
$U\subset \Bbb A^1$ is a coordinate disc, and $T =\{ 0\} \times Q$.

$M$ extends in a unique manner to a $\CD_{\Bbb A^1_Q/Q}$-module on
$\Bbb A^1_Q$ which is smooth outside $T$;  denote it also by $M$. So
we can assume that $X =\Bbb A^1_Q$. Set $\bar{X} :=\Bbb P^1_Q$; let
$\bar{q}: \bar{X} \to Q$ be the projection, so $X =\bar{X} \setminus
T^\infty$, $T^\infty_Q :=\{ \infty\}\times Q$.

 Let us extend $M$ to an $\CO_{\bar{X}}$-module $\bar{M}$ such that
the relative connection has logarithmic singularity at $ T^\infty$.
Such an $\bar{M}$ exists locally on $Q$. Replacing $\bar{M}$ by some
$\bar{M} (n T^\infty )$, we can assume that the eigenvalues of
$-t\partial_t$ in the fiber of $\bar{M}$ over $ T^\infty$ do not
meet $\Bbb Z_{\ge 0}$.  Let $dR_{\bar{X}/Q} (\bar{M}):=\CC one
(\bar{M}\buildrel\nabla\over\to \omega \bar{M}(T^\infty ))$ be the
relative de Rham complex of $\bar{M}$ with logarithmic singularities
at $T^\infty$. One has the usual quasi-isomorphisms
  \begin{equation} R\bar{q}_*(dR_{\bar{X}/Q}(\bar{M}))\iso q_* (dR_{X/Q}(M)) \iso   i^\cdot_T dR_{X/Q}(M) . \label{ 2.14.1}\end{equation}

Let  $\bar{L}\subset \bar{M}$, $\bar{L}_\omega \subset \omega
\bar{M}(T^\infty )$ be $\CO_{\bar{X}}$-submodules that equal $L$,
$L_\omega$ on $X$ and coincide with $\bar{M}$, $\omega
\bar{M}(T^\infty )$ outside $T$. Now $dR_{X/Q}(M)/\CC one (L\to
L_\omega) $ equals $ dR_{\bar{X}/Q  }(\bar{M})/\CC one (\bar{L}\to
\bar{L}_\omega )$, so (2.14.1) yields an exact triangle
\begin{equation}R\bar{q}_* \CC one (\bar{L}\to \bar{L}_\omega )\to
i^\cdot_T dR_{X/Q}(M)\to q_* dR_{X/Q}(M)/\CC one (L\to L_\omega).
\label{ 2.14.2}\end{equation} Its left term is $\CO_Q$-coherent, so
the other two are coherent simultaneously, q.e.d. \hfill$\square$

\section{The de Rham $\varepsilon$-lines: analytic theory}

{\it From now on we work in the analytic setting over $\Bbb C$ using
 the classical topology.}

\subsection{}\label{3.1.} Let $X$ be a smooth (not necessary compact) complex curve,
$T$ its finite subset. For a holonomic $\CD$-module $M$  we denote
by $B(M)$ the de Rham complex $dR (M)$ viewed as mere perverse $\Bbb
C$-sheaf on $X$, and set $H^\cdot_{\text{B}}(X,M):= H^\cdot
(X,B(M))$; thus one has an evident {\it period} isomorphism
$\rho: H^\cdot_{\text{B}}(X,M)\iso H^\cdot_{\text{dR}}(X,M)$. Here
is the principal result of this section:

\bigbreak\noindent\textbf{Theorem-construction.}\emph{ Let $M$, $M'$
be holonomic $\CD$-modules on $(X,T)$. Then every isomorphism $\phi
: B(M)\iso B(M')$ yields naturally an identification of the de Rham
factorization lines $\phi^\epsilon :
\CE_{\text{dR}}(M)\iso\CE_{\text{dR}}(M')$. The construction has
local origin, and is compatible with constraints from 2.8.  If $X$
is compact, then the next diagram of isomorphisms commutes:
\begin{equation}
\begin{matrix}
 \CE_{\text{dR}}(M)(X)& \buildrel{\phi^\epsilon}\over\lra & \CE_{\text{dR}}(M')(X)
 \\
 \eta_{\text{dR}} \downarrow \,\,\,\,&& \eta_{\text{dR}}\downarrow \,\,\,\,  \\
\det H^\cdot_{\text{dR}}(X,M) &
 &\det H^\cdot_{\text{dR}}(X,M')
 \\
\rho\downarrow & & \rho\downarrow
\\
 \det H^\cdot_{\text{B}} (X,M) & \buildrel{\phi}\over\lra &  \det H^\cdot_{\text{B}} (X,M').
\end{matrix} \label{ 3.1.1}
\end{equation}
}\bigbreak

The idea of the proof: By a variant of Riemann-Hilbert
correspondence, $B(M)$ amounts to the $\CD^\infty$-module
$M^\infty$. Thus what we need is to render  the story of \S 2  into
the analytic setting of $\CD^\infty$-modules, which is done using a
version of constructions from \cite{PS}, \cite{SW}.

An alternative proof of the theorem, which uses 2.13 and \S 4
instead of analytic Fredholm determinants,  is presented in 5.8.
Thus the reader can skip the rest of the section and pass directly
to \S 4.

\subsection{}\label{3.2.} {\it A digression on  $\CD^\infty$-modules and Riemann-Hilbert correspondence.} For the proofs of the next results, see
\cite{Bj} III 4, V 5.5, or \cite{Me}.

For a complex variety $X$ we denote by $\CD^\infty $ or
$\CD^\infty_X$ the sheaf of differential operators of infinite order
on $X$.  If $X$ is a curve and $U$ is an open subset  with a
coordinate function $t$, then  $\CD^\infty (U)$ consists of series
$\mathop\Sigma\limits_{n\ge 0} a_n \partial_t^n$, where $a_n$ are
holomorphic functions on $U$ such that for every $\epsilon >0$ the
series $\Sigma a_n \epsilon^{-n}n!$ converges absolutely on any
compact subset.
 $\CD^\infty_X$ is a sheaf of rings that acts on $\CO_X$ in
an evident manner;\footnote{By \cite{I}, $\CD^\infty_X$ coincides
with the sheaf of all $\Bbb C$-linear continuous endomorphisms of
$\CO_X$.} it contains $\CD_X$, and $\CD^\infty_X$ is a faithfully
flat $\CD_X$-module.

By Grothendieck and Sato, one can realize $\CD^\infty (U)$ as
$H^{\text{dim} X}_{\Delta (U)} (U\times U, \CO\boxtimes \omega )$
where $\Delta $ is the diagonal embedding. If $X$ is a curve and $U$
has no compact components, this means that
\begin{equation}\CD^\infty (U)= (\CO \boxtimes \omega )(U\times U
\setminus \Delta (U))/(\CO \boxtimes \omega )(U\times U). \label{
3.2.1}\end{equation} Here $k(x,y)\in (\CO \boxtimes \omega )(U\times
U \setminus \Delta (U))$ acts on $\CO (U)$ as $f\mapsto k(f)$,
$k(f)(x):= \Res_{y=x}k(x,y)f(y)$.

For a (left) $\CD$-module $M$ set $M^\infty :=  \CD^\infty
\mathop\otimes\limits_{\CD}M$. The embedding   $dR (M)\hra
dR(M^\infty )$ is a quasi-isomorphism, hence
$H^\cdot_{\text{dR}}(X,M)\iso H^\cdot_{\text{dR}}(X,M^\infty )$. If
$M$ is smooth, then $M\iso M^\infty$.

For a $\CD^\infty$-module $N$ a {\it $\CD$-structure} on $N$ is a
$\CD$-module $M$ together with a $\CD^\infty$-isomorphism $M^\infty
\iso N$. Our $N$ is said to be holonomic if it admits a
$\CD$-structure with holonomic $M$.

For a holonomic $\CD$-module $M$ the de Rham complex $dR(M)$ is a
perverse $\Bbb C$-sheaf, which we denote, as above, by $B(M)$; same
for a holonomic $\CD^\infty$-module. Therefore $B(M)=B(M^\infty )$.
The functor $B$ is an equivalence between the category of holonomic
$\CD^\infty$-modules and  that of perverse $\Bbb C$-sheaves. The
inverse functor assigns to a perverse sheaf $F$ the
$\CD^\infty$-module\footnote{Here  $\CD^\infty$ acts via the
$\CO_X$-factor, and $F\otimes^! G :=R \Delta^! F\boxtimes G$. }  $
\CO_X {\mathop\otimes\limits_{\Bbb C}}^! F[\dim X]$. Thus for a
holonomic $M$ the $\CD^\infty$-module $M^\infty$ carries the same
information as $B(M)$.

The  functor $M\mapsto M^\infty$  yields an equivalence between the
category of holonomic $\CD$-modules with regular singularities and
that of holonomic $\CD^\infty$-modules. Its inverse assigns to a
holonomic $\CD^\infty$-module $N$ its maximal $\CD$-submodule
$N^{\text{rs}}$ with regular singularities, so one has
\begin{equation}(N^{\text{rs}})^{\infty}=N. \label{ 3.2.2}\end{equation} Therefore
every holonomic $\CD^\infty$-module admits a unique $\CD$-structure
with regular singularities.

{\it Exercises.} Let $U$ be a coordinate disc, $t$ be the
coordinate, $j$ be the embedding $U^o :=U\setminus \{ 0\}\hra U$.

(i) Recall a description of indecomposable $\CD_U$-modules which are
smooth of rank $n$ on $U^o$ and have regular singularity at 0. For
$s\in \Bbb C$ denote by $M_{s,n}$  a $\CD$-module whose sections are
collections of functions $(f_i )=(f_1 ,\ldots ,f_n )$ having
meromorphic singularity at 0, and $\nabla_{\partial_t }((f_i ))=
(\partial_t (f_i)+ s f_i +f_{i-1})$.\footnote{ $M_{s,n}$ depends
only on $s$ modulo $\Bbb Z$-translation: one has $M_{s,n}\iso
M_{s-1,n}$, $(f_i)\mapsto (tf_i )$.} Let $M_{0,n}^1$ be a
$\CD$-submodule of $M_{0,n}$ formed by $(f_i )$ with $f_1$ regular
at 0. Consider an embedding $\CO_U \hra M_{0,n+1}$, $f \mapsto
(0,\ldots, 0, f)$; set $M_{0,n}^2 :=M_{0,n+1}/\CO$, $M_{0,n}^3
:=M^1_{0,n+1}/\CO$. E.g., $M^1_{0,1}=\CO_U$, and $M_{0,0}^2 =\delta$
(the $\delta$-function $\CD$-module). Then any indecomposable
$\CD_U$-module $M$ as above is isomorphic to either some $M_{s,n}$,
or one of $M_{0,n}^a$,  $a=1,2,3$.

Show that the corresponding $\CD^\infty$-module $M^\infty$ has the
same explicit description  with ``meromorphic singularity" replaced
by ``arbitrary singularity".

(ii) For $n>0$ let $E_{(n)}$ be a $\CD$-module of rank 1  generated
by $\exp (t^{-n})$, i.e., $E_{(n)}$ is generated by a section $e$
subject to the (only) relation $t^{n+1}\partial_t (e)=-ne$. Show
that there is an isomorphism of $\CD^\infty$-modules
\begin{equation} E_{(n)}^\infty \iso M_{0,1}^{2\infty}\oplus
(\delta^\infty )^{n-1} . \label{ 3.2.3}\end{equation} Here is an
explicit formula for (3.2.3). Let $g(z)$, $h_1 (z),\ldots
,h_{n-1}(z)$
 be  entire functions such that $(\partial_z -n z^{n-1})g(z)= z^{-2}(\exp (z^n )-1-z^n )$ and
$(\partial_z -n z^{n-1})h_i (z)= z^{i-1}$. Then (3.2.3) assigns $e$
a vector whose $(M_{0,1}^{2})^\infty$-component  is $(\exp (t^{-n}),
g(t^{-1}))$ and the  $\delta^\infty$-components are $h_i (t^{-1})\in
(M^2_{0,0})^\infty =\delta^\infty$.

\subsection{}\label{3.3.} {\it A digression on Fredholm determinants} (cf.~\cite{PS} 6.6).  Recall that a
{\it Fr\'echet space} is a  complete, metrizable,  locally convex
topological $\Bbb C$-vector space. The category $\CF$ of those is a
quasi-abelian (hence exact) Karoubian  $\Bbb C$-category. A morphism
$\phi : F\to F'$ is said to be  Fredholm if it is Fredholm as a
morphism of abstract vector spaces, i.e., if $\Ker\, \phi$ and
Coker$\,\phi$ have finite dimension. Then $\phi$ is a split
morphism, i.e., $\Ker\phi$, Im$\phi$ are direct summands of,
respectively, $F$ and $F'$,  and $F/\Ker \phi \iso \text{Im} \phi$.
Denote by $\fF \subset \CF$  the subcategory of  Fredholm morphisms
$\fF (F,F')\subset \Hom (F,F')$.

A Fredholm $\phi$ yields the determinant line $\lambda_{\phi} :=
\det (\Coker\, \phi)\otimes \det^{\otimes -1}(\Ker\, \phi )\in
\CL:=\CL_{\Bbb C}$ (see 1.2). Sometimes we denote $\lambda_\phi$ by
$\lambda (F'\buildrel{\phi}\over\to F)$ or, if $F=F'$, by $\lambda
(F)_\phi $. If $\phi$ is invertible, then $\lambda_\phi$ has an
evident trivialization; denote it by  $\det (\phi )\in
\lambda_\phi$.

 For any Fredholm $\phi$ one can find finite-dimensional
$F_0  \subset F$, $F'_0  \subset F'$  such that $\phi (F'_0 )\subset
F_0$ and the induced map $ F'/F'_0 \to F/F_0$ is an isomorphism
(equivalently, $F_0 +\phi (F')=F$, $F'_0 = \phi^{-1}(F_0 )$). Then
the exact sequence $0\to\Ker\phi \to F'_0 \to F_0 \to\Coker \phi \to
0$ yields a natural isomorphism \begin{equation}\lambda_\phi \iso
\det (F_0 )\otimes{\det}^{\otimes -1} (F'_0 ).\label{
3.3.1}\end{equation} If $\phi$ is invertible, then (3.3.1)
identifies $\det (\phi )\in\lambda_\phi$ with the usual determinant
of $\phi|_{F'_0} : F'_0 \iso F_0$ in $\Hom (\det F'_0 ,\det (F_0 ))=
\det (F_0 )\otimes{\det}^{\otimes -1} (F'_0 )$.

For  Fredholm $F''\buildrel{\phi'}\over\to F'\buildrel{\phi}\over\to
F$ there is  a canonical ``composition" isomorphism \begin{equation}
\lambda_\phi \otimes\lambda_{\phi'} \iso \lambda_{\phi\phi'}, \quad
a\otimes b\mapsto ab, \label{ 3.3.2}\end{equation} which satisfies
the associativity property. Therefore $\phi \mapsto \lambda_\phi$ is
a (central) {\it $\CL$-extension} $\fF^{\flat}$ of $\fF$ (see
e.g.~\cite{BBE}, Appendix to \S 1,   for terminology). To construct
(3.3.2), choose $F_0$, $F'_0$, $F''_0$ for $\phi$, $\phi'$ as above;
then (3.3.1) identifies the composition with an evident map $\det
(F_0 ) \otimes \det^{\otimes -1} (F'_0 )\otimes\det (F'_0
)\otimes\det^{\otimes -1} (F''_0 )\iso \det (F_0 ) \otimes
\det^{\otimes -1} (F''_0 )$. For invertible $\phi$, $\phi'$ one has
$\det (\phi )\det (\phi' )=\det (\phi \phi' )$.

Suppose   $F$, $F'$ are equipped with finite split filtrations
$F_\cdot$, $F'_\cdot$,  $\phi : F'\to F$  preserves the filtrations,
and  $\gr \phi : \gr F' \to \gr F$ is Fredholm. Then $\phi$ is
Fredholm, and there is a canonical isomorphism
\begin{equation}\lambda_{ \phi}\iso \otimes \lambda_{\gr_i \phi} .
\label{ 3.3.3}\end{equation} The identification is transitive with
respect to refinement of the filtration. If $\gr\phi$ is invertible,
it identifies $\det (\gr \phi )=\otimes\det (\gr_i \phi )$ with
$\det (\phi )$. For example, for a finite collection of Fredholm
morphisms $\{ \phi_\alpha \}$, every linear ordering of indices
$\alpha$ produces a filtration, hence an isomorphism
$\lambda_{\oplus \phi_\alpha }\iso \otimes\lambda_{\phi_\alpha}$; it
does not depend on the ordering.

Let $\CI^{\text{fin}}\subset \CI^{\text{tr}}\subset\CI^{\text{com}}$
be the two-sided ideals of finite  rank, nuclear, and compact
morphisms in $\CF$. We have the quotient categories $\CF /\CI^?$:
their objects are Fr\'echet spaces, and morphisms $\Hom_{/\CI^?}
(F,F')$ equal $ \Hom (F, F')/ \CI^? (F,F')$. A morphism $\phi$ is
Fredholm if and only if $\phi$ is invertible in either $\CF/\CI^?$.
Therefore the groupoids  Isom$(\CF/\CI^?)$ of isomorphisms in $\CF
/\CI^?$ are quotients of $\fF$ modulo the {\it $\CI^?$-equivalence}
relation  $\phi -\phi' \in\CI^?$.

{\it Exercise.} Let $G^? (F)\subset \Aut (F)$ be the (normal)
subgroup of automorphisms $\psi$ of $F$ that are
$\CI^{?}$-equivalent to $\id_F$. The next sequence is exact:
\begin{equation} 1\to G^{\text{?}}(F)/G^{\text{fin}}(F) \to \Aut_{/
\CI^{\text{fin}}}(F)\to \Aut_{/ \CI^{\text{?}}}(F)\to 1 \label{
3.3.4}\end{equation}

\bigbreak\noindent\textbf{Proposition.}\emph{ $\fF^{\flat}$ descends
naturally to an $\CL$-extension Isom$^\flat (\CF  /\CI^{\text{tr}})$
of the groupoid Isom$(\CF  /\CI^{\text{tr}})$. }\bigbreak

{\it Proof.}  We first descend $\fF^{\flat}$ to Isom$(\CF
/\CI^{\text{fin}})$, and then to Isom$(\CF  /\CI^{\text{tr}})$.

(i) To  descend $\fF^{\flat}$ to  Isom$(\CF  /\CI^{\text{fin}})$,
means to define for every $\CI^{\text{fin}}$-equivalent $\phi, \psi
\in \fF (F, F')$ a natural identification $\tau =\tau_{\phi,\psi}
:\lambda_\phi \iso \lambda_{\psi}$ which satisfies the transitivity
property and is compatible with composition.

The condition on $\phi$, $\psi$ means that we can find
finite-dimensional  $F_0 \subset F$, $F'_0 \subset F'$ such   that
$\phi (F'_0 ), (\phi -\psi) (F')\subset F_0$, and the map $ F'/F'_0
\to F/F_0$ induced by $\phi$ (or $\psi$) is an isomorphism. Then
$\tau$ is the composition $\lambda_\phi \iso \det (F_0
)\otimes\det^{\otimes -1} (F'_0 )\iso \lambda_{\psi}$ of
isomorphisms (3.3.1) for $\phi$, $\psi$. The construction  does not
depend on the choice of auxiliary datum, and satisfies the necessary
compatibilities.

(ii)  Recall that for  $\psi \in \End F$ that is
$\CI^{\text{tr}}$-equivalent to $\id_F$, its {\it Fredholm
determinant} $\det_\fF (\psi )\in\Bbb C$ is defined (see
e.g.~\cite{Gr2})  as the sum of a rapidly converging series
\begin{equation}{\det}_\fF (\psi ):=\mathop\Sigma\limits_{k\ge 0}
\tr \Lambda^k ( \psi - \id_F ), \label{ 3.3.5}\end{equation} where
$\Lambda^k ( \psi - \id_F )$ is the $k$th exterior power of $\psi -
\id_F$. If $\psi - \id_F $ is of finite rank, then the sum is
finite, and $\det_\fF (\psi )$ is the usual
determinant.\footnote{i.e.,  $\det (\psi|_{F_0})$ where $F_0$ is any
finite-dimensional subspace containing the image of $\psi - \id_F
$.}

The central $\Bbb C^\times$-extension $\Aut^\flat (F)$ of $\Aut (F)$
is trivialized by the section $\psi \mapsto \det (\psi )$. The
Fredholm determinant is multiplicative and invariant with respect to
the adjoint action of $\Aut (F)$. We get a trivialization $\psi
\mapsto \tau^{\text{an}} (\psi ):= \det_\fF^{-1} (\psi ) \det (\psi
)$ of the $\Bbb C^\times$-extension $G^{\text{tr}}(F)^\flat$ which
is invariant for the adjoint $\Aut (F)$-action.

Since for $\psi \in G^{\text{fin}}(F)$ one has $\tau^{\text{an}}
(\psi )= \tau_{\id_F ,\psi}\in \lambda_\psi$, our $\tau^{\text{an}}$
can be viewed as a trivialization  of the extension $ \Aut^\flat_{/
\CI^{\text{fin}}}(F)$ over the normal subgroup
$G^{\text{tr}}(F)/G^{\text{fin}}(F)$. It is invariant with respect
to the adjoint action of $\Aut_{/ \CI^{\text{fin}}}(F)$. Thus, by
(3.3.4),  $\tau^{\text{an}}$ defines a descent of $ \Aut^\flat_{/
\CI^{\text{fin}}}(F)$ to an $\CL$-extension   $\Aut_{/
\CI^{\text{tr}}}^\flat (F)$ of $\Aut_{/ \CI^{\text{tr}}}(F)$.

More generally, for every $F, F'\in\CF$, the set {Isom}$_{/
\CI^{\text{fin}}} (F,F')$ is  a\linebreak
$(G^{\text{tr}}(F')/G^{\text{fin}}(F'),
G^{\text{tr}}(F)/G^{\text{fin}}(F))$-bitorsor over Isom$_{/
\CI^{\text{tr}}} (F,F')$, and we define the $\CL$-extension Isom$_{/
\CI^{\text{tr}}} (F,F')^\flat $ as the quotient of Isom$_{ /
\CI^{\text{fin}}}^\flat (F,F')$ by the  $\tau^{\text{an}}$-lifting of
either  $G^{\text{tr}}(F)/G^{\text{fin}}(F)$- or 
$G^{\text{tr}}(F')/G^{\text{fin}}(F')$-action. \hfill$\square$

\medskip

{\it Remark.} The above constructions  are compatible with
constraint (3.3.3).

\subsection{}\label{3.4.}
For  a topological space $X$ whose topology has countable base, a
{\it Fr\'echet sheaf} on $X$  means a sheaf of Fr\'echet vector
spaces. A {\it Fr\'echet algebra} $\CA$  is a sheaf of topological
algebras which is a Fr\'echet sheaf; a {\it Fr\'echet $\CA$-module}
is a Fr\'echet sheaf equipped with a continuous (left) $\CA$-action.

The problem of finding Fr\'echet structures on a given $\CA$-module
$M$ is delicate. Here is a simple uniqueness assertion. Suppose that
$M$ satisfies the next condition:  those open subsets $U$ of $X$
that $M(U)$ is a finitely generated $\CA (U)$-module form a base of
the topology of $X$.

\bigbreak\noindent\textbf{Lemma.}\emph{   Every morphism of
$\CA$-modules $\phi: M\to N$ is continuous with respect to any
Fr\'echet structures on $M$, $N$. Thus $M$ admits at most one
Fr\'echet structure. }\bigbreak

{\it Proof.}  It suffices to check that the maps $\phi_U : M(U)\to
N(U)$ are continuous for all $U$ as above. Thus  there is a
surjective  $\CA (U)$-linear map $\pi_U : \CA (U)^n
\twoheadrightarrow M(U)$. The maps $\pi_U$ and $\phi_U \pi_U$ are
evidently continuous. Since $\CA (U)^n /\Ker (\pi_U ) \to M(U)$ is a
continuous algebraic isomorphism of Fr\'echet spaces, it is a
homeomorphism, and we are done.  \hfill$\square$

{\it Example.} Every locally free $\CA$-module of finite rank is  a
Fr\'echet $\CA$-module.

From now on our $X$ is a complex curve. The two basic examples of
Fr\'echet algebras on $X$ are $\CO_X$ and $\CD^\infty_X$. For an
open $U\subset X$  the topology on the space of holomorphic
functions $\CO (U)$  is  that of uniform convergence on compact
subsets of $U$. If $t$ is a coordinate function on $U$, then the
topology on $\CD^\infty (U)$ is given by a collection of semi-norms
$|| \Sigma a_n \partial_t^n ||_{K\epsilon}:=
\mathop{\text{max}}\limits_{x\in K} \Sigma |a_n
(x)|\epsilon^{-n}n!$; here   $K$ is any compact subset of $U$ and
$\epsilon$ is any small positive real number. Equivalently, one can
use (3.2.1): then $(\CO \boxtimes \omega )(U\times U)$ is a closed
subspace of $(\CO \boxtimes \omega )(U\times U \setminus \Delta
(U))$, and the topology on $\CD^\infty (U)$ is the quotient one.

\bigbreak\noindent\textbf{Proposition.}\emph{ Any  holonomic
$\CD^\infty$-module $N$  on $X$ admits a unique structure of a
Fr\'echet $\CD^\infty_X$-module. }\bigbreak

{\it Proof.}  Uniqueness: As follows easily from Exercise (i) in
3.2,  $N$ satisfies the condition of the previous lemma. Existence:
The problem is local, so it suffices to define some Fr\'echet
structure compatible with the $\CD^\infty (U)$-action on $N(U)$,
where $U$ is a disc and $N$ is smooth outside the center 0 of $U$.
Then $N=M^\infty$ where $M$ is a $\CD$-module with regular
singularities; we can assume that $M$ is indecomposable. If $M\simeq
M_{s,n}$ (see Exercise (i) in 3.2), then, by loc.~cit., $M^\infty (U
)\iso M^\infty (U^o )\simeq \CO (U^o )^n$, and we equip it with the
topology of $\CO (U^o )^n$. Otherwise $M^\infty (U )$ is a
subquotient of some $M_{0,n}^\infty (U)$, and we equip it with the
corresponding Fr\'echet structure.
 \hfill$\square$

{\it Question.} Can one find a less ad hoc proof  (that would not
use (3.2.2))? Is the assertion of the proposition remains true for
all perfect $\CD^\infty$-modules (or perfect $\CD^\infty$-complexes)
on $X$ of any dimension?\footnote{For a perfect $\CD^\infty$-complex
$N$,  \cite{PSch} define a natural ind-Banach structure on its
complex of solutions $R\CH om_{\CD^\infty}(N,\CO_X)$. It is not
clear if  this result helps to see the topology on $N$. }

 \subsection{}\label{3.5.} Let $K\subset X$ be a compact subset which
 does not contain a connected component of $X$;
denote by $j_K$ the embedding $X\setminus K \hra X$. Let  $E$ be  a
Fr\'echet $\CO_X$-module. Suppose that for some open neighborhood
$U$ of $K$, $E|_{U\setminus K}$ is a locally free $\CO_{U\setminus
K}$-module of finite rank. A {\it $K$-lattice} in $E$ is an
$\CO_X$-module $L$, which is locally free on $U$, together with an
$\CO_X$-linear morphism $L\to E$  such that $L|_{X\setminus K}\iso
E|_{X\setminus K}$. Then $L$ is a Fr\'echet $\CO_X$-module, and
$L\to E$ is a continuous morphism. Set $\Gamma (E/L):= H^0 R\Gamma
(X,\CC one (L\to E))=H^0 R\Gamma (U,\CC one (L\to E))$.

Shrinking $U$ if needed, we can assume  that the closure $\bar{U}$
of $U$ is compact with smooth  boundary $\partial \bar{U}$. We
denote by $\partial U$ a  contour  in $U\setminus K$  homologous to
the boundary of $ \bar{U}$ in $\bar{U}\setminus K$.

{\it Remarks.} (i) For all our needs it suffices to consider the
situation when $U$ is a disjoint union of several discs.

(ii) If Int$(K)\neq \emptyset$, then the morphism $L\to E$ need not
be injective. The map $L(U)\to E(U)$ is injective though, i.e.,
$H^{-1} R\Gamma (U,\CC one (L\to E))=0$.

{\it Example.} If $L$ is any locally free $\CO_X$-module of finite
rank, then $L$ is a $K$-lattice in $j_{K\cdot} L:= j_{K\cdot
}(L|_{X\setminus K})$.

\bigbreak\noindent\textbf{Proposition.}\emph{ (i) $L(U)$ is a direct
summand of the Fr\'echet space $E(U)$.  \\ (ii) If $H^1 (U,L)=0$
(which happens, e.g., if none of the connected components of $U$ is
compact), then $E(U)/L(U)\iso \Gamma (E/L).$
\\
(iii)  Let $(E',L')$ be a similar pair, and $\phi : E'\to E$ be any
morphism of Fr\'echet sheaves. Then the map $ E' (U)\to E(U)$,
viewed as a morphism in $\CF /\CI^{\text{tr}}$, sends the subobject
$L' (U)$  to $L (U)$. }\bigbreak

{\it Proof.} (i) We want to construct a left inverse to the morphism
of Fr\'echet spaces $L(U)\to E(U)$. It suffices to define a left
inverse to the composition $L(U)\to E(U)\to E(U\setminus
K)=L(U\setminus K)$, i.e., to the restriction map $L(U)\to
L(U\setminus K)$.

We can assume that $U$ is connected and non compact.\footnote{If $U$
is  compact, then $L(U)$ is finite dimensional, and the assertion
follows from the Hahn-Banach theorem.} Then one can find a Cauchy
kernel  on $U\times U$, which is a section $\kappa$  of $L\boxtimes
\omega L^* (\Delta )$ with residue at the diagonal equal to
$-\id_E$. The promised left inverse is  $f\mapsto \kappa (f)$,
$\kappa (f)(x)= \mathop\int\limits_{\partial U} \kappa (x,y)f(y)$.

(ii) Follows from the exact cohomology sequence.

(iii) We want to check that the composition $L' (U)\to
 E' (U)\to E(U)\to E (U)/   L (U) $ is nuclear. We can assume that $U$ is connected and non compact.
Choose an open $V\supset K$ whose closure $\bar{V}$ is compact and
lies in $U$. Our map equals the composition  $L' (U)\to L' (V)\to E'
(V)\to E(V)\to E(V)/L(V)\buildrel\sim\over\leftarrow E(U)/L(U)$ (for
$\buildrel\sim\over\leftarrow $, see  (ii)). We are done, since the
first arrow is nuclear (see e.g.~\cite{Gr1}). \hfill$\square$

\bigbreak\noindent\textbf{Corollary.}\emph{ (a) The isomorphism from
(ii) yields a natural Fr\'echet space structure on $\Gamma (E/L)$,
which does not depend on the auxiliary choice of $U$.
\\
(b) Every  $\phi$ as in (iii) yields naturally a morphism
$\phi_{E'/L',E/L}:\Gamma (E'/L')\to \Gamma (E/L)$ in $\CF
/\CI^{\text{tr}}$. In particular, the spaces  $\Gamma (E/L)$ for all
$K$-lattices $L$ in $E$ are canonically identified  as objects of
$\CF /\CI^{\text{tr}}$. \hfill$\square$ }\bigbreak

{\it Proof.} (a) follows since, for $U'\subset U$ as in (ii), the
restriction map $E(U)/L(U)\to E(U')/L(U')$ is a continuous algebraic
isomorphism, hence a homeomorphism, and $U$'s form a directed set.
The first assertion in (b) follows from (iii); for the second one,
consider $\phi =\id_E$. \hfill$\square$

\subsection{}\label{3.6.}
The set $\Lambda_K (E)$  of   $K$-lattices in $E$ has natural
structure of  an $\CL$-groupoid (see 2.2). Namely,
 by the corollary in 3.5,
for every $L,L' \in\Lambda_K (E)$  one has a  canonical
identification $\id_{E/L,E/L'} : \Gamma (E/L )\iso \Gamma (E/L')$ in
$\CF /\CI^{\text{tr}}$. We set
\begin{equation}\lambda_K (L /L' ):= \lambda_{\id_{E/L,E/L'}}. \label{ 3.6.1}\end{equation} The composition
\begin{equation}
\lambda_K (L /L' )\otimes \lambda_K (L' /L'')\iso \lambda_K (L /L''
) \label{ 3.6.2}\end{equation} comes from (3.3.2).  \hfill$\square$

Suppose  $\Lambda_K (E)$ is non-empty. We get an $\CL$-torsor $\CD
et_K (E):= \Hom_\CL (\Lambda_K (E),\CL )$ of {\it
 determinant theories on $E$ at $K$}. For $E_1$, $E_2$ we get an
 $\CL$-torsor $\CD et_K (E_1 /E_2 ):= \CD et_K (E_1 )\otimes \CD et_K (E_2 )^{\otimes -1}$ of
 {\it relative determinant theories on $E_1 / E_2 $ at $K$}.

 If $K=\sqcup K_\alpha$,
then a $K$-lattice $L$ amounts to a collection of
$K_\alpha$-lattices $L_\alpha$, and one has an evident canonical
isomorphism \begin{equation}\otimes \lambda_{K_\alpha} (L_\alpha
/L'_\alpha )\iso \lambda_K (L/L') \label{ 3.6.3}\end{equation}
compatible with composition isomorphisms (3.6.2). Thus one has a
canonical identification $\otimes \CD et_{K_\alpha }(E)\iso  \CD
et_{K }(E)$, $(\otimes\lambda_\alpha )(L)=\otimes \lambda_\alpha
(L_\alpha )$.

Below we fix a neighborhood $U$ of $K$ as in 3.5; we assume that it
has no compact components, so for every $K$-lattice $L$ the $\CO
(U)$-module $L(U)$ is free.

$\Lambda_K (E)$ carries a natural topology of compact convergence on
$U\setminus K$. Namely, to define a neighborhood of $L$ we pick an
$\CO (U)$-base $\{ \ell_i  \}$ of $L(U)$, a compact $C\subset
U\setminus K$ and a number $\epsilon >0$. The neighborhood is formed
by those $L'$ which admit a base $\{ \ell'_i \}$ of $L'(U)$, $
\ell'_i =\Sigma a_{ij}\ell_i$, with $(a_{ij})$  $\epsilon$-close to
the unit matrix on $C$.

The $\CL$-groupoid structure on $\Lambda_K (E)$ is continuous, i.e.,
$\lambda_K$ form naturally a line bundle on $\Lambda_K (E)\times
\Lambda_K (E)$, and the composition is continuous (cf.~\cite{PS}
6.3, 7.7). Namely, if $P$ is a closed linear subspace of $E(U)$,
then the subset $\Lambda^{P}_K (E) $ of those $K$-lattices   $L$
that $L(U)$ is complementary to $P$, i.e., $P\iso \Gamma (E/L)$, is
open in $\Lambda_K (E)$. For $L,L'\in \Lambda^{P}_K (E) $ the
morphism $\id_{E/L,E/L'}$ in $\CF /\CI^{\text{tr}}$ is represented
by $\id_P$, hence the restriction of $\lambda_K$ to $\Lambda^{P}_K
(E)\times \Lambda^{P}_K (E)$ is trivialized by the section $\delta^P
:=\det (\id_P )$. The topology on $\lambda_K$ is uniquely determined
by the condition that all these local trivializations are
continuous, and the composition is continuous.

    {\it Remarks.} (i) If $K'$ is a compact such that $K\subset K'\subset U$, then every $K$-lattice $L$ in $E$ is a $K'$-lattice, and
$\lambda_{K'}(L/L')=\lambda_K (L/L')$. Hence $\CD et_K (E) =\CD
et_{K'}(E)$.

(ii)  By (iii) of the proposition in 3.5, the subobjects $L(U)$,
$L'(U)$ of $E(U)$ coincide if viewed in $\CF /\CI^{\text{tr}}$. Let
$\phi_U :L(U)\iso L' (U)$ be the corresponding isomorphism in $\CF
/\CI^{\text{tr}}$. Then there is a canonical identification
\begin{equation} \lambda_K (L/L')\iso  \lambda_{\phi_U}^{\otimes -1}
\label{ 3.6.4}\end{equation} compatible with the composition maps.
Indeed, by (3.3.3) and Remark at the end of 3.3 applied to the
filtrations $L^{(\prime )} \subset E(U)$, one has
$\lambda_{\id_{E/L, E/L'}} \otimes \lambda_{\phi_U}\iso
\lambda_{\id_{E(U)}}=\Bbb C$.

(iii) Every $K$-lattice $L$ in $E$ can be viewed as a $K$-lattice in
$j_{K\cdot}E:= j_{K\cdot}j_K^\cdot E$ (via $L\to E\to j_{K\cdot
}E$). Then (3.6.4) shows that $\lambda_K (L/L')$ does not depend on
whether we consider $L$, $L'$ as $K$-lattices in $E$ or in
$j_{K\cdot }E$. Thus $\CD et_K (E)=\CD et_K (j_{K\cdot}E)$.

(iv) Let $S$ be an analytic space, $L_S$ be an $S$-family of
$K$-lattices in $E$.\footnote{I.e., $L_S$ is an $\CO_{X_S}$-module
together with a morphism $L\to E_S$ such that $L_S$ is locally free
on $U_S$, and $L\to E_S$ is an isomorphism off $K_S$.} Then the
pull-back of $\lambda_K$ to $S\times S$ is naturally a holomorphic
line bundle, so that the pull-back to $S$ of any local
trivialization $\delta^P$ is holomorphic.

(v) If $L$, $L'$ are meromorphically equivalent, then $\lambda_K (L
/L' )$ coincides with the relative determinant line from 2.3 (where
$P$ is a finite subset in $K$ such that $L$ equals $L'$ off $P$).
Indeed, in view of (3.6.2), it suffices to identify the lines in
case $L\supset L'$, where the identification comes from
$L(U)/L'(U)\iso \Gamma (U,L/L')$. If $L$, $L'$ vary holomorphically
as in (iv), then this identification is holomorphic.

(vi) For every $f\in \CO^\times (U\setminus K)$ the lines $\lambda_K
(fL/L)$ for  all $L\in\Lambda_K (j_{K\cdot}E)$ are canonically
identified. Namely, one defines the isomorphism $\lambda_K
(fL'/L')\iso \lambda_K (fL/L)$ as $\lambda_K (fL'/L')\iso \lambda_K
(fL'/fL)\otimes \lambda_K (fL/L)\otimes \lambda_K (L'/L)^{\otimes
-1} \iso \lambda_K (fL/L)$ where the first arrow is inverse to the
composition, and the second comes from the multiplication by $f$
identification $\lambda_K (L'/L)\iso \lambda_K (fL'/fL)$.

(vii) Let $g\in\CO^\times (U)$ be an invertible function. The
multiplication by $g$ automorphism of $E|_U$ preserves every
$K$-lattice. Let $g(L/L')\in\Bbb C^\times$ be the corresponding
automorphism of $\lambda_K (L/L')$.

{\it Example.} Suppose that $\lambda_K (L/L')$ has degree 0. Choose
a Fr\'echet isomorphism $\alpha : \Gamma (E/L) \iso \Gamma (E/L')$
which represents $\id_{E/L,E/L'}$. Then $g(L/L')$ is the Fredholm
determinant $\det_{\fF}( \alpha^{-1}g_{E/L'}\alpha g^{-1}_{E/L})$,
where $g_{E/L}$, $g_{E/L'}$ are multiplication by $g$ automorphisms
of $\Gamma (E/L)$, $ \Gamma (E/L')$.

Here is a formula for $g(L/L')$ (cf.~\cite{PS} 6.7, \cite{SW} 3.6).
Consider the line bundle $\det E|_{U\setminus K}$. Then
$L\in\Lambda_K (j_{K\cdot} E)$ yields a $K$-lattice $\det L|_U \in
\Lambda_K (j_{K\cdot}\det  E|_{U\setminus K})$. We can assume that
$U$ has no compact components. Then the line bundle $\det L|_U$ is
trivial; let $\theta_L$ be any its trivialization. For two lattices
$L$, $L'$ we get a function $\theta_L /\theta_{L'} \in \CO^\times
(U\setminus K)$. Consider the analytic symbol $\{ g, \theta_{L'}
/\theta_{L}\} \in H^1 (U\setminus K,\Bbb C^\times )$. Then (see 3.5
for the notation $\partial U$) \begin{equation}g(L/L')= \{ g,
\theta_{L'} /\theta_{L}\}(\partial U ).  \label{
3.6.5}\end{equation} To check (3.6.5), consider first the case when
$L$, $L'$ are meromorphically equivalent. Then $\theta_{L'}
/\theta_{L}$ is a meromorphic function on $U$, and both parts of
(3.6.5) evidently coincide with $g(\text{div}(\theta_{L'}
/\theta_{L}))$. The general case follows since for any $L$, $L'$ one
can find (possibly enlarging $K$, as in (ii)) an $L''$
meromorphically equivalent to $L$ which is arbitrary close to $L'$,
and both parts of (3.6.5) depend continuously on $L'$.

\subsection{}\label{3.7.}
Let $N$ be a holonomic $\CD^\infty$-module on $X$. By 3.4, it
carries a canonical Fr\'echet structure. For $K$  as in 3.5, 3.6,
let us define a relative determinant theory $\mu^\nabla_K =\mu
(N/\omega N)^\nabla_K \in \CD et_K (N/\omega N)$ (cf.~2.4).

If $L$, $L_\omega$ are  $K$-lattices in $N$,  $\omega N$ such that
$\nabla (L)\subset L_\omega$, then $\nabla$ yields a morphism of
sheaves $N/L\to\omega N/L_\omega $, and we denote by $\CC
(L,L_\omega )_{N,K}$ its cone.

{\it Example.} Let  $M^\infty \iso N$ be a $\CD$-structure on $N$
(see 3.2), and $P$ be any finite subset of $K$ such that $M$ is
smooth on $K\setminus P$. Every $P$-lattices $L$, $L_\omega$  in
$M$, $\omega M$ can be viewed as $K$-lattices in $N$, $\omega N$. If
$\nabla (L)\subset L_\omega$, then we get an evident morphism of
complexes of sheaves (see 2.4)
\begin{equation}\CC (L,L_\omega )_{M,P}\to \CC (L,L_\omega )_{N,K}.
\label{ 3.7.1}\end{equation}

\bigbreak\noindent\textbf{Proposition.}\emph{
 (3.7.1) is a quasi-isomorphism.
}\bigbreak

{\it Proof.}  Our complexes are supported on a finite set $P$, so it
suffices to check that $R\Gamma (X, (3.7.1))$ is a
quasi-isomorphism. Notice that $dR (L,L_\omega ):= \CC one
(L\buildrel\nabla\over\to L_\omega )$ is a subcomplex of both
$dR(M)$ and $dR(N)$, and  $\CC (L,L_\omega )_{M,P}= dR (M)/ dR
(L,L_\omega )$, $\CC (L,L_\omega )_{N,T'}= dR (N)/ dR (L,L_\omega )
$. Since $dR(M)\to dR(N)$ is a quasi-isomorphism (see 3.2), we are
done. \hfill$\square$

By the corollary in 3.5, for $K$-lattices $L$, $L_\omega$ in $N$,
$\omega N$ one has a $\CF /\CI^{\text{tr}}$-morphism
\begin{equation}\nabla_{ N/L, \omega N/L_\omega}: \Gamma ( N/L)\to \Gamma (\omega N/L_\omega ). \label{ 3.7.2}\end{equation}

\bigbreak\noindent\textbf{Corollary.}\emph{ (3.7.2) is a Fredholm
map. }\bigbreak

{\it Proof.} By loc.~cit., the validity of the assertion does not
depend on the choice of $L$, $L_\omega$. So we can assume to be in
the situation of Example, and we are done by 2.4 and the
proposition.  \hfill$\square$

We define $\mu^\nabla_K \in \CD et_K (N/\omega N)$ as a relative
determinant theory such that for any $L\in\Lambda_K (N)$, $L_\omega
\in\Lambda_K (\omega N)$ one has
\begin{equation}
\mu^\nabla_K (L/L_\omega ):= \lambda_{\nabla_{ N/L, \omega
N/L_\omega}} , \label{ 3.7.3}
\end{equation}
and the structure
isomorphisms $\lambda_K (L'/L)\otimes \mu^\nabla_K (L/L_\omega
)\otimes \lambda_K (L_\omega /L'_\omega )\iso \mu^\nabla_K
(L'/L'_\omega )$ are compositions (3.3.2) for
$\id_{N/L',N/L}\nabla_{ N/L, \omega N/L_\omega}\id_{\omega
N/L_\omega ,\omega N/L'_\omega}=\nabla_{ N/L', \omega N/L'_\omega}$.

\subsection{}\label{3.8.}  Let  $\nu$ be any invertible holomorphic 1-form defined
on $U\setminus K$, where $U$ as in 3.5. The multiplication by $\nu$
isomorphism $j_{K\cdot }N|_U \iso j_{K\cdot }\omega N|_U$ yields an
identification of the $\CL$-groupoids $\Lambda_K (j_{K\cdot }N)\iso
\Lambda_K (j_{K\cdot }\omega N)$, hence a relative determinant
theory (see Remark (iii) in 3.6) \begin{equation}\mu^\nu_K =\mu
(N/\omega N)^\nu_K \in \CD et_K (j_{K\cdot }N/j_{K\cdot }\omega
N)=\CD et_K (N/\omega N). \label{ 3.8.1}\end{equation} Set
\begin{equation} \CE_{\text{dR}}(N)_{(K,\nu )}:= \mu^\nabla_K
\otimes (\mu^{\nu}_K )^{ \otimes -1} \in\CL . \label{
3.8.2}\end{equation} Thus for every $L\in\Lambda_K (j_{K\cdot }N)$
one has a canonical isomorphism
\begin{equation}\CE_{\text{dR}}(N)_{(K,\nu )}\iso \mu^\nabla_K (L/\nu L ); \label{ 3.8.3}\end{equation} for
two lattices $L$, $L'$ the corresponding identification
$\mu^\nabla_K (L/\nu L )\iso \mu^\nabla_K (L'/\nu L')$ is
\begin{equation}\mu^\nabla_K (L/\nu L ) \buildrel\sim\over\leftarrow
\lambda_K (L/L')\otimes \mu^\nabla_K (L'/\nu L' )\otimes \lambda_K
(\nu L/\nu L')^{\otimes -1} \iso  \mu^\nabla_K (L'/\nu L' ),\label{
3.8.4}\end{equation} where the first arrow is the composition, the
second one comes from the multiplication by $\nu$ identification
$\lambda_K (L/L')\iso \lambda_K (\nu L/\nu L')$.

The  construction does not depend on the auxiliary choice of $U$.
When $\nu$ varies holomorphically, $\CE_{\text{dR}}(N)_{(K,\nu)}$
form a holomorphic line bundle on the parameter space (by (3.8.2)
and Remark (vi) in 3.6). If $K=\sqcup K_\alpha$, then (3.6.3) yields
a factorization  (here  $\nu_\alpha$ are the restrictions of $\nu$
to neighborhoods of $K_\alpha$)
\begin{equation}\otimes\,\CE_{\text{dR}}(N)_{(K_\alpha ,\nu_\alpha
)}\iso \CE_{\text{dR}}(N)_{(K,\nu)}.\label{ 3.8.5}\end{equation}

\subsection{}\label{3.9.}
The above constructions are compatible with those from 2.5.
Precisely, let $M$, $P$ be as in Example in 3.7, and suppose that
$\nu$ is meromorphic on $U$ with $D:=$ div$(\nu )$ supported on $P$.
Then $\Lambda_P (M)\subset \Lambda_K (N)$, $\Lambda_P (M(\infty
P))\subset \Lambda_K (j_{K\cdot}N)$ (see loc.~cit.). These
embeddings are naturally compatible with the $\CL$-groupoid
structures, so \begin{equation}\CD et_P (M) =\CD et_K (M), \quad \CD
et_P (M/\omega M)= \CD et_K (N/\omega N). \label{
3.9.1}\end{equation} By the proposition in 3.7, (3.7.1) provides an
identification \begin{equation}\mu^\nabla_P (M/\omega M)\iso
\mu^\nabla_K  (N/\omega N). \label{ 3.9.2}\end{equation} Joint with
an evident isomorphism $\mu^\nu_P (M/\omega M)\iso \mu^\nu_K
(N/\omega N)$, it yields
\begin{equation}\CE_{\text{dR}}(M)_{(D,c,\nu_P )} \iso
\CE_{\text{dR}}(N)_{(K,\nu )}. \label{ 3.9.3}\end{equation} If $\nu$
varies holomorphically, then (3.9.3) is holomorphic.

\subsection{}\label{3.10.}
{\it Proof of the theorem in 3.1.} By 3.2,  $\phi : B(M)\iso B(M')$
amounts to an isomorphism of $\CD^\infty$-modules $ M^\infty \iso
M^{\prime\infty}$. Threfore we can view $M$ and $M'$ as two
$\CD$-structures on a holonomic $\CD^\infty$-module $N$. Choose the
multiplicity of $T$ to be compatible with both $M$ and $M'$ (see
2.1). We want to define an isomorphism $\phi^\epsilon :
\CE_{\text{dR}}(M)\iso \CE_{\text{dR}}(M')$ in $
\CL^\Phi_{\text{dR}}(X,T)$.

Let $S$ be an analytic space, $(D,c,\nu_P )\in\fD^\diamond (S)$ (see
1.1). We work locally on $S$, so we have  $T^c \subset T$. Choose a compact $K$, its open
neighborhood $U$, and  a meromorphic $\nu$ on $U_S$ such that $K$,
$U$ satisfy the conditions from 3.5, $P\subset K_S$, $D=-$
div$(\nu)$, and $\nu_P =\nu|_P$. We define $\phi^\epsilon $ at
$(D,c,\nu_P )$ as the composition  $ \CE_{\text{dR}}(M)_{(D,c,\nu_P
)}\iso \CE_{\text{dR}}(N )_{(K,\nu )}
\buildrel{\sim}\over\leftarrow\CE_{\text{dR}}(M')_{(D,c,\nu_P )};$ here
$\iso$ are (3.9.3). Equivalently, choose $L\in\Lambda_P (M)$,
$L'\in\Lambda_P (M')$; then $\phi^\epsilon$ is the composition $
\CE_{\text{dR}}(M)_{(D,c,\nu_P )}\iso \mu^\nabla_K (L/\nu L)\iso
\mu^\nabla_K (L'/\nu
L')\buildrel{\sim}\over\leftarrow\CE_{\text{dR}}(M')_{(D,c,\nu_P )}$,
the first and the last arrows are compositions of (3.9.2) and
(2.5.6), the middle one is (3.8.4).

The last description shows that $\phi^\epsilon$ does not depend on
the auxiliary choice of $\nu$. Indeed, $\nu$ is defined up to
multiplication by an invertible function $g$ on $U$ which equals 1
on $P$. By (3.8.4), replacing $\nu$ by $g\nu$ multiplies the
isomorphism $\mu^\nabla_K (L/\nu L)\iso \mu^\nabla_K (L'/\nu  L')$
by $g(L/L')$ (see Remark (vii) in 3.6). Since $T$ is compatible with
$M$ and $M'$, the 1-form $d\log (\theta_{L'}/\theta_{L})$ (see
loc.~cit.) has pole at $P$ of order $\le$ the multiplicity of $P$,
so (3.6.5) implies that  $g(L/L')=1$.

The construction is  compatible with factorization, so we have
defined $\phi^\epsilon$ as isomorphism in $\CL_\CO^\Phi (X,T)$. One
has $M|_{X\setminus T} = M'|_{X\setminus T}$, and
$\phi^\epsilon|_{X\setminus T}$ is the corresponding evident
identification. By the corollary in 1.12, this implies that
$\phi^\epsilon$ is horizontal, i.e., it is an isomorphism in
$\CL^\Phi_{\text{dR}}(X,T)$.

The compatibility of $\phi^\epsilon$ with constraints from 2.8 is
evident. Finally, the commutativity of (3.1.1)  follows from the
next proposition:

\subsection{}\label{3.11}

\bigbreak\noindent\textbf{Proposition.}\emph{ Suppose that $X$ is
compact, $N$ is a holonomic $\CD^\infty$-module smooth on
$X\setminus K$, and $\nu$ is a holomorphic invertible 1-form  on
$X\setminus K$. Then there is a canonical isomorphism
\begin{equation}\eta_{\text{dR}} : \CE_{\text{dR}}(N)_{(K,\nu )}\iso
\det R\Gamma_{\!\text{dR}}(X,N). \label{ 3.11.1}\end{equation} If
$\nu$ is meromorphic and $M$ is a $\CD$-structure on $N$, then the
next diagram of isomorphisms commutes (see 1.4 and 2.7 for the left
column, the top arrow is (3.9.2)):
\begin{equation}
\begin{matrix}
 \CE_{\text{dR}}(M)_\nu & \iso & \CE_{\text{dR}}(N)_{(K,\nu)}
 \\
 \eta_{\text{dR}} \downarrow \,\,\,\,& & \eta_{\text{dR}}\downarrow \,\,\,\,  \\
 \det H^\cdot_{\text{dR}}(X,M) & \iso
 &\det H^\cdot_{\text{dR}}(X,N).
\end{matrix} \label{ 3.11.2}\end{equation}
}\bigbreak

{\it Proof} (cf.~2.7). For $L\in\Lambda_K (j_{K\cdot }N)$ set
$\lambda (L):=\det R\Gamma (X,L)$. Then $\lambda$ is a determinant
theory on $j_{K\cdot }N$ at $K$ in an evident way. Replacing $N$ by
$\omega N$, we get $\lambda_\omega \in \CD et_K (j_{K\cdot }\omega
N)$, hence $\lambda \otimes\lambda_\omega^{\otimes -1}\in  \CD et_K
(j_{K\cdot }N/j_{K\cdot }\omega N)$. One has an isomorphism
\begin{equation}\mu^\nu_K \iso \lambda
\otimes\lambda_\omega^{\otimes -1}, \label{ 3.11.3}\end{equation}
namely, $\mu^\nu_K (L/L_\omega ):=\lambda_K (\nu L/L_\omega
)=\lambda_\omega (\nu L)\otimes \lambda_\omega (L_\omega )^{\otimes
-1}\iso\lambda (L)\otimes \lambda_\omega (L_\omega )^{\otimes -1}=(
\lambda \otimes\lambda_\omega^{\otimes -1})(L/L_\omega )$ where
$\iso$ comes from the isomorphism $\nu^{-1} : \nu L\iso L$.

For $K$-lattices $L$ in $N$, $L_\omega $ in  $N_\omega$ such that
$\nabla (L)\subset L_\omega$, set $dR (L,L_\omega ):=\CC one
(L\buildrel{\nabla}\over\to L_\omega )$. Since $dR(N)/dR (L,
L_\omega )= \CC (L,L_\omega )$, our $dR(N)$ carries a 3-step
filtration with successive quotients $L_\omega$, $L$, $\CC
(L,L_\omega )$. Applying $\det R\Gamma$, we get an isomorphism
\begin{equation}\det R\Gamma (X,\CC (L,L_\omega ))\otimes \lambda
(L)^{\otimes -1}\otimes \lambda (L_\omega )\iso  \det
R\Gamma_{\!\text{dR}}(X,N). \label{ 3.11.4}\end{equation} To get $
\eta_{\text{dR}}$, we combine (3.11.4) with (3.11.3) (and (3.7.3)).
The construction does not   depend on the auxiliary choice  of $L$,
$L_\omega$. \hfill$\square$

\section{The Betti $\varepsilon$-line}

We present a construction from \cite{B} in a format adapted for
the current subject. {\it In 4.2--4.5  $X$ is considered as a mere
real-analytic surface. }

\subsection{}\label{4.1.}
Let $\CL$ be any Picard groupoid. For a (non-unital) Boolean
algebra\footnote{Recall that a Boolean algebra is the same as a
commutative $\Bbb Z/2$-algebra each of whose elements is idempotent;
the basic Boolean operations are $x\cap y = xy$, $x\cup y = x+y
+xy$; elements $x$, $y$ are said to be disjoint if $x\cap y=0$. The
Boolean algebras we meet are already realized as  Boolean algebras
of subsets of some set.} $\CC$, an {\it $\CL$-valued measure}
$\lambda$ on $\CC $ is a rule that assigns to every $S\in \CC $ an
object $\lambda (S)\in \CL$, and to every finite collection $\{
S_\alpha \}$ of {\it pairwise disjoint} elements of $\CC$ an
identification $\otimes  \lambda (S_\alpha )\iso \lambda (\cup
S_\alpha )$ (referred to as {\it integration}); the latter  should
satisfy an evident transitivity property. Such $\lambda$ form
naturally a Picard groupoid $\CM (\CC,\CL )$.

{\it Remarks.} (i)  For an abelian group $A$ denote by $\CM (\CC
,A)$ the group of $A$-valued measures on $\CC $. Then $\pi_1 (\CM
(\CC ,\CL ))= \CM (\CC ,\pi_1 (\CL ))$, and one has a map $\pi_0
(\CM (\CC ,\CL ))\to \CM (\CC,\pi_0 (\CL))$ which assigns to
$[\lambda ]$ a $\pi_0 (\CL )$-valued measure $|\lambda |$, $|\lambda
|(S):= [\lambda (S)]$ (see 1.1 for the notation).

(ii) Let $\CI\subset \CC$ be an ideal. Then $\CM (\CC/\CI,\CL)$
identifies naturally with the Picard groupoid of pairs $(\lambda
,\tau )$ where $\lambda \in\CM (\CC,\CL )$ and $\tau $ is a
trivialization of its restriction $\lambda|_{\CI}$  to $\CI$, i.e.,
an isomorphism $1_{\CM (\CI,\CL )}\iso \lambda|_{\CI}$ in $\CM
(\CI,\CL )$.

(iii) Suppose $\CC$ is finite, i.e., $\CC= (\Bbb Z/2)^T :=$ the
Boolean algebra of all subsets of a finite set $T$. Then an
$\CL$-valued measure $\lambda$ on $\CC$ is the same as a collection
of objects $\lambda_t = \lambda (\{ t\})$, $t\in T$. Thus $\CM
(T,\CL ):= \CM (\Bbb Z/2^T ,\CL )\iso \CL^T$.

\subsection{}\label{4.2.}
For an open $U\subset X$ we denote by $\CC (U)$ the (non-unital)
Boolean algebra of relatively compact subanalytic subsets of $U$.
For $U'\subset U$ one has $\CC(U')\subset \CC(U)$, and $U \mapsto
\CM (\CC (U),\CL )$ is a sheaf of Picard groupoids on $X$.

For a commutative ring $R$, let $\CL_R$ be the Picard groupoid of
$\Bbb Z$-graded super $R$-lines. Its objects are pairs $L=(L,\deg
(L))$ where $L$ is an invertible $R$-module, $\deg (L)$ a locally
constant $\Bbb  Z$-valued function on $\Spec R$; the commutativity
constraint is ``super" one. Every perfect $R$-complex $F$ yields a
graded super line $\det F \in\CL_R$. For a finite filtration
$F_\cdot$ on $F$ by perfect subcomplexes, one has a canonical
isomorphism $\det F \iso \otimes\det \gr_i F$; it satisfies
transitivity property with respect to refinement of the filtration.
For a finite collection $\{ F_\alpha \}$, every linear ordering of
$\alpha$'s yields a filtration on $\oplus F_\alpha$, and the
corresponding isomorphism $\det \oplus F_\alpha \iso \otimes\det
F_\alpha$ is independent of the ordering; thus $\det$ is a symmetric
monoidal functor.

Let $F=F_U$ be a perfect constructible complex of $R$-sheaves on
$U$. Then for every locally closed subanalytic subset $i_C : C\hra
U$ the $R$-complex $R\Gamma (C, Ri^!_C F )$ is perfect, so we have
$\det R\Gamma (C, Ri^!_C F)\in\CL_R$.

Suppose we have a finite closed subanalytic
filtration\footnote{I.e., each $C_{\le i}$ is a subanalytic subset
closed in $C$.} $C_{\le\cdot}$ on $C$ (therefore $C_i := C_{\le
i}\setminus C_{\le i-1}$ are locally closed and form a partition of
$C$). It yields a finite filtration\footnote{A filtration on an
object $C$ of a derived category is an object of the corresponding
filtered derived category identified with $C$ after the forgetting
of the filtration.} on $R\Gamma (C, Ri^!_C F )$ with $\gr_i R\Gamma
(C, Ri^!_C F ) =  R\Gamma (C_i , Ri^!_{C_i } F) $, hence an
identification \begin{equation}\otimes\, \det  R\Gamma (C_i ,
Ri^!_{C_i }F)  \iso  \det R\Gamma (C, Ri^!_C F). \label{
4.2.1}\end{equation} It satisfies transitivity property with respect
to refinement of the filtration.

\bigbreak\noindent\textbf{Lemma.}\emph{ There is a unique (up to a
unique isomorphism) pair $(\lambda (F), \iota )$, where $\lambda
(F)\in \CM (\CC (U) ,\CL_R )$,  $\iota$ is a datum of isomorphisms
$$\iota_C: \lambda (F)(C)\iso \det R\Gamma (C, Ri^!_C F)$$ defined for
any locally closed $C$, such that for every filtration
$C_{\le\cdot}$ on $C$ as above,  $\iota$ identifies (4.2.1) with the
integration $\otimes\, \lambda (F)(C_i )\iso \lambda (F)(C)$.
}\bigbreak

{\it Proof.} Suppose we have a compact subanalytic subset of $U$
equipped with a subanalytic stratification whose strata $C_\alpha$
are smooth and connected. The strata generate  a Boolean subalgebra
$\CC(\{ C_\alpha \})$ of $\CC(U)$;  call a subalgebra of such type
{\it nice}. Every finite subset of $\CC(U)$ lies in a nice
subalgebra; in particular, the set of nice subalgebras is directed.
To prove the lemma, it suffice to define the restriction of
$(\lambda (F),\iota )$ to every nice $\CC(\{ C_\alpha \})$; their
compatibility is automatic.

By   Remark (iii) in 4.1, $\lambda (F)|_{\CC(\{ C_\alpha \})}$ is
the measure defined by  condition\linebreak $\lambda (F)(C_\alpha )=\det
R\Gamma (C_\alpha ,Ri^!_{C_\alpha} F)$. For a locally closed  $C$ in
$\CC(\{ C_\alpha  \})$ one defines $\iota_C$ using (4.2.1) for a
closed filtration on $C$ whose layers are strata of increasing
dimension; its independence of the choice of filtration follows
since $\det$ is a symmetric monoidal functor in the way described
above. The compatibility with (4.2.1) is checked  by induction by
the number of strata involved. \hfill$\square$

{\it Example.} Suppose  $C', C'' \in \CC(U)$ are such that $C'$,
$C''$, and $C:=C'\cup C''$ are locally closed, $C'\cap
C''=\emptyset$. By the lemma, there is a canonical isomorphism
\begin{equation}\det R\Gamma (C' ,Ri^!_{C'} F)\otimes\det R\Gamma (C'' ,Ri^!_{C''}
F)\iso \det R\Gamma (C ,Ri^!_{C} F). \label{ 4.2.2 }\end{equation}
If, say, $C'$ is closed in $C$, then this is  (4.2.1) for the
filtration $C'\subset C$. To construct (4.2.2) when neither $C'$ nor
$C''$ are closed (e.g.~$C=X$ is a torus, and $C'$, $C''$ are
non-closed annuli), consider a 3-step closed filtration $C'\subset
\bar{C}'\subset C$, where $\bar{C}'$ is the closure of $C'$ in $C$;
set $P:= \bar{C}'\setminus C' = \bar{C}' \cap C''$, $Q:= C\setminus
\bar{C}' =C''\setminus P$. By (4.2.1), $\det R\Gamma (C' ,Ri^!_{C'}
F)\otimes\det R\Gamma (P ,Ri^!_{P} F) \iso \det R\Gamma (\bar{C}'
,Ri^!_{\bar{C}'} F)$, $\det R\Gamma (P ,Ri^!_{P} F)\otimes \det
R\Gamma (Q ,Ri^!_{Q} F)\iso \det R\Gamma (C'', Ri^!_{C''}F)$, and
$\det R\Gamma (\bar{C}' ,Ri^!_{\bar{C}'} F)\otimes \det R\Gamma (Q
,Ri^!_{Q} F)$ $\iso \det R\Gamma (C ,Ri^!_{C} F)$. Combining them,
we get (4.2.2).

\subsection{}\label{4.3.} Let $U\subset X$ be an open subset,
and $\CN=\CN_U \subset T_U$ be a continuous family of proper cones
in the tangent bundle (so for each $x\in U$ the fiber $\CN_x $ is a
proper closed sector with non-empty interior in the tangent plane
$T_x$). For an open $V\subset U$ we denote by $\CN_V$ the
restriction of $\CN$ to $V$.

One calls $C\in \CC (U)$ an {\it $\CN$-lens} if it satisfies the
next two conditions:

(a) Every point in $U$ has a neighborhood $V$ such that $C\cap V =
C_1 \setminus C_2$ where $C_1$, $C_2$ are closed subsets of $V$ that
are invariant with respect to some family of proper cones $\CN'_V
\Supset \CN_V$.\footnote{Here $\Supset$ means that Int$\,\CN'_V
\supset \CN_V \setminus \{ 0\}$, and $\CN'_V$-invariance of $C_i$
means that every $C^1$-arc $\gamma : [0,1]\to V$ such that $\gamma
(0)\in C_i$ and $\frac{d}{dt}\gamma (t) \in \CN'_V \setminus \{ 0\}$
for every $t$, lies in $C_i$.}

(b) There is a $C^1$-function $f$ defined on a neighborhood $V$
 of the closure $\bar{C}$ such that for every $x\in V$ and a non-zero $\tau \in \CN_x$ one has $\tau (f)>0$.

Let $\CI (U,\CN )\subset  \CC (U) $ be  the Boolean subalgebra
generated by all $\CN$-lenses.

 Basic properties of lenses (see \cite{B} 2.4, 2.7):
(i) Every $\CN$-lens $C$ is locally closed, and Int$\, C$  is dense
in $C$; the intersection of two $\CN$-lenses is an $\CN$-lens.
\newline (ii) Every point in $U$ admits a base of neighborhoods
formed by $\CN$-lenses. \newline (iii) Suppose we have an $\CN$-lens
$C$ and a (finite) partition $\{ C_\alpha \}$ of $C$ with $C_\alpha
\in \CI (U,\CN )$. Then there exists a finer partition $\{ C_1
,\ldots ,C_n \}$ of $C$ such that $C_i$ are $\CN$-lenses and each
subset $C_{\le i}:=C_1 \cup C_2 \cup \ldots \cup C_i$ is closed in
$C$.

 {\it Exercise.} Every $C\in \CC(U)$ that satisfies (a) lies in $\CI(U,\CN )$.

Suppose $F$ from 4.2 is locally constant (say,  a local system of
finitely generated projective $R$-modules). Then for any $\CN$-lens
$C$ one has $R\Gamma (C, Ri^!_C F)=0$ (see \cite{B} 2.5). Let
$\tau_C : 1 \iso \lambda (F)(C)$ be the corresponding trivialization
of the determinant.

\bigbreak\noindent\textbf{Proposition.}\emph{ The restriction of
$\lambda (F)$ to $\CI (U,\CN )$ admits  a unique trivialization
$\tau_\CN : 1_{\CM (\CI(U,\CN ),\CL_R )} \iso \lambda
(F)|_{\CI(U,\CN )}$ such that for every $\CN$-lens $C$ the
trivialization $\tau_{\CN C}$ coincides with $\tau_C$. }\bigbreak

{\it Proof.} By (iii) above, every finite subset of $ \CI(U,\CN )$
lies in the Boolean subalgebra generated by a finite subset of
pairwise disjoint $\CN$-lenses. This implies uniqueness. To show
that $\tau_\CN$ exists, it suffices to check the next assertion: For
any $\CN$-lens $C$ and any finite partition $\{ C_\alpha \}$ of $C$
by $\CN$-lenses the integration $\otimes \lambda (F)(C_\alpha )\iso
\lambda (F)(C)$ identifies $\otimes\tau_{C_\alpha}$ with $\tau_C$.

Choose $\{ C_1 ,\ldots ,C_n \}$ as in (iii) above. Since $C_{\le
\cdot}$ is a closed filtration,  the integration $\otimes \lambda
(F)(C_i )\iso \lambda (F)(C)$ identifies $\otimes \tau_{C_i}$ with
$\tau_C$ (see the lemma in 4.2), and   for each $\alpha$ the
integration $\otimes \lambda (F)(C_i \cap C_\alpha )\iso \lambda
(F)(C_\alpha )$ identifies $\otimes \tau_{C_i \cap C_\alpha}$ with
$\tau_{C_\alpha}$. The partition $\{ C_i \}$ is finer than $\{
C_\alpha \}$, so we are done by the transitivity of integration.
\hfill$\square$

\subsection{}\label{4.4.} Let $K$ be a compact subset of $X$, $W$ be an open subset
that contains $K$, $U:= W\setminus K$. For $\CN =\CN_U$ as above,
let   $\CC (W, \CN )$ be  the set of $C\in \CC(W)   $ that satisfy
the next two conditions:

(a) For every $C'\in \CI(U,\CN )$ one has $C\cap C' \in \CI(U,\CN
)$.

(b) One has Int$(C)\cap K= \bar{C} \cap K$.

Then  $\CC (W, \CN )$ is a Boolean subalgebra of $\CC(W)$, and $
\CI(U,\CN )$ is an ideal in it.

{\it Exercise.} Let $K'$ be a subset of $K$ which is open and closed
in $K$. Choose an open relatively compact subset $V$ of $W$ such
that $V\cap K=\bar{V}\cap K=K'$, and $C' \in \CI(U,\CN)$ such that
$C'\supset \partial V := \bar{V}\setminus V$. Then $C:= V\setminus
C' \in \CC(W,\CN )$ and $C\cap K= K'$.

Denote by $\CC[K]$ the Boolean algebra of subsets  of $K$ which are
open and closed in $K$.
 By (b), we have a morphism of Boolean algebras
 $ \CC (W, \CN )
\to \CC[K]$, $C\mapsto C\cap K$. It yields an identification
\begin{equation} \CC (W, \CN )/ \CI(U,\CN )\iso \CC[K] . \label{
4.4.1}\end{equation}

Let $F=F_W$ be a perfect constructible complex of $R$-sheaves on $W$
whose restriction to $U$ is locally constant. By the proposition in
4.3, we have a trivialization $\tau_\CN$ of the restriction of
$\lambda (F)$ to $\CI(U,\CN )$. By Remark (ii) in 4.1, (4.4.1),  the
pair $(\lambda (F)|_{\CC(W,\CN )},\tau_N )$ can be viewed as a
measure $\CE (F)_{\CN }\in \CM (\CC[K], \CL_R )$. If $K$ is finite,
then, by Remark (iii) in 4.1, it amounts to a collection of lines
$\CE (F)_{(b,\CN )}:=\CE (F)_{\CN b} $,
  $b\in K$.

If $C\in \CC(W,\CN )$ is locally closed, then we have
identifications \begin{equation}\CE (F)_{\CN}(C\cap K)\iso \lambda
(F)(C)\iso \det R\Gamma (C, Ri^!_C F). \label{ 4.4.2}\end{equation}
In particular, if $X$ is compact and $W=X$, then $X\in \CC(X, \CN
)$, and
\begin{equation} \CE (F)_{\CN }(K ) \iso  \det R\Gamma (X, F). \label{
4.4.3}\end{equation} If $K$ is finite, this is a product formula
\begin{equation} \mathop\otimes\limits_{b\in K} \CE (F)_{(b,\CN )} \iso
\det R\Gamma (X, F). \label{ 4.4.4}\end{equation}

\subsection{}\label{4.5.}
Suppose we have another datum of $W' \supset U'$, $\CN' =\CN'_{U'}$
as above, such that $W' \subset W$, $U' \subset U$, and $\CN'
\supset \CN_{U'}$. Then $\CI(U',\CN'  )\subset \CI(U,\CN )$,
$\CC(W', \CN' )\subset \CC(W, \CN )$, and (4.4.1) identifies the
morphism of Boolean algebras $\CC(W', \CN' )/ \CI(U',\CN' )\to
\CC(W,\CN )/\CI(U,\CN )$ with a  morphism $r : \CC[K']\to   \CC[K]$,
$Q\mapsto r(Q):=Q\setminus U= Q\cap K$. Since $\tau_{\CN'}$ equals
the restriction of $\tau_\CN$ to  $\CI(U',\CN' )$, one has
\begin{equation}\CE (F)_{\CN' }=r^* \CE (F )_{\CN } . \label{
4.5.1}\end{equation}

{\it Remarks.} (i) Taking for $W' $ a small neighborhood of a
component $K'$ of $K$, $U'= W'\cap U$,  $\CN' = \CN|_{U'}$, we see
that $\CE (F)_{\CN }$ has local nature with respect to $K$.

(ii) By (4.5.1), $\CE (F)_\CN (W' \cap K)$ depends only on the
restriction of $\CN$ to $U'$.

\subsection{}\label{4.6.} Suppose now $X$ is a complex curve, $T\subset X$ a finite
subset, $F$ is a constructible sheaf on $X$ which is smooth on
$X\setminus T$. Let us define a constructible factorization $R$-line
$\CE_{\text{B}}(F)$ on $(X,T)$ (see 1.15).

Let $S$ be an analytic space, $(D,c,\nu_P )\in\fD^\diamond (S)$. Let
us define a local system of $R$-lines $\CE (F)_{(D,c,\nu_P )}$ on
$S$. Consider  a datum $(W,K, \CN, \nu_S )$, where $ W$ is an open
subset of $X$, $K$ a compact subset of $W$, $\CN =\CN_U$ is a
continuous family of proper cones in the tangent bundle to
$U:=W\setminus K$ (viewed as a real-analytic surface, see 4.3),
$\nu_S$ is an $S$-family  of meromorphic 1-forms on $W$. We say that
our datum is {\it compatible} if $P=P_{D,c}\subset K_S$, div$(\nu
)=-D$,  $\nu|_P =\nu_P$, and the 1-forms Re$(\nu_s )$ are negative
on $\CN$. As in 4.4, every compatible datum yields the $R$-line $
\CE (F)_\CN (K)$.

\bigbreak\noindent\textbf{Lemma.}\emph{ Locally on $S$ compatible
data exist; the lines $\CE (F)_\CN (K)$ for all compatible data
are naturally identified. }\bigbreak

{\it Proof.} The existence statement is evident. Suppose that we fix
an open subset $W_0$ of $X$ and an $S$-family of meromorphic forms
$\nu_S$ on $W_0$ such that $P\subset W_{0S}$, div$(\nu_S )=D$, and
$\nu_P =\nu|_P$. Let us consider compatible data with $W\subset W_0$
and the above $\nu_S$. The identification of the lines for these
data comes from 4.5. Thus our line depends only on $\nu_S$; in fact,
by Remark (i) in 4.5, on the germ of $\nu_S$ at $P$. If we move
$\nu_S$, it remains locally constant.  Since the space of germs of
$\nu_S$ is contractible, we are done. \hfill$\square$

Locally on $S$, we define $\CE (F)_{(D,c,\nu_P )}$ as $\CE (F)_\CN
(K)$ for a compatible datum. The factorization structure is evident.
For $X$ compact, we have, by (4.4.3), a canonical identification
\begin{equation}\eta : \CE (F)(X)\iso \det R\Gamma (X,F) . \label{
4.6.1}\end{equation}

{\it Exercise.} Check that $\CE$ satisfies the constraints from 2.8.

{\it Remark.} Suppose $X$ is compact and a rational form $\nu$ has
property that Re$(\nu )$ is exact, Re$(\nu )=df$. Then the
isomorphism $\eta : \CE (F)_\nu \iso \det R\Gamma (X,F)$ can be
computed using Morse theory: indeed, if $a < a'$ are non-critical
values of $f$, then $f^{-1}( (a',a])\in C(\CN )$ for $\CN$
compatible with $\nu$.

In \S5 we apply this construction to $F=B(M)$, the de Rham complex
of a holonomic $\CD$-module $M$, and write $\CE_{\text{B}} (M):= \CE
(B(M))$. We use the same notation for the corresponding de Rham
factorization line (in the analytic setting).

\subsection{}\label{4.7.}
For $b\in X$ and a meromorphic $\nu$ on a neighborhood of $b$, $v_b
(\nu )=-\ell $, set $\CE (F)_{(b,\nu )}:=\CE (F)_{(\ell b, \nu )}$.
Let $F_b^{(! )}:= Ri^!_b F=R\Gamma_{\{b\}} (X,F)=R\Gamma_c (X_b
,F)$, $F_b^{(*)}:=i_b^* F = R\Gamma (X_b ,F)$ be the fibers of $F$
at $b$ in !- and $*$-sense (here $X_b$ is a small open disc around
$b$). Let $t$ be a local parameter at $b$.

\bigbreak\noindent\textbf{Lemma.}\emph{ There are canonical
identifications
\begin{equation}
\CE (F)_{(b,t^{-1}dt)} =\det F_b^{(! )} ,\quad\CE (F)_{(b,-t^{-1}dt)} =\det F_b^{(*)}  .\label{ 4.7.1}
\end{equation}
}\bigbreak

 {\it Proof.} Let
$X_b$ be a small disc $|t|< r$. Let $W$  be the open disc of radius
$r'$,  $K$ be the  closed disc of radius $r''$, $r'' <r<r'$. If
$\CN$ is a sufficiently tight cone around the Euler vector field
Re$(t\partial_t )$, then $X_b \in C(W,\CN )$  and $\bar{X}_b \in
C(W,-\CN )$. The data $( W,K, \CN, -t^{-1}dt )$ and $( W,K, -\CN,
t^{-1}dt )$ are compatible, and we are done. \hfill$\square$

Thus if $F$ is the $*$-extension at $b$, i.e., $F^{(!)}_b =0$, then
$\CE (F)_{(b,t^{-1}dt)}$ is canonically trivialized; if $F$ is the
$!$-extension at $b$, i.e., $F^{(*)}_b =0$, then $\CE
(F)_{(b,-t^{-1}dt)}$ is canonically trivialized. Denote these
trivializations by $1_b^! \in \CE (F)_{(b,t^{-1}dt)}$, $1_{b}^* \in
\CE (F)_{(b,-t^{-1}dt)}$.

{\it Exercise.} Suppose $X=\Bbb P^1$ and $F$ is smooth outside 0,
$\infty$. Then the composition $\det F^{(!)}_0 \otimes \det
F^{(*)}_\infty \iso \CE (F)_{(0,t^{-1}dt )}\otimes \CE (F)_{(\infty
,t^{-1}dt )} \buildrel{\eta}\over\to \det R\Gamma (\Bbb P^1 ,F)$
comes from the standard triangle $R\Gamma_{\{0\}} (\Bbb P^1 ,F)\to
R\Gamma (\Bbb P^1 ,F)\to R\Gamma (\Bbb P^1\setminus \{ 0\} ,F)=
F^{(*)}_\infty$. In particular, if $F$ is $*$-extension at 0 and
$!$-extension at $\infty$, then  $\eta (1_0^! \otimes 1^*_\infty )
=1 $ := the trivialization of $\det R\Gamma (\Bbb P^1 ,F)$ that
comes since $R\Gamma (\Bbb P^1 ,F)=0$.

For $x\in X\setminus T$ one has $F_x^{(!)}=F^{(*)}_x (-1)[-2]$,  so
(4.7.1)  yields a natural identification \begin{equation}\CE
(F)^{(1)}_{X\setminus T}\iso\det F(-1)_{X\setminus T}.\label{
4.7.2}\end{equation} If $R=\Bbb C$, then the Tate twist   acts as
identity. If $M$ is a holonomic $\CD$-module, then $B(M)_{X\setminus
T}=M^\nabla_{X\setminus T}[1]$, and (4.7.2) can be rewritten as
\begin{equation}\CE_{\text{B}}(M)^{(1)}_{X\setminus T}\iso (\det
M_{X\setminus T})^{\otimes -1}. \label{ 4.7.3}\end{equation}

{\it Remark.} For any  $\ell\in \Bbb Z$ we have the local system of
lines $\CE (F)_{(b,zt^{-\ell}dt)}$, $z\in\Bbb C^\times$. A simple
computation (or a reference to the compatibility property in 1.11)
together with (4.7.2) shows that its monodromy  around $z=0$ equals
$(-1)^{\text{rk}(F)\ell} m_b$ where rk$(F)=\deg \det F_{X\setminus
T}$ is the rank of $F$ and $m_b$ is the monodromy of $\det F$ around
$b$. Thus (4.7.1) provides  two descriptions of this local system
for $\ell =1$ (using the fibers at $z=\pm 1$). For a relation
between them, see below.

\subsection{}\label{4.8.} We are in the situation of 4.7.
Suppose $R$ is a field, outside singular points our $F$ is a local
system of rank 1 placed in degree $-1$, and $F^{(*)}_b =0$. Let $m$
be the monodromy of $F$ around $b$; suppose $m\neq 1$. Then
$F^{(!)}_b$ vanishes as well, so we have $1_b^! \in \CE
(F)_{(b,t^{-1}dt)}$, $1_{b}^* \in \CE (F)_{(b,-t^{-1}dt)}$.

\bigbreak\noindent\textbf{Proposition.}\emph{ The (counterclockwise)
monodromy from $t^{-1}dt$ to $-t^{-1}dt$ identifies $1^!_{b}$ with
$(1-m)^{-1} 1^*_{b}$. }\bigbreak

{\it Proof.} Consider an annulus around $b$ (which lies in $U$), and
cut it like this:

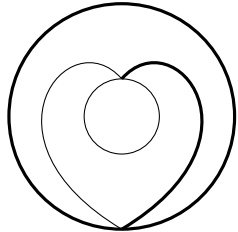
\begin{center}
\begin{tikzpicture}
\draw (0,0) circle(5mm); \draw[very thick] (0,0) circle(15mm);
\draw[very thick] (0,0.5) .. controls (0.7,1.3) and (2,-0.3) ..
(0,-1.5); \draw (0,0.5) .. controls (-0.7,1.3) and (-2,-0.3) ..
(0,-1.5);
\end{tikzpicture}
\end{center}

Let $D_+$ be the larger open disc,  $D_-$ the smaller one,
$\bar{D}_\pm$ be their closures. Let $\heartsuit$ be a constructible
set obtained from the closed heart figure by removing  the right
part of its boundary together with the lower vertex; the upper
vertex stays in $\heartsuit$. Set $A:=\bar{D}_+ \setminus \heartsuit
$,  $B:=\heartsuit \setminus D_-\in \CI(U,\CN_- )$, $C:=A\cup B=
\bar{D}_+ \setminus D_-$.

Suppose $\epsilon >0$ is small; let $\theta_\pm$ be the real part of
the complex vector field $\exp( i\pi /2 \pm i\epsilon )t\partial_t
$, and  $\CN_\pm$ be a tight cone around $\theta_\pm$. Then $\CN_+$
is compatible with $\exp (i\alpha )t^{-1}dt$ for  $\alpha\in
[0,\pi/2]$, and $\CN_-$ is compatible with $\exp (i\beta )t^{-1}dt$
for  $\beta\in [\pi/2,\pi]$. Our $\theta_\pm $ are transversal to
the lines of the drawing. Then  $A \in \CI(U,\CN_+ )$ and $B\in
\CI(U,\CN_- )$ (each of them is the union of two lenses), hence
$\heartsuit \in \CC (W,\CN_\pm )$.

Thus  the monodromy in the statement of the proposition is inverse
to the composition $\lambda (F)(D_- )\iso \lambda (F)(\heartsuit
)\iso \lambda (F)(\bar{D}_+ )$ where the first arrow is
multiplication by  $\tau_{\CN_-}\in \lambda (F)(B)$, the second one
is multiplication by $\tau_{\CN_+}\in \lambda (F)(A)$ (we use
tacitly the integration). Notice that $R\Gamma_C (X,F)$ is acyclic;
let $\tau$ be the corresponding trivialization of $\lambda (F)(C)$.
Since the multiplication by $\tau$ map $ \lambda (F)(D_- )\iso
\lambda (F)(\bar{D}_+ )$ sends $1^*_{b}$ to $1^!_{b}$, the
proposition can be restated as $\tau_{\CN_+}
\tau_{\CN_-}=(1-m)\tau$.

Consider the chain complex $(P,d)$ that computes $R\Gamma_C (X,F)$
by means of the cell decomposition of the drawing. The graded vector
space $P$   is a direct sum of rank 1 components $P_\alpha$ labeled
by the cells. Set $P_A :=\mathop\oplus\limits_{\alpha \in
A}P_\alpha$, $P_B :=\mathop\oplus\limits_{\beta\in B} P_\beta$.
Since $A$ is the image of a 2-simplex with one face removed,  $P_A$
carries a differential $d_A$ such that $(P_A ,d_A)$ is the chain
complex of the simplex modulo the face. One defines $d_B$ in a
similar way. Both $(P_A , d_A )$ and $(P_B ,d_B )$ are acyclic, and
the corresponding trivializations of $\det P_A = \lambda (F)(A)$,
$\det P_B =\lambda (F)(B)$ equal $\tau_{\CN_+} $, $\tau_{\CN_-}$.

Let $P' =P'_A \oplus P'_B$ be sum of $P_\alpha$'s for $\alpha$ in
the boundary of $\heartsuit$. Then $P'$ is a subcomplex with respect
to both $d$ and $d_A \oplus d_B$, on $P/P'$ the differentials $d$
and $d_A \oplus d_B$ coincide, and the complex $P/P'$ is acyclic.
We see that $P'$ is acyclic with respect to both $d|_{P'}$ and $(d_A
\oplus d_B )|_{P'}$, and $\tau$, $\tau_{\CN_+}\tau_{\CN_-}$ are the
trivializations of $\det P'$ that correspond to these differentials.
Our $P'$ sits in degrees $0$, $1$. Choose base vectors $e_A \in
P_A^{\prime 0}$, $f_A \in P_A^{\prime 1}$, $e_B \in P_B^{\prime 0}$,
$f_B \in P^{\prime 1}_B$ such that $d_A (e_A )=f_A$, $d_B (e_B
)=f_B$, and $d(e_A )= f_A - f_B$. Then $d(e_B )= -m f_A +f_B$.
Therefore $\tau_{\CN_+}\tau_{\CN_-}/\tau = 1-m$, q.e.d.
\hfill$\square$

\subsection{}\label{4.9.} Let $b$, $\nu$, $\ell$ be as in the beginning of 4.7; suppose $\ell \neq 1$.
Denote by $\nu_b$ the principal term of $\nu$ at $b$. Let $f$  be
any holomorphic function defined near $b$ and vanishing at $b$  such
that  $\nu_b =(df)_b$ if $\ell <1$, and $\nu_b = (d(f^{-1}))_b$ if
$\ell >1$. For a small  $z\in\Bbb C$, $z\neq 0$, the set $f^{-1}(z
)$ is finite of order $|\ell-1|$. For a finite subset $Z$ of $X$,
denote by $F^{(!)}_Z$, $F^{(*)}_Z$  the direct sum of !-,
resp.~$*$-fibers of $F$ at points of $Z$. Let $\epsilon $ be a small
positive real number.

\bigbreak\noindent\textbf{Proposition.}\emph{  For $\ell <1$, one
has canonical identifications \begin{equation}\CE (F)_{(b,\nu )}\iso
\det F_b^{(!)} \otimes (\det F^{(!)}_{f^{-1}(-\epsilon )})^{\otimes
-1}\iso \det F^{(*)}_b \otimes ( \det F^{(*)}_{f^{-1}(\epsilon
)})^{\otimes -1} .
 \label{ 4.9.1}\end{equation}
For $\ell >1$, one has \begin{equation}\CE (F)_{(b,\nu )} \iso \det
F^{(*)}_b \otimes \det F^{(!)}_{g^{-1}(-\epsilon )}\iso  \det
F^{(!)}_b \otimes \det F^{(*)}_{g^{-1}(\epsilon )} .  \label{
4.9.2}\end{equation} }\bigbreak

{\it Proof.}  Since $\CE (F)_{(b,\nu )}$ depends only on the
principal term of $\nu$ at $b$, we can assume that $\nu$ equals $df$
or $d(f^{-1})$. Set $\lambda :=\lambda (F)$.

Case $\ell = 0$: Then $f$ is a local coordinate at $b$. Let $Q$ be
an open  romb around $b$ with vertices at $f =\pm \epsilon, \pm
i\epsilon$, $I$ and $I'$ be the parts of its boundary where Re$(f
)\le 0$, resp.~Re$(f )> 0$; set $C:=\bar{Q}\setminus I=Q\cup I'$.
Then one has $\CE (F)_{(b,\nu )}= \lambda (C)\iso  \lambda
(\bar{Q})\otimes \lambda (I)^{\otimes -1}\iso \lambda (Q)\otimes
\lambda (I')$. This yields (4.9.1) due to the next identifications:
\newline (a) $\lambda (\bar{Q})\iso \det R\Gamma_{\bar{Q}}
(X,F)\buildrel\sim\over\leftarrow \det F^{(!)}_b$ and $\lambda (Q)
\iso \det R\Gamma (Q,i^*_Q F) \iso \det F^{(*)}_b$; \newline (b)
$\lambda (I)\iso\det R\Gamma_I (X,F)$, and
$R\Gamma_{f^{-1}(-\epsilon)}(X,F)\iso R\Gamma_I (X,F)$; \newline (c)
$\lambda (I')\iso \det R\Gamma (I', Ri^!_{I'} F)  \iso \det R\Gamma
(I',i^*_{I'}F [-1])\iso ( \det F^{(*)}_{f^{-1}(\epsilon )})^{\otimes
-1}$, where the second isomorphism  comes from  $i^*_{I'}F[-1] \iso
R i^!_{I'} i_{Q!} i^*_{Q} F\iso R i^!_{I'} F$.

Case $\ell =2$: Then $f$ is a local coordinate at $b$. Define $Q$,
etc., as above. Then
 $\CE (F)_{(b,\nu )}=\lambda (Q \cup I)= \lambda (\bar{Q}\setminus I' )$, so (a)--(c) yield (4.9.2).

If $\ell \neq 1$ is arbitrary, then  $f$ is a $|\ell -1|$-sheeted
cover of a neighborhood of $b$ over a coordinate disc. The
projection formula compatibility $\CE (F)_{(b,\nu )} \iso \CE (f_*
F)_{(0,dt)}$ for $\ell <1$ and $\CE (F)_{(b,\nu )} \iso \CE (f_*
F)_{(0,d(t^{-1}))}$ for $\ell >1$ reduces the assertion to the cases
of $\ell $ equal to $0$ and $ 2$, and we are done.  \hfill$\square$

\section{The torsor of $\varepsilon$-periods.}

\subsection{}\label{5.1.}
We  consider triples $(X,T,M)$ where $X$ is a complex curve, $T$ its
finite subset, and $M$ is a holonomic $\CD$-module  on $(X,T)$
(i.e., a $\CD$-module on $X$ smooth off $T$). For us, a {\it weak
theory of $\varepsilon$-factors} is a rule $\CE$ that assigns to
every such triple a de Rham factorization line $\CE (M)$ on $(X,T)$
(in complex-analytic sense). Our $\CE$ should  be functorial with
respect to isomorphisms of triples, and have local nature, i.e., be
compatible with pull-backs by open embeddings. We  ask that:

(i) For a nice flat family  $(X/Q, T, M)$ (see 2.12) with reduced
$Q$ the factorization lines $\CE (M_q )$ vary holomorphically, i.e.,
we have $\CE (M)\in \CL^\Phi_{\text{dR}/Q} (X/Q ,T)$. If the family
is isomonodromic, then $\CE (M)\in \CL^\Phi_{\text{dR}} (X/Q ,T)$.

(ii)  $\CE (M)$ is multiplicative  with respect to finite
filtrations of $M$'s: for a finite filtration $M_\cdot$ on $M$ there
is a natural isomorphism $\CE (M)\iso \otimes\, \CE (\gr_i M)$.

(iii)  (projection formula) Let $\pi : (X',T')\to  (X,T)$ be a
finite morphism of pairs \'etale over $X\setminus T$, so, as in
Remarks (i), (ii) in 1.2, (ii) in 1.5, we have a morphism $\pi_* :
\CL^\Phi_{\text{dR}} (X',T')\to \CL^\Phi_{\text{dR}} (X,T)$. Then
for any $M'$ on $(X',T')$ one has a natural identification $\CE
(\pi_* M')\iso \pi_* \CE (M')$ compatible with composition of
$\pi$'s.

 (iv) (product formula)  For compact $X$
there is a natural  identification (see 1.4 for the notation) $ \eta
=\eta (M) : \CE (M) (X) \iso \det R\Gamma_{\!\text{dR}}(X,M)$.

The constraints  should be pairwise compatible in the evident sense.
(i) should be compatible with the base change.  (ii) should be
transitive with respect to refinements of the filtration, and the
isomorphism $\CE (\oplus M_\alpha )\iso \otimes \CE (M_\alpha )$
should not depend on the linear ordering of the indices $\alpha$
(which makes $\CE$ a symmetric monoidal functor). (iii) should be
compatible with the composition of $\pi$'s.

Weak theories of $\varepsilon$-factors form naturally a groupoid
which we denote by $ {}^w \Bbb E$. Its key objects are
$\CE_{\text{dR}}$ and $\CE_{\text{B}}$.

Replacing $\det R\Gamma_{\!\text{dR}}(X,M)$ in (iv) by the trivial
line $\Bbb C$ and leaving the rest of the story unchanged, we get a
groupoid ${}^w \Bbb E^0$. It has an evident Picard groupoid
structure, and $ {}^w \Bbb E$ is naturally a $ {}^w \Bbb
E^0$-torsor. Below we denote by (i)$^0$--(iv)$^0$ the structure
constraints in ${}^w \Bbb E^0$ that correspond to (i)--(iv) above.

\subsection{}\label{5.2.}
{\it Compatibility with quadratic degenerations of $X$.} Let us
formulate an important property of $\CE\in{}^w \Bbb E$.

Suppose we have data 2.13(a),(b) in the analytic setting; we follow
the notation of loc.~cit. In 2.13 we worked in the formal scheme
setting. Now the whole story of (2.13.1)--(2.13.7) makes sense
analytically: we have a proper family of curves $X$ over a small
coordinate disc $Q$ and an $\CO_X$-module $M$ equipped with a
relative connection $\nabla$, etc., so that the picture of 2.13
coincides with the formal completion of the present one at $q=0$. To
construct $X$, notice that $t_\pm$ from 2.13(a) converge on some
true neighborhoods $U_\pm$ of $b_\pm$. Suppose that $t_\pm$ identify
$U_\pm$ with coordinate discs of radii $r_\pm$ and $U_+ \cap U_-
=U_\pm \cap (T\cup |D|)= \emptyset$. Then $Q$ is the coordinate disc
of radius $r_+ r_-$. Let $W$ be an open subset  of $Y\times Q$
formed by those pairs $(y,q)$ that if $y\in U_+$, then $r_- |t_+
(y)|^2 > |q| r_+ $, and if $y\in U_-$, then $r_+ |t_- (y)|^2
>|q|r_-$. Our
$X$ is the union of two open subsets $V:=U_+ \times U_-$ and $W$: we
glue $(y,q)\in W$ such that $y\in U_\pm$ with $(y_+ ,y_- )\in V$
such that either $y_+$ or $y_-$ equals $y$ and $t_+ (y_+ )t_- (y_-
)= q$.  The projection $q :X\to Q$ is $(y,q)\mapsto q$ on $W$ and
$q(y_+ ,y_- )= t_+ (y_+ )t_- (y_- )$ on $V $.
  Set $K_{0Q} := X\setminus W = \{ (y_+ ,y_- )\in V : r_+^{-1}|t_+ (y_+ )|= r_-^{-1}|t_- (y_- ) |\}$.
The formal trivializations of $L$ at $b_\pm$ from 2.13(b) converge
on $U_\pm$, and we define $M$ by gluing $M_{b_0} \otimes \CO_{V}$
and the pull-back $N_W$ of $N$ by the projection $W\to Y$.

Let us define an analytic version of (2.13.7), which is an
isomorphism of $\CO_Q$-lines \begin{equation}\det  Rq_{\text{dR}*}
M\iso \det R\Gamma_{\!\text{dR}} (Y, N)\otimes\CO_Q . \label{
5.2.1}\end{equation} Let $i$, $j$ be the embeddings $K_{0Q} \hra
X\hookleftarrow W$. Since $V\setminus K_{0Q}  $ is disjoint union of
two open subsets $V_\pm$, $V_+ := \{ (y_+ ,y_- ): r_+^{-1}|t_+ (y_+
)|>r_-^{-1}|t_- (y_- )|\}$, the complex $\CF :=i^* j_* dR_{W/Q}(N_W
)$ is the direct sum of the two components $\CF_\pm$. Both maps $i^*
dR_{X/Q}(M) \to \CF_\pm$ are quasi-isomorphisms, so $i_* \CF_- \iso
\CC one (dR_{X/Q}(M)\to j_* dR_{W/Q}(N_W ))$, hence $ \det
Rq_{\text{dR}*} M\iso (\det Rq|_{W\text{dR}*} N_W )\otimes (\det
Rq|_{K_{0Q} *} \CF_- )^{\otimes -1}$. Let $N_X \subset j_* N_W$ be
$*$-extension from $V_+$ side and $!$-extension from $V_-$ side.
Then $i^* dR_{X/Q} (N_X )=\CF_+ \subset \CF_+ \oplus \CF_-$, hence
$\det Rq_{\text{dR}*}  N_X \iso (\det  Rq|_{W\text{dR}*} N_W
)\otimes (\det Rq|_{K_{0Q} *} \CF_- )^{\otimes -1}$. Thus we get a
canonical identification $\alpha :  \det  Rq_{\text{dR}*} M\iso \det
Rq_{\text{dR}*}N_X $. Now $N_W$, hence $N_X$, are $\CD$-modules,
i.e., they carry an absolute flat connection, so $
Rq_{\text{dR}*}N_X$, $Rq|_{W\text{dR}*} N_W$ carry a natural
connection. It is clear from the topology of the construction that
the cohomology are smooth, hence constant, $\CD_Q$-modules. Since
the fiber of $ Rq_{\text{dR}*}N_X$ at $q=0$ equals
$R\Gamma_{\!\text{dR}} (Y, N)$, we get (5.2.1).

For any $\CE\in {}^w \Bbb E$ we have an $\CO_Q$-line $\CE
(M)_{\nu_Q}:=  \CE (M)_{(D_Q ,1_T, \nu_Q )}$ and a natural
isomorphism $\CE (M)_{\nu_Q}\iso \CE (N)_{(D,1_T ,\nu
)}\otimes\CO_Q$, cf.~(2.13.4). There is a canonical isomorphism $\CE
(N)_{(D,1_T ,\nu )}\iso \CE (N)_\nu$ defined in the same way as
(2.13.9) (using $\eta$ on $\Bbb P^1$). Since $X$ is smooth over $Q^o
:= Q\setminus \{ 0\}$,  we have  $\eta : \CE (M)_{\nu_Q}|_{Q^o}\iso
\det  Rq_{\text{dR}*} M|_{Q^o}$. Thus comes a diagram
\begin{equation}
\begin{matrix}
 \CE (M )_{\nu_Q}|_{Q^o}& \buildrel{\eta}\over\lra
 &\det Rq_{\text{dR}*}M|_{Q^o}
  \\   \downarrow\,\, \,\,\,\,&&\,\,\,\,\,\,\,\,\, \downarrow    \\
 \CE(N)_{\nu }\otimes \CO_{Q^o}& \buildrel{\eta}\over\lra
 &\det R\Gamma_{\!\text{dR}}(Y,N)\otimes \CO_{Q^o}.
\end{matrix} \label{ 5.2.2}
\end{equation}

We say that $\CE$ is a {\it theory of $\varepsilon$-factors} if
(5.2.1) commutes for all data 2.13(a),(b). Such $\CE$ form a
subgroupoid $\Bbb E$ of $  {}^w \Bbb E$ called the
 {\it $\varepsilon$-gerbe}; see 5.4 for the reason.

  For $\CE^0 \in
 {}^w \Bbb E^0 $ there is a similar diagram

\begin{equation}
\begin{matrix}
 \CE^0 (M )_{\nu_Q}|_{Q^o}  \\  &\!\!\!\!\!\! \searrow  \\
  \downarrow&&\CO_{Q^o}\\
  &\!\!\!\!\!\!\nearrow\\
 \CE^0 (N)_{\nu }\otimes \CO_{Q^o}&
\end{matrix} \label{ 5.2.3}\end{equation}

Those $\CE^0$ for which (5.2.3) commutes for every datum 2.13(a),(b)
form a Picard subgroupoid $\Bbb E^0$ of  $ {}^w \Bbb E^0$. Our $\Bbb
E$ is an $\Bbb E^0$-torsor.

\bigbreak\noindent\textbf{Proposition.}\emph{   $\CE_{\text{dR}}$
and $ \CE_{\text{B}}$ are theories of
$\varepsilon$-factors.}\bigbreak

{\it Proof.}  Compatibility of  $\CE_{\text{dR}}$ with quadratic
degenerations follows from  2.13. Namely, the construction from
loc.~cit., spelled analytically as above, provides
$\eta_{\text{dR}}: \CE_{\text{dR}} (M)_{\nu_Q}\iso \det
Rq_{\text{dR}*} M$ over the whole $Q$ (not only on $Q^o$), and the
proposition in 2.13 says that our diagram commutes on the formal
neighborhood of $q=0$. Hence it commutes everywhere, q.e.d.

Let us treat $\CE_{\text{B}}$.  Let $K_\infty$ be a compact
neighborhood of  $T\cup |D|$ in $Y$ that does not intersect $U_\pm$,
$K:= K_\infty \cup \{ b_+ ,b_- \}$. Let $\CN$ be a continuous family
of proper cones in the tangent bundle to $Y\setminus K$ such that
Re$(\nu )$ is negative on it. Set $K_{\infty Q}:=  K_{\infty }\times
Q$, $K_Q := K_{\infty Q} \sqcup K_{0Q} \subset X$; let $\CN_W$ be
the pull-back of $\CN$ by the projection $W\to Y$. Then $(X , K_Q
,\CN_W ,\nu_Q )$ form a $Q$-family of compatible data as in 4.6.
Consider  isomorphisms $\eta$ of (4.7.1) for $M$ and $N_X$. Our
$Q$-family is constant near $K_{\infty Q}$, so
$\CE_{\text{B}}(M)_\CN (K_{\infty Q} )= \CE_{\text{B}}(N_X )_\CN
(K_{\infty Q} )= \CE_{\text{B}}(N)_\CN (K_{\infty})\otimes\CO_Q$.
Let $C$ be a locally closed subset of $V$ which consists of those $
(y_+ ,y_- )$ that $2 |t_+ (y_+ )|/r_+ - |t_- (y_- )|/r_- \le 1$ and
$2 |t_- (y_- )|/r_- - |t_+ (y_+ )|/r_+ < 1$; set
$Rq_{\text{dR}*}^{(C)}(?):= Rq|_{C *} Ri_C^! dR_{X/Q}(?)$. If $\CN$
is sufficiently tight, then $C\in C(X, \CN_W )$ (see the proof of
the lemma in 4.7), thus $\CE_{\text{B}}(? )_\CN (K_{\infty Q} )=\det
Rq_{\text{dR}*}^{(C)} (?)$.

Since the construction of (5.2.1) was local at $K_{0Q}$, the
 composition of $\CE_{\text{B}} (M)\buildrel{\eta_{\text{B}}}\over\lra
\det Rq_{\text{dR}*} M \iso \det R\Gamma_{\!\text{dR}}
(Y,N)\otimes\CO_Q =\det Rq_{\text{dR}*} N_X$ in (5.2.2) can be
rewritten as  $\CE_{\text{B}}(N)_\CN (K_{\infty})\otimes \det
Rq_{\text{dR}*}^{(C)}(M) \iso \CE_{\text{B}}(N)_\CN
(K_{\infty})\otimes \det  Rq_{\text{dR}*}^{(C)}(N_X
)\buildrel{\eta_{\text{B}}}\over\lra Rq_{\text{dR}*} N_X$. Here
$\iso$ comes from the identification \begin{equation}\alpha_C : \det
Rq_{\text{dR}*}^{(C)} (M) \iso \det Rq_{\text{dR}*}^{(C)}  (N_X
)\label{ 5.2.4}\end{equation} defined by the same construction as
(5.2.1) with
 $Rq_{\text{dR}*}$ replaced by $Rq_{\text{dR}*}^{(C)}$.

 The composition
 $\CE_{\text{B}} (M)\iso \CE_{\text{B}} (N)\otimes\CO_Q \iso \det R\Gamma_{\!\text{dR}} (Y,N)\otimes\CO_Q =\det
Rq_{\text{dR}*} N_X$ in (5.2.2) equals   $\CE_{\text{B}}(N)_\CN
(K_{\infty})\otimes \det Rq_{\text{dR}*}^{(C)} (M)\iso
\CE_{\text{B}}(N)_\CN (K_{\infty})\iso \CE_{\text{B}}(N)_\CN
(K_{\infty})\otimes \det Rq_{\text{dR}*}^{(C)}  (N_X
)\buildrel{\eta_{\text{B}}}\over\lra Rq_{\text{dR}*} N_X$. Here
$\iso $ come since  $Rq_{\text{dR}*}^{(C)} (M)=Rq_{\text{dR}*}^{(C)}
(N_X )=0$, hence their determinant lines are trivialized (notice
that the trivialization of $Rq_{\text{dR}*}^{(C)} (N_X )$ is
horizontal, and at $q=0$ it equals the Betti version of (2.13.8) due
to Exercise in 4.7).

We see that commutativity of (5.2.2) means that $\alpha_C$
identifies the above trivializations. To see this, consider the open
subspace $V_- \subset V$, and the corresponding 2-step filtration
$j_{V_- !} M|_{V_-} \subset M|_V$, and notice that $N_X |_V=\gr
M|_V$. The assertion follows now from the construction of
$\alpha_C$. \hfill$\square$

\subsection{}\label{5.3.}
Let $(X,T,M)$ be as in 5.1. For $\CE\in{}^w \Bbb E$, $b\in X$,  and
a meromorphic form $\nu$ on a neighborhood of $b$, $v_b (\nu
)=-\ell$, we write $\CE (M)_{(b,\nu )}:= \CE (M)_{(\ell b ,\nu )}$.

{\it Remark.} If $b$ is a smooth point of $M$, then  $\CE
(M)_{(b,\nu )}$ does not depend on whether we view $b$ as a point of
$T$ or $X\setminus T$ (by 5.1(iii) with $\pi = \id_X$,
$T=T'\cup\{b\}$).

Let   $\delta (M)_{b,\nu }\in\Bbb Z$ be  the degree of $\CE
(M)_{(b,\nu ) }$, and $\mu (M)_{b,\nu} \in\Bbb C^\times$ be the
value of $\mu \in\Aut (\CE (M))$ (see 1.15) at $(b,\nu )$.

\bigbreak\noindent\textbf{Lemma.}\emph{ (i) One has   $\delta
(M)_{b,\nu }=\dim (B(M)_b^{(!)}) + (1-\ell )\text{rk}(M)$.
\\
(ii) For  $c\in\Bbb C^\times$ the multiplication by $c$ automorphism
of $M$ acts on $\CE (M)_{(b,\nu )}$ as multiplication by
$c^{\delta(M)_{b,\nu}}$.
\\
(iii) One has $\mu (M)_{b,\nu} = (-1)^{ \ell\,\text{rk}(M)}m_b
(M)^{-1}$ where $m_b (M)$ is the monodromy  of $\det
M_{X\setminus T}$   around $b$.
\\
(iv) For smooth $M$, there is an isomorphism $\CE (M)^{(1)}\iso
(\det M)^{\otimes -1}$ compatible with constraint 5.1(ii) and
pull-backs by open embeddings. In particular, this is an isomorphism
of symmetric monoidal functors.
\\
(v) Suppose we have $M$, $M'$ over discs $U$, $U'$ which have
 regular singularity and are either $*$- or $!$-extension at $b$, $b'$.
Let $\phi : U\to U'$ be any open embedding, $\phi (b)=b'$, and
$\tilde{\phi} : M \to \phi^* M'$ be any its lifting. Then the
isomorphism $\CE (\tilde{\phi }): \CE (M)_b^{(1)}\iso  \CE
(M')_{b'}^{(1)}$ does not depend on the choice of $(\phi
,\tilde{\phi})$.}\bigbreak

{\it Proof.} (i) Let $t$ be a local coordinate at $b$, $t(b)=0$.  We
can view $t$ as an identification of a small disc $X_b$ at $b$ with
a neighborhood of $0\in\Bbb P^1$. Let us extend $M|_{X_b}$ to a
$\CD$-module $M^{(t)}$ on $\Bbb P^1$ which is smooth outside $\{ 0,
\infty \}$, and is the $*$-extension with regular singularities at
$\infty$. Such $M^{(t)}$ is unique.

By continuity,  $\delta (M)_{b,\nu }$ is the same for all $\nu$
with fixed $\ell$. We can assume that $\nu$ is meromorphic on $\Bbb
P^1$ with $\text{div} (\nu )\subset \{ 0,\infty\}$, so $v_0 (\nu
)=-2-v_\infty (\nu )$. By 5.1(iv), one has $\delta (
M^{(t)})_{0,\nu} +\delta (M^{(t)})_{\infty,\nu}=
\chi_{\text{dR}}(\Bbb P^1 ,M^{(t)})=\dim (B (M^{(t)})_b^{(!)})$.
Thus the assertion for $(X_b ,M,\nu )$ amounts to that for $(\Bbb
P^1_\infty ,M^{(t)}  ,\nu )$, i.e., we are reduced to the case when
$M$ is the $*$-extension with regular singularities. By 5.1(ii), it
suffices to treat the case of $\text{rk} (M)=1$; then, by
continuity, it suffices to consider $M=\CO_X$ (the trivial
$\CD$-module). By 5.1(iv) applied to $\Bbb P^1$ and $ t^{-1} dt$, we
see that $\delta (\CO_{\Bbb P^1} )_{0, t^{-1} dt} + \delta
(\CO_{\Bbb P^1} )_{\infty,t^{-1}dt}=-2$, hence, since
$v_{\infty}(t^{-1}dt )=v_0 (t^{-1}dt)$, one has $\delta (\CO_{\Bbb
P^1} )_{b,t^{-1}dt} = -1$. By factorization, $\delta (\CO_{\Bbb P^1}
)_{0, t^{-\ell} dt} = \ell \delta (\CO_{\Bbb P^1} )_{0, t^{-1}
dt}=-\ell$, q.e.d.

(ii) The $\Bbb C^\times$-action on $\CE (M)_{(b,\nu )}$, which comes
from the action of homotheties on $M$, is a holomorphic character of
$\Bbb C^\times$. Thus $c$ acts  as multiplication by $c^{\delta'
(M)_{b,\nu }}$ for some $\delta' (M)_{b,\nu }\in\Bbb Z$. The
argument of (i) works for $\delta$ replaced by $\delta'$, so
$\delta$ and $\delta'$ are given by the same formula, q.e.d.

(iv) By (i), $\CE^{(1)}(\CO_X )$ is a de Rham line of degree $-1$,
which has local origin. Thus there is a line $E$ of degree $-1$ and
an isomorphism $E\otimes \CO_X \iso \CE^{(1)}(\CO_X )$ compatible
with the pull-backs by open embeddings of $X$'s; such a datum is
uniquely defined.

The set of isomorphisms $\alpha : E\iso \Bbb C [1]$ identifies
naturally with the set of isomorphisms of symmetric monoidal
functors $\alpha_{\CE} : \CE (M)^{(1)}\iso (\det M)^{\otimes -1}$
(where $M$ is smooth) that are compatible with the pull-backs by
open embeddings. Namely, $\alpha_{\CE}$ is a unique isomorphism that
equals $\alpha \otimes\id_{\CO_X}$ for $M=\CO_X$. To see this,
notice that a matrix $g \in \text{GL}_n (\Bbb C)\iso \Aut (\CO_X^n
)$ acts on $\CE^{(1)}(\CO^n_X )$ as multiplication by $\det
(g)^{-1}$ (which follows from (ii) and 5.1(ii)).

(iii) Use (iv) and the compatibility property from 1.11.

(v) Let us show that $\Aut (M)$ acts trivially on $\CE (M)_b^{(1)}$.
Pick any $g\in\Aut (M)$. Since $M_{U\setminus \{b\}}$ admits a
$g$-invariant filtration with successive quotients of rank 1, we are
reduced, by 5.1(ii),  to the case of $M$ of rank 1. Here $g$ is
multiplication by some $c\in\Bbb C^\times$, and we are done by (ii)
(since, by (i), $\CE (M)_{b}^{(1)}$ has degree 0).

Thus for given $\phi$ the isomorphism $\CE (\tilde{\phi })$ does not
depend on the choice of $\tilde{\phi}$. The space of $\phi$'s is
connected, so it suffices to show that the map $\phi \mapsto \CE
(\tilde{\phi })$ is  locally constant. If $\phi$ varies in a disc
$Q$, then we can find $\tilde{\phi}$  which is an isomorphism of
$\CD$-modules on $U\times Q$, hence our map is horizontal (see
5.1(i)), q.e.d. \hfill$\square$

\subsection{}\label{5.4.}
For $\CE\in \Bbb E$ and $(X,T,M)$ as in 5.1 the canonical
automorphism $\mu$ of  $\CE (M)$ (see 1.15) is evidently compatible
with  constraints 5.1(i)--(iv), i.e., $\mu$ is an automorphism of
$\CE$. Here is the main result of this section:

\bigbreak\noindent\textbf{Theorem.} \emph{$ \Aut (\CE )$ is an
infinite cyclic group generated by $ \mu$. All objects of $\Bbb E$
are isomorphic. Thus $\Bbb E$ is a $\Bbb Z$-gerbe. }\bigbreak

We call  $\rho^\varepsilon : \CE_{\text{dR}}\iso
\CE_{\text{B}}$ an {\it $\varepsilon$-period} isomorphism. By the theorem,
$\varepsilon$-period isomorphisms form a  $\Bbb Z$-torsor $E_{\text{B}
/\text{dR}}$  referred to as   the  {\it  $\varepsilon$-period}
torsor.

Since $\Bbb E$ is an $\Bbb E^0$-torsor, the theorem can be
reformulated as follows. By (iii) of the lemma in 5.3, every $\CE^0
\in \Bbb E^0$ carries a natural automorphism $\mu^0$ that acts on
$\CE^0 (M)_{(b,\nu )}$ as multiplication by  $\mu (M)_{b,\nu
}\in\Bbb C^\times$.\footnote{ $\mu^0$ does {\it not} equal the
canonical automorphism $\mu$ of $\CE^0 (M)$ (which is identity by
5.5(i)).}

\bigbreak\noindent\textbf{Theorem$^\prime$.} \emph{The map $\Bbb Z
\to\pi_1 (\Bbb E^0 )$, $1\mapsto \mu^0$, is an isomorphism, and
 $\pi_0 (\Bbb E^0 )=0$. }\bigbreak

The proof  occupies the rest of the section.

\subsection{}\label{5.5.}  Pick any $\CE^0 \in {}^w \Bbb E^0$. Then for  $(X,T,M)$ as in 5.1 one has:

\bigbreak\noindent\textbf{Lemma.}\emph{ (i) The factorization line
$\CE^0 (M) \in \CL^{\Phi}_{\text{dR}}(X,T)$ is trivial.
\\
(ii) Every automorphism of $M$ acts trivially on $\CE^0 (M)$.
\\
(iii)  For $M$ of rang 0, the factorization line $\CE^0 (M)$ is
canonically trivialized. The trivialization has local nature and is
compatible with  constraints 5.1(i)$^0$--(iv)$^0$; it is uniquely
defined by this property.}\bigbreak

{\it Proof.} (i) One checks that the de Rham line $\CE^0 (M)^{(1)}$
on $X\setminus T$ is trivial by modifying the argument in the proof
of 5.3(iv) in the evident manner (or one can use  5.3(iv) directly,
noticing that $\CE^0$ is the ratio of two objects of ${}^w \Bbb E$).
 Similarly, $\CE (M)$ has zero degree by 5.3(i). Now use the theorem in 1.6.

(ii) Let us show that $g\in\text{Aut}(M)$ acts trivially on $\CE^0
(M)_{(b,\nu )}$. Let $M^{(t)}$ be as in the proof of 5.3(i); $g$
acts on it. We can assume that $\nu$ is meromorphic on $\Bbb P^1$
with div$(\nu)\subset\{ 0,\infty \}$.  The action of $g$ on $\CE^0
(M^{(t)}  )(\Bbb P^1 )=\CE^0 (M)_{(0,\nu )}\otimes\CE^0
(M^{(t)})_{(\infty ,\nu )}$  is trivial by 5.1(iv)$^0$, it suffices
to check that it acts trivially on $\CE^0 (M^{(t)})_{(\infty ,\nu
)}$. Thus we are reduced to the situation when our $M$ is the
$*$-extension with regular singularities. By constraint 5.1(ii)$^0$,
it suffices to consider the case when the monodromy of $M$ around
$b$ is  multiplication by a constant. Then $\Aut (M)$ is generated
by the diagonal matrices, and by 5.1(ii)$^0$ we are reduced to the
case when $M$ has rank 1, where we are done by  5.3(ii).

(iii) is left to the reader.  \hfill$\square$

{\it Remarks.} (i) By (i) of the lemma, the degree 0 lines $\CE^0
(M)_{(b,\nu )}$ for all $\nu$ with fixed $ v_b (\nu )=-\ell$ are
canonically identified; we denote this line by $\CE^0 (M)_{(b,\ell
)}$. By (ii) of loc.~cit., it depends only on the isomorphism class
of $M$, and by (iii) there is a canonical identification $\CE^0
(M)_{(b,\ell )}\iso \CE^0 (j_{b*} M)_{(b,\ell )}$.

(ii) Suppose we have $M$, $M'$ on discs $U$, $U'$ which have regular
singularity at $b^{(\prime )} \in U^{(\prime )}$. Let $\phi : U\to
U'$  be an open embedding, $\phi (b)=b'$, and $\tilde{\phi} : M \to
\phi^* M'$ be any its lifting. Then the isomorphisms $\CE^0
(\tilde{\phi }): \CE^0 (M)_{(b,\ell )}\iso  \CE^0 (M')_{(b' ,\ell
)}$ do not depend on the choice of $(\phi ,\tilde{\phi})$. To see
this, we can assume that $M$, $M'$ are $*$-extensions at $b$, $b'$,
and then repeat the second part of the proof of 5.3(v).

\subsection{}\label{5.6.} For $m\in\Bbb C^\times$ let $M_m$ be a $\CD$-module of
rank 1 on a disc $U$, which has regular singularity at $b\in U$ with
the monodromy   $m$ and is $*$-extension at $b$. By the remark in
5.5, the degree 0 line $\CG_{(m,\ell )}:=\CE^0 (M_m )_{(b,\ell )}$
depends only on $m$ and $\ell $. By 5.1(i)$^0$, $\CG_{(m,\ell )}$
form  a holomorphic line bundle $\CG =\CG (\CE^0 )$ over $\Bbb G_m
\times \Bbb Z$. The factorization structure on $\CE^0 (M_m )$
provides, by 5.5(i),  canonical isomorphisms
\begin{equation}\otimes\CG_{(1,\ell_\alpha )} \iso
\CG_{(1,\Sigma\ell_\alpha )}, \quad \CG_{(1,\ell
)}\otimes\CG_{(m,\ell' )} \iso \CG_{(m,\ell +\ell' )}. \label{
5.6.1}\end{equation} Suppose we have a finite collection $\{
(m_\alpha ,\ell_\alpha )\}$ with $\Pi m_\alpha =1$, $\Sigma
\ell_\alpha =2$. Then for any choice of a subset $\{ b_\alpha \}
\subset \Bbb P^1$ one can find a $\CD$-module $M$ on $\Bbb P^1$ of
rank 1 which is smooth off $\{ b_\alpha \}$ and is $*$-extension
with  regular singularity at $b_\alpha$ with monodromy $m_\alpha$,
and a rational form $\nu$ with div$(\nu )=-\Sigma \ell_\alpha
b_\alpha$. Writing $\CE (M)(\Bbb P^1 )= \CE (M)_\nu$ in 5.1(iv)$^0$,
we get \begin{equation}\eta : \otimes \CG_{(m_\alpha ,\ell_\alpha
)}\iso \Bbb C. \label{ 5.6.2}\end{equation} It does not depend on
the auxiliary choices (for $\eta$ is locally constant, and the datum
of $\{ b_\alpha \}$, $\nu$ forms a connected space; $M$ is unique up
to an isomorphism).

Now assume that $\CE^0 \in\Bbb E$. Consider the holomorphic $\Bbb
G_m$-torsor $\CG^\times = \CG^\times (\CE^0 )$ over $\Bbb G_m
\times\Bbb Z$ that corresponds to $\CG$.

\bigbreak\noindent\textbf{Lemma.}\emph{  $\CG^\times $ has a unique
structure of holomorphic commutative group $\Bbb G_m$-extension of
$\Bbb G_m \times \Bbb Z$ such that for any $ g_i \in \CG^\times_{
(m_i ,\ell_i )}$, $1\le i\le n$,  and $g\in \CG^\times_{ ((m_1
\ldots m_n )^{-1} ,2-\ell_1 -\ldots -\ell_n  )}$ one has $\eta (g
\otimes ( g_1 \cdots g_n ))=\eta (g \otimes g_1 \otimes\ldots
\otimes g_n )$.
 }\bigbreak

{\it Proof.} The above formula defines commutative $n$-fold product
maps $ \CG^{\times n} \to \CG^\times$, $(g_1 ,\ldots ,g_n )\mapsto
g_1 \cdots g_n$, which lift the $n$-fold products on $\Bbb G_m
\times\Bbb Z$. We need to check the associativity property, which
says that for any $g_1 ,\ldots ,g_n \in\CG^\times$ and any $k$,
$1<k< n$, one has $g_1 \cdots g_n =  (g_1 \cdots g_k ) g_{k+1}\cdots
g_n $.

For $m\in\Bbb C^\times$ set $\CG^\times_m := \sqcup
\CG^\times_{(m,\ell )} \subset \CG^\times$. The maps $(\CG_1^\times
)^n \to \CG_1^\times$, $\CG_1^\times \times\CG_m^\times \to
\CG_m^\times$ coming from the arrows in (5.6.1) are evidently
associative and commutative, i.e., they define a commutative group
structure on $\CG_1^\times$  and   a $\CG^\times_1$-torsor structure
on  $\CG^\times_m $. Thus $\CG^\times$ is a $\CG^\times_1$-torsor
over $\Bbb G_m$.

Since (5.6.2) comes from a trivialization of $\CE (M)$, the above
maps are restrictions of the multiple products maps in $\CG^\times$.
Moreover, the $n$-fold product  on $\CG^\times  $ is compatible with
the $\CG^\times_1$-action on $\CG^\times$: for $h\in  \CG^\times_1$
one has $(h g_1 )g_2 \cdots g_n =h (g_1 ,\ldots ,g_n )$. So, while
checking the associativity, we have a freedom to change $g_i$ in its
$ \CG^\times_1$-orbit. Thus we can assume that $g_i
\in\CG^\times_{(m_i,\ell_i )}$ are such that $\ell_1 +\ldots +\ell_k
= 1$. Then one can find a quadratic degeneration picture as in 5.2
such that  $\tilde{T}=  b_{1Q}\sqcup \ldots \sqcup b_{nQ}$,
$\tilde{M}$ of rank 1 has regular singularities at $b_{iQ}$ with
monodromy $m_i$, div$(\tilde{\nu})=-\Sigma \ell_i b_{iQ}$; the fiber
$\tilde{X}_1$ is $\Bbb P^1$, and $\tilde{X}_0$ is the union of two
copies of $\Bbb P^1$ with $\{ b_1 ,\ldots , b_k \}$ in the first
copy and $\{ b_{k+1},\ldots ,b_n \}$ in the second.  The
compatibility with quadratic degeneration yields  the promised
associativity, q.e.d. \hfill$\square$

\subsection{}\label{5.7.} Let $\pi^{(n)}: U' \to U$, $\pi^{(n)} (b')=b$, be a
degree $n$  covering of a disc completely ramified at $b$. Then
$\pi^{(n)}_* M_{m'} $ is isomorphic to $\mathop\oplus\limits_{
m^{n}=m'} M_{m}$, and $v_{b'} (\pi^{(n)*} \nu) +1= n (v_{b
}(\nu)+1)$. Therefore 5.1(iii)$^0$, 5.1(ii)$^0$ yield a canonical
isomorphism \begin{equation} \mathop\otimes\limits_{m^{
n}=m'}\CG_{(m,\ell )}\iso \CG_{(m', n(\ell -1)+1)}. \label{
5.7.1}\end{equation} For example, if $n=2$, $m'=1$, $\ell =1$, then
(5.6.2), with $\CG_{(1,1)}$ factored off, is a trivialization of $
\CG_{(-1,1)}$, which we denote by $e_{(-1,1)}\in\CG_{(-1,1)}$.
Notice that
\begin{equation}\eta (e_{(-1,1)}^{\otimes 2})=1 . \label{ 5.7.2}\end{equation} This
follows from compatibility of $\eta$ with $\pi_*$ applied to a
covering  $\Bbb P^1 \to\Bbb P^1$, $t\mapsto t^2$, the trivial
$\CD$-module $\CO_{\Bbb P^1}$ on the source, and the form
 $t^{-1}dt$ on the target.

Let  Ext$(\Bbb G_m ,\Bbb G_m )$ be the Picard groupoid of
holomorphic commutative group extensions of $\Bbb G_m $ by $\Bbb
G_m$. One has $\pi_0 (\text{Ext}(\Bbb G_m ,\Bbb G_m )=0$ and $\pi_1
(\text{Ext}(\Bbb G_m ,\Bbb G_m )=\Hom (\Bbb G_m ,\Bbb G_m )= \Bbb
Z$.

The quotient of $\Bbb G_m \times \Bbb Z$ modulo the subgroup
generated by $(-1 ,1)$ identifies  with $\Bbb G_m$ by the projection
$(m,\ell )\mapsto (-1)^\ell m$. Thus the quotient $\bar{\CG}^\times
(\CE^0 )$ of $\CG^\times (\CE^0 )$ modulo the subgroup generated by
$e_{(-1,1)}$ is an object of  Ext$(\Bbb G_m ,\Bbb G_m )$.
\begin{equation}\bar{\CG}^\times : \Bbb E^0 \to \text{Ext}(\Bbb G_m
,\Bbb G_m ) \label{ 5.7.3}\end{equation} is a Picard functor. It
assigns to $\mu^0 \in\Aut (\CE^0 )$ (see 5.4) the generator $-1$ of
$ \Bbb Z =\Aut (\bar{\CG}^\times (\CE^0 ))$. Therefore we can
reformulate the theorem from 5.4 as follows:

\bigbreak\noindent\textbf{Theorem.}\emph{  $\bar{\CG}^\times $ is an
equivalence of Picard groupoids. }\bigbreak

Let us define a Picard functor \begin{equation}\text{Ext}(\Bbb G_m
,\Bbb G_m )\to \Bbb E^0 \label{ 5.7.4}\end{equation} right inverse
to (5.7.3). We need to assign to an extension $\bar{\CG}^\times $ an
object $\CE^0 = \CE^0 (\bar{\CG}^\times )$ of $\Bbb E^0$. Suppose we
have $(X,T,M)$ as in 5.1. For $b\in T$ let $m(M)_b$ be the monodromy
of $\det M|_{X\setminus T}$ around $b$; for $c\in 2^T$ we denote by
$m(M)_c$ the product of $m(M)_b$ for $b\in T^c$  (see 1.1 for the
notation). Then
\begin{equation}
\CE^0 (M)_{(D,c,\nu )}:= \bar{\CG}_{(-1)^{\deg
(D)\text{rk}(M)}m(M)_c}. \label{ 5.7.5}
\end{equation}
Here $\bar{\CG}$ is the degree 0 line that corresponds to the $\Bbb
G_m$-torsor $\bar{\CG}^\times$. The factorization structure comes
from the product in $\bar{\CG}^\times$. Constraints 5.1(i)$^0$,
5.1(ii)$^0$ are evident. The identification $ \CE^0 (\pi_*
M')\iso\pi_* \CE (M') $ of 5.1(iii)$^0$ comes since both lines are
fibers of $\bar{\CG}$ over the same point of $\Bbb C^\times$. To see
this, it suffices to consider the situation of (5.7.1): there the
assertion is clear since $\mathop\Pi\limits_{m^n = m'} (-1)^\ell m =
(-1)^{n(\ell -1 )+1} m'$. Finally, for compact $X$ one has $m(M)_1
=1$ and $\deg (D)$ is even, hence $\CE^0 (M)_{(D,c,\nu
)}=\bar{\CG}_1=\Bbb C$, which is 5.1(iv)$^0$. The constraints are
mutually compatible by construction. We leave it to the reader to
check that $\CE^0$ is compatible with quadratic degenerations of
$X$, so we have defined (5.7.4). Due to an evident identification
$\bar{\CG}^\times (\CE^0 )=\bar{\CG}^\times$,  (5.7.3) is left
inverse to (5.7.4), so the theorem amounts to the next statement:

\bigbreak\noindent\textbf{Theorem$'$.}\emph{ For any $\CE^0 \in\Bbb
E^0$ there is a natural isomorphism $\iota: \CE^0 \iso \CE^0
(\bar{\CG}^\times (\CE^0 )). $ }\bigbreak

\subsection{}\label{5.8.}
The next step  reduces us to the setting of $\CD$-modules with
regular singularities. For a holonomic $\CD$-module $M$ we denote by
$M^{\text{rs}}$ the holonomic $\CD$-module with regular
singularities such that $B(M^{\text{rs}})=B(M)$, or, equivalently,
$M^\infty = M^{\text{rs}\infty}$, see 3.2. The functor $M\mapsto
M^{\text{rs}}$ sends nice families of $\CD$-modules to nice families
(as follows from 2.14), it is exact, comutes with $\pi_*$, and  one
has an evident identification $R\Gamma_{\!\text{dR}}(X,M)\iso
R\Gamma_{\!\text{dR}}(X,M^{\text{rs}})$. Thus for any theory of
$\varepsilon$-factors $\CE$ the rule $M\mapsto {}^r \!\CE (M):= \CE
(M^{\text{rs}})$ is again a theory of $\varepsilon$-factors. Clearly
$r$ is an endofunctor of $\Bbb E$. The same formula defines an
endofunctor $r^0$ of $\Bbb E^0$. It is naturally a Picard
endofunctor, and $r$ is a companion $\Bbb E^0$-torsor endofunctor:
one has ${}^{r^0}\!\CE^0 \otimes{}^{r^0} \!\CE^{\prime 0}\iso
{}^{r^0}\!(\CE^0 \otimes\CE^{\prime 0})$, ${}^{r^0}\!\CE^0
\otimes{}^r\!\CE\iso {}^r \!(\CE^0 \otimes\CE)$.

\bigbreak\noindent\textbf{Proposition.}\emph{ The endofunctor $r$ of
$\Bbb E$ is naturally isomorphic to $\id_{\Bbb E}$. Namely, there is
a unique  $\kappa : \id_{\Bbb E} \iso r$ such that $\kappa_{\CE ,M}
: \CE (M)\to \CE (M^{\text{rs}})$ is the identity map if $M$ has
regular singularities. Same is true for $\Bbb E$ replaced by $\Bbb
E^0$. Here $\kappa^0 : \id_{\Bbb E^0} \iso r^0 $ is an isomorphism
of Picard endofunctors, and $\kappa$  is an isomorphism of  the
companion $\Bbb E^0$-torsor endofunctors.
 }\bigbreak

{\it Proof.} (i) Let us define a canonical isomorphism of
factorization lines \begin{equation}\kappa =\kappa_{\CE ,M} :  \CE
(M)\iso \CE (M^{\text{rs}}).\label{ 5.8.1}\end{equation}

One has  $M|_{X\setminus T}= M^{\text{rs}}|_{X\setminus T}$, and
over $X\setminus T$ our $\kappa$ is the identity map.  It remains to
define $\kappa^{(1)} : \CE (M)^{(1)}_b \iso \CE
(M^{\text{rs}})^{(1)}_b$ for $b\in T$ (see 1.6 and  1.15).

Pick a local parameter $t$ at $b$. As in the proof of 5.3(i),
 $M$ yields a $\CD$-module
 $M^{(t)}$ on $\Bbb P^1$ with $\CE (M)^{(1)}_b = \CE (M^{(t)})_{(0,t^{-1}dt )}$. Ditto for
 $M^{\text{rs}}$. Since $M^{\text{rs}(t)}= M^{(t)\text{rs}}$ equals
$M^{(t)}$ outside 0,  constraints 5.1(iv) for $M^{(t)}$ and
$M^{\text{rs}(t)}$ yield isomorphisms $ \CE (M^{(t)})_{(0,t^{-1}dt
)}\otimes \CE (M^{(t)})_{(\infty ,t^{-1}dt )}\iso\det
R\Gamma_{\!\text{dR}} (\Bbb P^1 , M^{(t)})= \det R\Gamma_{\!\text{dR}}
(\Bbb P^1 , M^{\text{rs}(t)}) \buildrel\sim\over\leftarrow \CE
(M^{\text{rs}(t)})_{(0,t^{-1}dt )}\otimes \CE (M^{(t)})_{(\infty
,t^{-1}dt )}$. Factoring out $\CE (M^{(t)})_{(\infty ,t^{-1}dt )}$,
we get  $\kappa^{(1)}$.

  It remains to show that $\kappa^{(1)}$ does not depend on the auxiliary choice of $t$.
The space of local parameters $t$ is connected, so we need to check
that $\kappa$ is locally constant with respect to it. Let $t_s$ be a
family of local parameters at $b$ that are defined on the same disc
$X_b$ and depend holomorphically on $s\in S$; then $t_s$ identify
$X_{b  S}$ with a neighborhood $U$ of $\{ 0\} \times S$ in $\Bbb
P^1_S$. Let $M_U^{(t)}$ be the pull-back of $M$ by the projection
$U\to X_b$, $(v,s)\mapsto t_s^{-1}(v)$. This is a holonomic
$\CD_U$-module; let $M^{(t)}$ be a holonomic $\CD$-module on $\Bbb
P^1_S$ which equals $M_U^{(t)}$ on $U$, is smooth outside $\{
0,\infty\}\times S$, and is the $*$-extension with regular
singularities at $\{\infty\} \times S$. The restriction of $M^{(t)}$
to any fiber equals $M^{(t_s )}$, i.e., $M^{(t_s )}$ form a nice
isomonodromic family. The identifications $\CE (M)_b^{(1)}\iso  \CE
(M^{(t_s )})_{(0,t^{-1}dt )}$ are  horizontal, ditto for
$M^{\text{rs}}$. We are done the compatibility of $\eta$ with the
Gau\ss-Manin connection.

The same construction (with $R\Gamma_{\!\text{dR}}(\Bbb P^1, \cdot )$
replaced by $\Bbb C$) yields for $\CE^0 \in\Bbb E^0$ a canonical
isomorphism \begin{equation} \kappa^0 :  \CE^0 (M)\iso \CE^0
(M^{\text{rs}}).\label{ 5.8.2}\end{equation}

(ii)  $\kappa^0$ and $\kappa$ are evidently compatible with the
tensor product of $\CE^0$'s and $\CE$'s, and with constraints
5.1(i)--(iii) and 5.1(i)$^0$--(iii)$^0$. It remains to check
compatibility with constraint (iv). We treat the setting of $\Bbb
E^0$  (which suffices, say, since $r^2 =r$).  Suppose we have
$(X,T,M)$ with compact $X$. We need to prove that the composition
$\CE^0 (M)(X)\buildrel{\kappa^0}\over\to \CE^0
(M^{\text{rs}})(X)\buildrel{\eta (M^{\text{rs}})}\over\lra \Bbb C$
equals $\eta (M)$.

Let $b\in T$ be a point where $M$ has non-regular singularity, and
 $M^{\text{rs}_b}$ be the $\CD$-module which equals $M$ outside of $b$ and $M^{\text{rs}}$ near $b$. Let $\kappa^0_b :
 \CE^0 (M)\iso \CE^0 (M^{\text{rs}_b})$ be equal to $\kappa^0$ near $b$ and the identity morphism off $b$.

\bigbreak\noindent\textbf{Lemma.}\emph{ The composition $\CE^0
(M)(X)\buildrel{\kappa^0_b}\over\lra \CE^0
(M^{\text{rs}_b})(X)\buildrel{\eta (M^{\text{rs}_b})}\over\lra \Bbb
C$ equals $\eta (M)$.}\bigbreak

The lemma implies the proposition: since
$(M^{\text{rs}_b})^{\text{rs}}= M^{\text{rs}}$ and the composition
$\CE^0 (M) \buildrel{\kappa^0_b}\over\lra \CE^0 (M^{\text{rs}_b})
\buildrel{\kappa^0}\over\lra \CE (M^{\text{rs}})$ equals $\kappa^0$,
we are done by  induction by the number of points of $T$ where $M$
has non-regular singularity.

{\it Proof of Lemma.} Let  $\nu$ be a rational form  on $X$ such
that $\Res_b \,\nu  =1$. Let $t_b$ be a local parameter  at $b$ such
that $t_b^{-1}dt_b =\nu$.

Consider a datum 2.13(a) with $Y= \Bbb P^1 \sqcup X$, $b_+ = b\in
X$, $b_- =\infty \in \Bbb P^1$ and $\nu_Y$ equal  to $t^{-1}dt$ on
$\Bbb P^1$ and $\nu$ on $X$. Let $t_+$ be the parameter $t_b$, $t_-$
be the parameter $t^{-1}$ at $\infty$. The corresponding family of
curves $X'$ as defined in 5.2 (it was denoted by $X$ there) over
$Q=\Bbb A^1$ is the blow-up of $X\times \Bbb A^1$ at $(b,0)$. We
have a datum 2.13(b) with $N$ equal to $M^{(t_b )}$ on $\Bbb P^1
\setminus \{ \infty \}$ and to $M$ on $X\setminus \{ b\}$ (this
determines $N$ since it is !-extension with regular singularities at
$b_-$ and $*$-extension with regular singularities at $b_+$). Let
$L$ be any $t_\pm \partial_{t_\pm}$-invariant $b_\pm$-lattice in $N$
such that the eigenvalues of $\pm t_\pm \partial_{t_\pm}$ on
$L_{b_\pm}$ and their pairwise differences do not contain non-zero
integers.  Then the spectra of the $\pm t_\pm \partial_{t_\pm}$
actions on $L_{b_\pm}$ coincide, and there is a canonical
identification $\alpha : L_{b_+}\iso L_{b_-}$ characterized by the
next property: Consider the $t_\pm \partial_{t_\pm}$-invariant
embeddings $L_{b_-} \subset \Gamma (\Bbb P^1 \setminus \{ 0\}, L)$
and $L_{b_+} \subset \Gamma (U, L)$  as in 2.13(b). By the
definition of $M^{(t_b )}$, its sections over a punctured
neighborhood of 0 are identified with sections of $M$ over
$U\setminus \{ b\}$; by this identification the subspaces
$L_{b_\pm}$ correspond to one another, and $\alpha$ is the
corresponding isomorphism. Let $M'$ be the corresponding family of
$\CO_X$-modules with relative connection on $X/Q$ (which was denoted
by $M$ in 5.2).

At $q=1$ our $M'$ equals $M$, so the top  arrow in (5.2.3) at $q=1$
equals $\eta (M)$. And the composition of its  lower  arrows equals
the composition from the statement of our lemma. Since $\CE^0
\in\Bbb E^0$, the diagram commutes; we are done. \hfill$\square$

\subsection{}\label{5.9.}
Let us turn to the proof of Theorem$'$ in 5.7.  For $\CE^0 \in\Bbb
E^0$ let $\bar{\CG}^\times =\bar{\CG}^\times$ be the corresponding
extension of $\Bbb G_m$ (see (5.7.3)) and $\CE^{0\prime}\in \Bbb
E^0$ be the object defined by  $\bar{\CG}^\times$ (see (5.7.4)). We
want to define a natural isomorphism $\iota
:\CE^0\iso\CE^{0\prime}$.

For $(X,T,M)$ as in 5.1 we define a canonical isomorphism of
factorization lines \begin{equation}\iota : \CE^0 (M)\iso
\CE^{0\prime}(M)  \label{ 5.9.1}\end{equation}  as follows. Due to
isomorphism (5.8.2), we can assume that $M=M^{\text{rs}}$. Now
$\iota$ is a unique isomorphism of local nature which is compatible
with constraints 5.1(ii)$^0$ and 5.5(iii), and is the identity map
for $M= M_m$ (see 5.6). Indeed, we can assume that $X$ is a disc,
$T=\{ b\}$, and, by 5.5(iii), that  $M$ is $*$-extension at $b$.
Then $M = \oplus M^{(m)}$, where  the monodromy around $b$ acts on
$M^{(m)}$ with eigenvalues $m$. By compatibility with 5.1(ii)$^0$,
we can assume that $M=M^{(m)}$. Pick any filtration on $M$ with
successive quotients of rank 1 and define $\iota$ as the composition
$\CE^0 (M) \iso \otimes\CE^0 (\gr_i M)= \otimes\CE^{0\prime} (\gr_i
M)\iso \CE^{0\prime} (M)$ where $\iso$ are constraints 5.1(ii)$^0$.
The choice of the filtration is irrelevant by 5.1$^0$ (for the space
of filtrations is connected).

Our $\iota$ is  compatible with constraints 5.1(i)$^0$, 5.1(ii)$^0$;
its compatibility with  5.1(iii)$^0$ will be checked   in 5.13. We
 treat  5.1(iv)$^0$ first; this takes 5.10--5.12. For $(X,T,M)$ with
compact $X$ let $\xi (X,M)\in\Bbb C^\times$ be the ratio of $\eta
(M)$ for $\CE^0$ and the composition of $\eta (M)$ for
$\CE^{0\prime}$ with $\iota$.  We want to show that $\xi (X,M)\equiv
1$.

\subsection{}\label{5.10.} By 5.8 and the  construction of $\iota$, it suffices to
consider $M$ with regular singularities. By 5.5(iii), $\xi (X,M)$
depends only on $M|_{X\setminus T}$. Therefore, by compatibility
with 5.1$^0$,  $\xi (X,M)$ depends only on the purely topological
datum of (the isomorphism class of) a punctured oriented surface
$X\setminus T$ (that can be replaced by a compact surface with
boundary) and a local system on it.

The compatibility with quadratic degenerations implies that if $Y$
is obtained from $X$ by cutting along a disjoint union of embedded
circles and $N= M|_Y$, then $\xi (X,M)=\xi (Y,N)$. Here is an
application:

\bigbreak\noindent\textbf{Lemma.}\emph{  (i) If $M$ admits a
filtration such that $\gr_i M$ are $\CD$-modules of rank 1, then
$\xi (X,M)=1$.
\\
(ii)  For every $(X,T, M)$ with $X$ connected one can find $(X',T,
M')$ such that the restriction of $M$ to a neighborhood of $T$ is
isomorphic to that of $M'$, $X'$ is connected of any given genus $g
\ge g(X)$, and  $\xi (X,M)=\xi (X',M')$.
\\
(iii) For every $(X,T, M)$ one can find $(X',T, M')$ with $X'$
connected such that $\xi (X,M)=\xi (X',M')$ and for every $b\in T$
the restriction of $M'$ to a neighborhood of $b$ is isomorphic to
that of $M$ plus a direct sum of copies of a trivial
$\CD$-module.}\bigbreak

{\it Proof.} (i) By compatibility with 5.1(ii)$^0$, we can assume
that $M$ is a $\CD$-module of rank 1. Our assertion
 is true if  $X$ has genus 0 by the construction. An arbitrary $X$ can be cut into a union of genus 0 surfaces, and we are done.

(ii) Consider $Y= X\sqcup Z$ where $Z$ is a compact smooth connected
curve of genus $g-g(X)$; let $N$ be a $\CD_Y$-module such that $N|_X
=M$ and $N|_Z$ is a trivial $\CD$-module of the same rank  as $M$.
Pick $x\in X\setminus T$, $z\in Z$, cut off small discs around $x$,
$z$ and connect their boundaries by a tube. This is $X'$. Take for
$M'$ any extension of $M_Y$ (restricted to the complement of the cut
discs) to a local system on $X'$. Since $\xi (Z,N|_Z )=1$ by (i),
one has $\xi (Y,N)=\xi (X,M)$, hence $\xi (X',M')=\xi (X,M)$.

(iii) Let us construct $(X',M')$. First, add to $M$ on different
components of $X$ appropriate number of copies of the trivial
$\CD$-module to assure that the rank of $M$ is constant; this does
not  change  $\xi (X,M)$ by (i). Let $X_1 ,\ldots ,X_n$ be the
connected components of $X$. On each $X_i \setminus T$, choose a
pair of distinct points $x_i$, $y_i$. Cut off small discs around
$x_1 ,\ldots ,x_{n-1}$ and $y_2 ,\ldots ,y_n$, and connect the
boundary circle at $x_i$ with that at $y_{i+1}$ by a tube. This is
our $X'$. Take for $M'$ any extension of $M$ to a local system on
$X'$.
  \hfill$\square$

\subsection{}\label{5.11}

\bigbreak\noindent\textbf{Proposition.}\emph{ For $X$ connected,
$\xi (X,M)$ depends only on the datum of conjugacy classes of local
monodromies of $M$ (the rank of $M$ is fixed). }\bigbreak

{\it Proof.} According to \cite{PX1}, \cite{PX2}, the action of the
mapping class group on the moduli of unitary local systems of given
rank with fixed cojugacy classes of local monodromies, is ergodic
(provided that the genus of the Riemann surface is $>1$). As in
Theorem 1.4.1 in \cite{G}, this implies that for connected $X$ with
$g(X)>1$ our $\xi (X,M)$ depends only on $g(X)$, the rank of $M$,
and the datum of conjugacy classes of local monodromies of $M$
(indeed, $\xi$ is invariant with respect to the action of the
mapping class group by the compatibility with 5.1$^0$, and is
holomorphic; by the ergodicity, its restriction to the real points
of the moduli space of local systems is constant, and we are done).
Use 5.10(ii)  to eliminate the  dependence on $g(X)$ (and the
condition on $g(X)$). \hfill$\square$

\subsection{}\label{5.12.}
For any $(X,M)$, let $Sp (M)$ be the datum of other than 1
eigenvalues (with multiplicity) of the direct sum of local
monodromies. We write it as an element $\Sigma n_i z_i$ ($z_i$ are
the eigenvalues, $n_i$ are the multiplicities) of the quotient  of $
Div (\Bbb C^\times )$ modulo the subgroup of divisors supported at
$1\in\Bbb C^\times$.

\bigbreak\noindent\textbf{Lemma.}\emph{ $\xi (X,M)$ depends only on
$Sp(M)$. }\bigbreak

{\it Proof.} (i) By  5.10(iii), it suffices to check this assuming
that $X$ is connected, and by 5.10(i) we can assume that the rank of
$M$ is fixed. By 5.11,   it suffices to find for any  $(X,M)$ some
$(X' ,M')$ such that $\xi (X,M)=\xi (X',M')$, $Sp(M)=Sp(M')$, and
each local monodromy of $M'$ has at most one eigenvalue different
from 1. Take any $b \in T$; let $m_{b}$ be the local monodromy at
$b$.  Then one can find a  local system $K_{(b)}$ on $\Bbb P^1$ with
ramification at $\infty$ and $n$ other points, $n= \text{rk}( M)$,
such that its local monodromy at $\infty$ is conjugate to
$m_{b}^{-1}$, and $K_{(b)}$ admits a flag of local subsystems such
that each $\gr_i K_{(b)}$ has rank 1 and ramifies at $\infty$ and
only one other point. Cut a small disc around $b$ in $X$ and that
around $\infty$ in $\Bbb P^1$, and connect the two boundary circles
by a tube; we get a surface  $X'_{(b)}$. Let $M'_{(b)}$ be a local
system on it that extends $M$ and $K_{(b)}$. By  5.10(i), $\xi (\Bbb
P^1 ,K_{(b)})=1$, so  $\xi (X'_{(b)},M'_{(b)})=\xi (X,M)$ by the
compatibility with quadratic degenerations.  Repeating this
construction for each point of $T$, we get $(X',M')$.
\hfill$\square$

The lemma  implies that $\xi (X,M)=1$. Indeed, if $Sp (M)=\Sigma n_i
z_i$, then $\Pi z_i^{n_i}=1$ (for the product does not change if we
replace $M$ by $\det M$, where it equals 1 by the Stokes formula).
Therefore one can find a $\CD$-module $M'$ of rank 1 on $\Bbb P^1$
with $Sp (M')=Sp(M)$. Since $\xi (X,M')=1$, we are done by the
lemma.

\subsection{}\label{5.13.} It remains to check that  $\iota$ of (5.9.1) is
compatible with 5.1(iii)$^0$. We want to show that for $\pi : X'\to
X$ and a $\CD$-module $M'$ the diagram

\begin{equation}
\begin{matrix}
 \CE^0 (\pi_* M' )& \buildrel{\iota}\over\lra
 &\CE^{0\prime} (\pi_* M' )
  \\   \downarrow\,\, \,\,\,\,&&\,\,\,\,\,\,\,\,\, \downarrow    \\
 \pi_* \CE^0 (M')& \buildrel{\pi_* \iota}\over\lra
 &  \pi_* \CE^{0\prime} (M'),
\end{matrix} \label{ 5.13.1}
\end{equation}
where the vertical arrows are constraints 5.1(iii)$^0$ for $\CE^0$,
$\CE^{0\prime}$, commutes. For $b\in X$ let $\psi (M',\pi ,b)$ be
the ratio of the morphisms $\CE^0 (\pi_* M' )^{(1)}_b \to \pi_*
\CE^{0\prime} (M')^{(1)}_b$ that come from the two sides of the
diagram. We want to show that $\psi (M',\pi ,b) \equiv 1$ (see 1.6).
It is clear that $\psi (M',\pi ,b)=1$ if $\pi$ is unramified at $b$.

Our $\psi$ has $X$-local nature and it is multiplicative with
respect to disjoint unions of $X'$, so it suffices to consider the
case when $X$, $X'$ are discs and $\pi =\pi^{(n)}$ is ramified of
index $n$ at $b$. Choosing a local coordinate $t$ at $b$, we
identify $X'$ and $X$ with neighborhoods of 0 in $\Bbb P^1$ so that
our covering is the restriction of $\pi :  \Bbb P^1 \to  \Bbb P^1$,
$t\mapsto t^n$, to $X$. Let us extend $M'$ to a $\CD$-module on
$\Bbb P^1$, which we again denote by $M'$, such that it is smooth
outside $0$ and $1$. We know that  5.1(iii)$^0$ is compatible with
5.1(iv)$^0$ for both $\CE^0$ and $\CE^{0\prime}$. Since $\iota$ is
compatible with 5.1(iv)$^0$ and $\pi$ is ramified only at 0 and
$\infty$, we know that $\psi (M',\pi, 0) \psi (M',\pi, \infty )=1$.
Since $M'$ is smooth over $\infty$, this means that $\psi (M',
\pi^{(n)},b)\psi ( \CO_{X'},\pi^{(n)},b)^{\text{rk}(M)}=1$. We are
reduced to the case  $M' =\CO_{X'}$.

Set $\psi_n :=\psi (\CO_{X'}, \pi^{(n)},b)$. By above, $\psi_n^2
=1$, i.e., $\psi_n = \pm 1$. By the construction of $\CE^{0\prime}$
(see 5.7), one has $\psi_2 =1$. Due to compatibility with the
composition, one has $\psi_{mn}=\psi_n \psi_m^n = \psi_m \psi_n^m$,
i.e, $\psi_n^{m-1}=\psi_m^{n-1}$. For $m=2$ we get $\psi_n \equiv
1$, q.e.d. \hfill$\square$

\section{The  $\Gamma$-function.}

\subsection{}\label{6.1.}
Let us describe explicitly the $\varepsilon$-period map
$\rho^\varepsilon = \rho^\varepsilon (M):  \CE_{\text{dR}}(M)\iso
\CE_{\text{B}} (M)$ for   $\rho^\varepsilon \in  E_{\text{B}
/\text{dR}}$ (see 5.4).

By 5.8, $\rho^\varepsilon (M)$ equals the composition $
\CE_{\text{dR}}(M) \buildrel{\kappa}\over\to \CE_{\text{dR}} (
M^{\text{rs}})\iso \CE_{\text{B}} (M^{\text{rs}}) = \CE_{\text{B}}
(M)$, where $\kappa$ is
 the canonical isomorphism of (5.8.1) and $\iso$ is $\rho^\varepsilon ( M^{\text{rs}})$.
A different construction of the same $\kappa$ was given in (3.1.1)
in terms of certain analytical Fredholm determinant (a version of
$\tau$-function).

 From now on we assume that $M$ has regular singularities.

By (1.6.3), Example in 1.6, and Remark in 1.15,  $\CE_? (M)$ amounts
to a datum $(\CE_? (M)^{(1)}_{X\setminus T}, \{\CE_? (M)^{(1)}_b
\})$.
 Thus $\rho^\varepsilon $ is completely determined by the
isomorphisms $ \CE_{\text{dR}}(M)^{(1)}_{X\setminus T}\iso
\CE_{\text{B}} (M)_{X\setminus T}^{(1)}$ and  $
\CE_{\text{dR}}(M)^{(1)}_b \iso \CE_{\text{B}}(M)^{(1)}_b$, $b\in
T$.

Let us write a formula for $\rho^\varepsilon =\rho^\varepsilon_b :
\CE_{\text{dR}}(M)^{(1)}_b \iso \CE_{\text{B}}(M)^{(1)}_b$,  $b\in
T$. Below  $t$ is a local parameter at $b$, and $i_b$, $j_b$ are the
embeddings $\{b\}\hra X_b \hookleftarrow X_b^o :=  X_b \setminus \{
b\}$.

If $M$ is supported at $b$, then $\CE_{\text{dR}}(M)^{(1)}_b
=\CE_{\text{B}}(M)^{(1)}_b =\det R\Gamma_{\!\text{dR}}(X,M)$ and
$\rho^\varepsilon_b$ is the identity map. Thus for arbitrary $M$ one
has $\CE_? (M)^{(1)}_b = \CE_? (j_{b*}M)^{(1)}_b \otimes \det
R\Gamma_{\! \text{dR}\,\{ b\}}(X,M)$ and
\begin{equation}\rho^\varepsilon_b (M)= \rho^\varepsilon_b (j_{b*}
M)\otimes \id_{\det R\Gamma_{\! \text{dR}\,\{b\}}(X,M)}.  \label{
6.1.1}\end{equation} So it suffices to define $\rho^\varepsilon_b$
for $M=j_{b*} M$. Then we have a canonical trivialization $1^!_{b}$
of  $\CE_B (M)^{(1)}_b := \CE_B (M)_{(b,t^{-1}dt)}$, see 4.7.

Let $L$  be a $t\partial_t$-invariant $b$-lattice in $M$. Denote by
$\Lambda (L)$ the spectrum of the  operator $t\partial_t$ acting on
on $L_b =L/t L$. Suppose that it does not contain  positive
integers.  Then the complex $\CC (L,\omega L(b))$ (see 2.4) is
acyclic; here $L(b):= t^{-1}L$, i.e., $\omega L(b)= t^{-1}dtL$.
Denote by $\iota (L)_{zt^{-1}dt}$ the corresponding trivialization
of  $ \CE_{\text{dR}}(M)_{(b, zt^{-1}dt)} \iso  \det \CC  (L,\omega
L(b))$, $z\neq 0$ (see (2.5.6)); it  does not depend on the choice
of $t$. If $L' \supset L$ is another  lattice, then $\iota
(L)_{zt^{-1}dt}/\iota (L')_{zt^{-1}dt} $ is the determinant of the
action of $z^{-1}t\partial_t $ on $L'/L$ (see 2.5). In particular,
we have a trivialization $\iota (L)_{t^{-1}dt}$ of $
\CE_{\text{dR}}(M)^{(1)}_b$.

We write $\rho^\varepsilon_b (\iota (L)_{t^{-1}dt})=
\gamma^!_{\rho^\varepsilon} (L)1^!_{ b}$. For  example, if
$M=M^\lambda_t$ is the $\CD$-module $M^\lambda_t$ generated by
$t^\lambda$, $t\partial_t (t^\lambda )=\lambda t^\lambda$, where
$\lambda \in\Bbb C\setminus \Bbb Z_{>0}$, and $L=L^\lambda_t$ is the
lattice generated by $t^\lambda$, then $\Lambda (L^\lambda_t )=\{
\lambda\}$, and we write $\gamma^!_{\rho^\varepsilon}(\lambda ):=
\gamma^!_{\rho^\varepsilon}(L^\lambda_t)$.

\bigbreak\noindent\textbf{Theorem.}\emph{ (i) One has
$\gamma^!_{\rho^\varepsilon} (L)= \mathop\Pi\limits_{\lambda  \in
\Lambda (L)}\gamma^!_{\rho^\varepsilon}(\lambda)$.
\\
(ii) For $a \in\Bbb Z$ one has\footnote{Recall that
$E_{\text{B}/\text{dR}}$ is a $\Bbb Z$-torsor.}
$\gamma^!_{\rho^\varepsilon +a}(\lambda)= (-1)^a \exp (-2\pi
i\lambda a)\gamma^!_{\rho^\varepsilon}(\lambda)$.
\\
(iii) For one $\rho^\varepsilon$ in $E_{\text{B}/\text{dR}}$   one
has
\begin{equation}\gamma^!_{\rho^\varepsilon} (\lambda)=
(2\pi )^{-1/2} (1-\exp (2\pi i \lambda ))\Gamma(\lambda)  , \label{
6.1.2}\end{equation} where  $\Gamma$ is the Euler $\Gamma$-function
and $(2\pi)^{1/2}$ is the positive square root. }\bigbreak

 {\it Proof.} (i) By above, for $L'\supset L$
one has  $\gamma^!_{\rho^\varepsilon}
(L)/\gamma^!_{\rho^\varepsilon}  (L')=\det (t\partial_t ; L'/L)$.
Therefore the validity of (i) does not depend on the choice of $L$.
So we can assume that $(M,L)$ is a successive extension of some
$(M^\lambda_t ,L^\lambda_t )$, and we are done since all our objects
are multiplicative with respect to extensions.

(ii) By  5.3(iii), $\mu (M^\lambda_t )$ acts on $\CE (M^\lambda_t
)^{(1)}_b$ as multiplication by   $-\exp (-2\pi i  \lambda )$.

 (iii) The claim follows  from (ii) and the next lemma:

\bigbreak\noindent\textbf{Lemma.}\emph{ (i) The function  $
\gamma^!_{\rho^\varepsilon}(\lambda )$  is holomorphic and
invertible for $\lambda\in \Bbb C \setminus \Bbb Z_{>0}$, and
satisfies the next relations: (a)
$\gamma^!_{\rho^\varepsilon}(\lambda +1)= \lambda
\gamma^!_{\rho^\varepsilon}(\lambda )$; (b) For every positive
integer $n$ one has $\gamma^!_{\rho^\varepsilon}(\frac{\lambda}{ n}
) \gamma^!_{\rho^\varepsilon}(\frac{\lambda +1}{n})\cdots
\gamma^!_{\rho^\varepsilon} (\frac{\lambda
+n-1}{n})=n^{\frac{1}{2}-\lambda}\gamma^!_{\rho^\varepsilon}(\lambda
)$.
\\ (ii) Any function $\gamma$ that satisfies the properties from (i)  equals one
of the functions $\gamma_a (\lambda )= (2\pi )^{-1/2} (-1)^a
\exp(2\pi i \lambda a) (1-\exp (2\pi i\lambda ) )\Gamma (\lambda )$
for some integer $a$. }\bigbreak

{\it Proof of Lemma.} (i) We check (b); the rest is clear. Let $\pi
: X' \to X$ is a covering of a disc completely ramified of index $n$
at $b$, so for a parameter $t'$ at $b'$ one has $\pi^* (t)=t^{\prime
n}$, hence $\pi^* ( t^{-1}dt)=n t^{\prime -1}dt'$. Let $M'$ be a
$\CD$-module on $X'$ which is the $*$-extension with regular
singularities at $b'$. Consider  isomorphisms \begin{equation}
\CE_{\text{dR}}( \pi_* M' )_{(x,t^{-1}dt)}\buildrel\alpha\over\lra
\CE_{\text{dR}}(M')_{(x',m t^{\prime -1}dt')}
\buildrel\beta\over\lra \CE_{\text{dR}}(M')_{(x',t^{\prime -1}dt')}
\label{ 6.1.4}\end{equation} where $\alpha$ is the projection
formula identification and $\beta$ is  the
$\nabla^\varepsilon$-parallel transport  along the interval
$[m,1]t^{\prime -1}dt'$. By the construction of $1_b^!$, the Betti
version of $\beta\alpha$ transforms $1^!_b$ to $1^!_{b'}$.

Suppose $M'$ equals $M^{\lambda}_{t'}$ for some $\lambda\in\Bbb C$.
Then $\pi_* M'$ equals $M_{t}^{\lambda/n }\oplus M_{t}^{(\lambda
+1)/n}\oplus\ldots\oplus M_{t}^{(\lambda+n-1)/n}$. If
$L'=L_{t'}^{\lambda} \subset M_{t'}^{\lambda}$, then $\pi_* L'  =
L_{t}^{\lambda/n }\oplus L_{t}^{(\lambda +1)/n}\oplus\ldots\oplus
L_{t}^{(\lambda +(n -1))/n}$. It is clear that $\alpha$ sends $\iota
(\pi_* L')_{t^{-1}dt}$ to $\iota (L')_{n t^{\prime -1}dt'}$. Since
$\iota (L')$ is horizontal for the connection $\nabla_0$ of (2.11.3)
(with $\ell $ and $n$ in loc.~cit.~ equal to 1),  $\beta$ sends
$\iota (L')_{nt^{\prime -1} dt'}$ to $n^{\frac{1}{2} -\lambda }\iota
(L')_{t^{\prime -1}dt'}$. Now $\rho^\varepsilon (\iota
(L')_{t^{\prime -1}dt'})= \gamma^!_{\rho^\varepsilon}(\lambda
)1^!_{b'}$ and $\rho^\varepsilon (\iota (\pi_* L')_{t^{-1}dt})=
\gamma^!_{\rho^\varepsilon}(\pi_* L')1^!_b
=\gamma^!_{\rho^\varepsilon}(\frac{\lambda}{ n} )
\gamma^!_{\rho^\varepsilon}(\frac{\lambda +1}{n})\cdots
\gamma^!_{\rho^\varepsilon}(\frac{\lambda +n-1}{n})1^!_b$, and we
are done since $\rho^\varepsilon$ is compatible with  5.1(iii).

(ii) Denote by $E$ the set of functions $\gamma$ that satisfy
properties from (i). Let $E^0$ be the set of functions $e(\lambda )$
which are invertible and holomorphic on the whole $\Bbb C$ and
satisfy the relations (a) $e(\lambda +1)=e(\lambda )$; (b)
$e(\frac{\lambda}{ n} ) e(\frac{\lambda +1}{n})\cdots e
(\frac{\lambda +n-1}{n})=e(\lambda )$ for any positive integer $n$.
Then $E^0$ is a group with respect to multiplication, and $E$ is an
$E^0$-torsor.

Notice that the function $ \mu (\lambda ):= -\exp (-2\pi i\lambda )$
belongs to $E^0$, and $\{ \lambda_a \}_{a\in\Bbb Z}$ is a $\mu^{\Bbb
Z}$-torsor.  Recall that $\Gamma (\lambda)$ is holomorphic and
invertible for $\lambda \in\Bbb C \setminus \Bbb Z_{\le 0}$, and
satisfies the next relations:  (a) $\Gamma (\lambda+1
)=\lambda\Gamma (\lambda)$; (b)  $\Gamma (\frac{\lambda }{n})\ldots
\Gamma ( \frac{\lambda+n-1}{n})= (2\pi )^{\frac{n-1}{2}}
n^{\frac{1}{2} - \lambda}\Gamma (\lambda)$ for any positive integer
$n$. This implies that $ \gamma_a$ belong to $E$. To prove the
lemma, it remains to check that $\mu$ generates $E^0$.

Pick any $e\in E^0$; let $a$ be the index of the holomorphic map $e:
\Bbb C /\Bbb Z \to \Bbb C^\times$. Let us check that $e\mu^a \equiv
1$. Indeed, $e\mu^a$ has index 0, so $e\mu^a (\lambda )=\exp
f(\lambda )$ for some holomorphic $f : \Bbb C /\Bbb Z \to \Bbb C$.
Notice that for any $n\in\Bbb Z_{>0}$ the function  $\lambda \mapsto
f(\frac{\lambda}{ n} ) +f (\frac{\lambda +1}{n})+\ldots +
f(\frac{\lambda +n-1}{n})- f(\lambda )$ takes values in $2\pi i \Bbb
Z$, hence constant. Consider the coefficients of the Laurent series
$f(\lambda )= \Sigma b_m \exp (2\pi i m \lambda )$. The above
property implies that $(n-1)b_{\pm n}=0$ for $n> 1$, i.e., $b_m =0$
for $|m|>1$. The case $n=2$ shows that $b_{\pm 1}=0$. Finally the
fact that $(n-1)b_0 \in 2\pi \Bbb Z$ for any $n>0$ implies that $b_0
\in 2\pi i \Bbb Z$, and we are done. \hfill$\square$

\bigbreak\noindent\textbf{Corollary.}\emph{ For $\rho^\varepsilon$
as in (6.1.2) the isomorphism  $\rho^\varepsilon :
\CE_{\text{dR}}(M)_{X\setminus T}^{(1)}\iso \CE_{\text{B}}
(M)_{X\setminus T}^{(1)}$ equals the composition $\CE_{\text{dR}}
(M)^{(1)}_{X\setminus T}\buildrel{(2.6.1)}\over\lra  (\det
M_{X\setminus T})^{\otimes -1}  \buildrel{(4.7.3)}\over\lra
\CE_{\text{B}} (M)^{(1)}_{X\setminus T}$ multiplied by $((2\pi
)^{1/2} i )^{\text{rk}(M)}$. Replacing $\rho^\varepsilon$ by
$\rho^\varepsilon +a$ multiplies it by $(-1)^{\text{rk}
(M)a}$.}\bigbreak

{\it Proof.} Suppose $M$ is smooth at $b$. Compatibility with
5.1(iii) implies, as in Remark in 5.3, that  $\rho^\varepsilon_b :
\CE_{\text{dR}}(M)^{(1)}_b \iso \CE_{\text{B}} (M)_b^{(1)}$ does not
depend on whether $b$ is viewed as a point of $T$ or not. The exact
sequence $0\to M\to j_{b*} j_b^* M \to i_{b*} M_b \to 0$ shows that
$R\Gamma_{ \text{dR}\, b} (X,M)=M_b [-1]$, hence  $\CE_? (M)^{(1)}_b
=  \CE_? (j_{b*}  M )^{(1)}\otimes (\det M_b )^{\otimes -1}$. The
isomorphisms $\CE_? (M)^{(1)}_b \iso (\det M_b )^{\otimes -1}$ come
from  trivializations of $ \CE_? (j_{b*}  M )^{(1)}$, which are
$\iota (M)_{t^{-1}dt}$ in the de Rham and $1^!_{b}$ in the Betti
case. Since $\gamma^!_{\rho^\varepsilon} (M)=
((2\pi)^{1/2}i)^{\text{rk} (M)}$ by the theorem, we are done.
\hfill$\square$

The corollary together with the theorem completely determines
$\rho^\varepsilon (M)$.

\subsection{}\label{6.2.}
Here is another explicit formula for   $$\rho^\varepsilon :
\CE_{\text{dR}}(M)_{(b, -t^{-1}dt )}\iso \CE_{\text{B}} (M)_{(b,
-t^{-1}dt )}.$$

Recall that $\CE_{\text{B}} (M)_{(b,-t^{-1}dt)}\iso \det
R\Gamma_{\!\text{dR}}(X_b, M)$ by (4.7.1); here $X_b$ is a small open
disc at $b$. Let $L$ be a $t\partial_t$-invariant $b$-lattice in  $
M$ such that $\Lambda (L)$ does not contain non-positive integers,
$L_\omega$ be the $\CO$-submodule of $\omega M$ generated by $\nabla
(L)$ (this is a $b$-lattice). Then the projection
\begin{equation}\Gamma (X_b ,dR (M)) \twoheadrightarrow \CC (L,
L_\omega )  \label{ 6.2.1}\end{equation} is a quasi-isomorphism.
Together with  isomorphism  $r_{L,-t^{-1}dt}:
\CE_{\text{dR}}(M)_{(b,-t^{-1}dt)}\iso\det \CC (L,  L_\omega )$
from(2.5.6), it yields an identification $e(L): \CE_{\text{dR}}
(M)_{(b,-t^{-1}dt)}\iso \CE_{\text{B}}(M)_{(b,-t^{-1}dt)}$. Thus
$\rho^\varepsilon_{(b,-t^{-1}dt)} =\gamma_{\rho^\varepsilon}^* (L)
e(L)$ for some $\gamma^*_{\rho^\varepsilon} (L)\in\Bbb C^\times$.

\bigbreak\noindent\textbf{Proposition.}\emph{ One has
$\gamma_{\rho^\varepsilon}^* (L)=\mathop\Pi\limits_{\lambda
\in\Lambda (L)}\gamma_{\rho^\varepsilon}^* (\lambda)$, and for
$\rho^\varepsilon$ as in (6.1.2)
\begin{equation}\gamma^*_{\rho^\varepsilon}(\lambda ) = (2\pi
)^{-1/2}   \exp (\pi i (\lambda -1/2) )\Gamma (\lambda ) . \label{
6.2.2}\end{equation}}\bigbreak

{\it Proof.} If $L'\supset L$ is another lattice as above, then
$e(L') /e(L)$ is the determinant of the action of $-t\partial_t$ on
$tL'/tL$, so the validity of the assertion does not depend on the
choice of $L$. It is compatible with filtrations, and holds for $M$
supported at $b$, so we can assume that $M$ has rank 1 and is the
$*$-extension at $b$. Thus $\Lambda (L)=\{ \lambda \}$; by
continuity, it suffices to consider the case of $\lambda\notin \Bbb
Z$. Then the complexes in (6.2.1) are acyclic.  The corresponding
trivializations of  $\CE_{\text{dR}} (M)_{(b,-t^{-1}dt)}$ and
$\CE_{\text{B}}(M)_{(b,-t^{-1}dt)}$ are $\iota (L)_{-t^{-1}dt}$ from
6.1 and $1^*_b$ from 4.7; by construction,
$e(L)(\iota(L)_{-t^{-1}dt} )=1^*_b$.

By (2.11.3), the counterclockwise monodromy from $t^{-1}dt$ to
$-t^{-1}dt$ sends $\iota (L)_{t^{-1}dt}$ to $\exp (\pi i (\lambda
-1/2))\iota (L)_{-t^{-1}dt}$. According to 4.8, the same monodromy
sends $1^!_b$ to $(1-\exp (-2\pi i \lambda ))^{-1} 1^*_b$. Since
$\rho^\varepsilon$ is horizontal, one has
$\gamma^*_{\rho^\varepsilon }(\lambda )=(1-\exp (-2\pi i \lambda
))^{-1}$ $\gamma^!_{\rho^\varepsilon }(\lambda )\exp (\pi i (1/2
-\lambda ))= \exp (\pi i(\lambda - 1/2))(1-\exp (2\pi i \lambda
))^{-1}  \gamma^!_{\rho^\varepsilon }(\lambda )$, and we are done by
(6.1.2).  \hfill$\square$

{\it Example.}  If $M$ is smooth at $b$ and  $L=tM$, then
isomorphism (6.2.1) is $\Gamma (X_b , M^\nabla )\iso M_b$, $m\mapsto
m_b$, and $\gamma^*_{\rho^\varepsilon} (tM) =((2\pi )^{-\frac{1}{2}}
i )^{\text{rk} (M)}e(L)$.

{\it Exercise.} Deduce  Euler's reflection formula $\Gamma
(\lambda)\Gamma (1-\lambda)= \pi \sin^{-1} (\pi \lambda )$ from
(6.1.2), (6.2.2), the lemma in 2.7, and Exercise in 4.7.

\subsection{}\label{6.3.} Let us write down a formula for the factors $[\rho^\varepsilon_{(O,\nu)}]$ from 0.3.

Recall that we have $X$, $M$ and $\nu$ defined over a subfield $k$
of $\Bbb C$, and $B(M)$ is defined over a subfield $k'$. The  finite
set of singular points of $M$ and $\nu$ is defined then over $k$; it
is partitioned by $\Aut (\Bbb C/k)$-orbits.  Let $O$ be such an
orbit. The $\Bbb C$-line $\CE_{\text{dR}}(M)_{(O,\nu
)}=\mathop\otimes\limits_{x\in O}\CE_{\text{dR}}(M)_{(x,\nu )}$ is
defined  over $k$ by \S 2,\footnote{The group $\Aut (\Bbb C /k )$
acts on $\CE_{\text{dR}}(M)_{(O,\nu )}$  by transport of structure;
its fixed points is the $k$-structure on $\CE_{\text{dR}}(M)_{(O,\nu
)}$.} and  $\CE_{\text{B}}(M)_{(O,\nu )}=
\mathop\otimes\limits_{x\in O}\CE_{\text{B}}(M)_{(x,\nu )}$ is
defined over $k'$ by \S 4.    Computing $\rho^\varepsilon :
\CE_{\text{dR}}(M)_{(O,\nu )}\iso \CE_{\text{B}}(M)_{(O,\nu )}$,
$\rho^\varepsilon \in E_{\text{B}/\text{dR}}$, in  $k$- and
$k'$-bases, we get a number whose class  $[\rho^\varepsilon_{(O,\nu
)}]$  in $ \Bbb C^\times /k^{\prime \times} k^\times$ does not
depend on the choice of the bases and the choice of
$\rho^\varepsilon$ in $E_{\text{B}/\text{dR}}$. Let us compute
$[\rho^\varepsilon_{(O,\nu)}]$ explicitly assuming that $M$ has
regular singularities.

For $b\in O$ let $k_b \subset \Bbb C$ be its field of definition;
let $X_b$ be a small disc around $b$. Choose an auxiliary datum on
the de Rham side: it is $(t, L, u, v)$, where $t$ is a parameter at
$b$, $L$ is a $t\partial_t$-invariant $b$-lattice in $M$, $u$ is a
non-zero vector in $\det L_b$, and $v$ is a non-zero vector in $\det
\CC (L,L_\omega )$ (see 6.2); we assume that $(t, L, u, v)$ are
defined over $k_b$.  Let $\Lambda_b$ be the spectrum  (with
multiplicities) of $t\partial_t$ acting on the fiber $L_x$; we
assume that $\Lambda_b$ does not contain non-positive integers. An
auxiliary datum on the Betti side is  $( \phi , w)$, where $\phi$ is
a non-zero horizontal section of $\det M$ over the half-disc
Re$(t)>0$, which is defined over $k'$ (with respect to the Betti
$k'$-structure on the sheaf of horizontal sections), $w$ is a
non-zero vector in $\det R\Gamma_{\!\text{dR}}(X_b , M)$ defined over
$k'$.

The data yield numbers: The leading term of $\nu$ at $b$ is
$\alpha_b t^{-\ell} dt $, $\alpha_b \in k^\times_b$; let $r_b \in
k_b$ be the trace of $t\partial_t$ acting on $L_b$. Notice that
$m_b := \exp (- 2\pi r_b )$ is the monodromy of $\det M_b^\nabla $
around $b$, so $m_b\in k^{\prime\times}$. Then the section $t^{r_b}
\phi$ on the half-disc extends to an invertible holomorphic section
of $\det L$ on $X_b$; set $\beta_b := (t^{r_b} \phi )_b /u \in\Bbb
C^\times$. Let
 $\delta_b \in \Bbb C^\times$ be the ratio of $v$ and the image of $w$ by the determinant of (6.2.1).

Let us compute the numbers $\alpha_b$, $\beta_b$, $\delta_b$ and the
spectrum $\Lambda_b$ for each $b\in O$ using Galois-conjugate de
Rham side data. Set $n:= \text{rk}(M)$.

\bigbreak\noindent\textbf{Proposition.} \emph{One has
\begin{equation} [\rho^\varepsilon_{(O,\nu )}]= \mathop\Pi\limits_{b\in O} (2\pi
)^{-\frac{n\ell }{2}} i^{n\frac{\ell (\ell -1 )}{2}} m_b^{\frac{\ell
-1}{2}}\alpha_b^{\frac{n\ell }{2}-r_b} \beta_b^{\ell
-1}\delta_b\mathop\Pi\limits_{\lambda \in\Lambda_b} \Gamma (\lambda
). \label{6.3.1}\end{equation}}\bigbreak

{\it Proof.} For the sake of clarity, we do the computation assuming
that $b$ is a $k$-point, leaving the general case to the reader.

As follows from (2.11.3), the validity of formula does not depend on
$\alpha_b$. Notice that the class of $\alpha_b^{-r_b}:= \exp (-r_b
\log (\alpha_b ))$ in $\Bbb C^\times /k^{\prime\times}$ is well
defined: adding $2\pi i$ to the logarithm multiplies the exponent by
$m_b \in k^{\prime\times}$.

 If $\ell =1$ and $\alpha_b =-1$, then the formula follows from (6.2.2).

To finish the proof, it remains to check that
$$[\rho^\varepsilon_{(b,-t^{-\ell -1}dt )}]= (2\pi
)^{-\frac{n}{2}}i^n \beta [\rho^\varepsilon_{(b,t^{-\ell }dt )}].$$
Consider a family of forms $\nu_x := t^{-\ell}(x-t)^{-1}dt$. Then
$\CE_{\text{dR}}(M)_{(x, 1_b , \nu_x )}=\mu^\nabla (t^{\ell
-1}(t-x)L/L_\omega )= \det \CC (L,L_\omega )\otimes \lambda
(L/t^{\ell -1}(x-t)L)^{\otimes -1}=  \CC (L,L_\omega )\otimes
\lambda (L/t^{\ell -1} L)^{\otimes -1}\otimes \lambda (t^{\ell
-1}L/t^{\ell -1}(x-t)L)^{\otimes -1}$. We fix a non-zero $l$ in
$\lambda (L/t^{\ell -1}L )$ defined over $k$. Any local
trivialization $g$ of $\det L$ yields then a trivialization $e (g)_x
:= v\otimes l^{-1}\otimes t^{n(1-\ell )}g^{-1}_x$  of
$\CE_{\text{dR}}(M)_{(x, 1_b , \nu_x )}$; if $g$ is defined over
$k$, then so is $e (g)$.

The leading terms of $\nu_x$ at $t=0$ and $t=x$ are
$x^{-1}t^{-\ell}dt$ and $x^{-\ell}(x-t)^{-1}dt$. Applying (2.11.2)
to the $x$-lattice $(x-t)L$ and (2.11.3) to the $b$-lattice $t^{\ell
-1}L$, we see that $e(t^{r_b} \phi )_x =v\otimes ( x^{-\frac{ n(\ell
-2) }{2}  - r_b} \ell^{-1} )\otimes ( x^{-\frac {n\ell
}{2}}\phi^{-1})$ is a horizontal (with respect to $x$) section of
$\CE_{\text{dR}}(M)_{(x, 1_b , \nu_x )}$. Since the value at $b$ of
$\beta_b t^{n(1-\ell )} (t^{r_b}\phi )^{-1}$ is a generator of $\det
(t^{\ell -1}L/t^\ell L)^{\otimes -1}$ defined over $k$, we see that
$\beta_b e(t^{r_b} \phi )_b \in \CE_{\text{dR}}(M)_{(b, -t^{-\ell
-1}dt)}   $ is defined over $k$. If $s$ a horizontal section of
$\CE_{\text{B}}(M)_{(x,1_b ,\nu_x )}$ over $X_b$ defined over $k'$,
then $ \rho^\varepsilon (e(t^{r_b} \phi )) /s$ is a constant
function. Its value at $x=0$, i.e., at $b$, belongs to
$\beta_b^{-1}[\rho^\varepsilon_{(b,-t^{-\ell -1}dt )}]$. By
factorization and  Example in 6.2, its value   at $x=1$ belongs to
$[\rho^\varepsilon_{(b,t^{-\ell }dt )}] (2\pi )^{-\frac{n}{2}}i^n$,
and we are done. \hfill$\square$

 \bigskip

{\bf Notation.}   $a_\psi$ 2.9; $\CC (U)$ 4.2; $\CC (W,\CN )$ 4.4;
$\CC (L,L_\omega )$ 2.4; $\fD$, $\fD^\diamond$, $(D,c,\nu_P )$,
$(D,c,\nu )$ 1.1; $|D|$ 1.1; $dR (L,L_\omega )$ 2.7; $\CD
et_{P/S}(E)$, $\CD et_{P/S}(E_1/E_2 )$ 2.3;  $\Div (X)$ 1.1; ${}^w
\Bbb E$, ${}^w \Bbb E^0$ 5.1;  $\Bbb E$, $\Bbb E^0$ 5.2; $\CE$ 1.2;
$\CE^{(\ell )}$ 1.6; $E_{\text{dR} /\text{B}}$ 5.4; $e$, $e_L$ 2.5;
$\CE_B (M)$ 4.6; $\CE_{\text{dR}}(M)$ 2.5; $\CI (U,\CN) $ 4.3;
$G^\flat$ 2.2; $Hom_{\CL}(J_1 ,J_2 )$ 2.2; $j_{T*}M$ 2.4; $K(D)$,
$K(D)^\times_{D,c}$ 1.1; $\CK^{(\ell )}$ 1.1; $K^\times_{T_b}$ 1.1;
$\CL_?$, $\CL_k$, $\CL_\CO$, $\CL_{\text{dR}}$ 1.2;
$\CL^{\text{inv}}_{\text{dR}} (K^\times_{T_b})$ 1.11;
$\CL^0_{\text{dR}}$ 1.7; $\CL_{\text{dR}} (X,T)$, $\CL_{\text{dR}}
(X,T)^{(\ell )} $, $\CL_{\text{dR}} (X\setminus T)^{(\ell )} $  1.6;
$\CL^\Phi_?$ 1.5; $\CL^\natural_{\text{dR}} (X,T;K)$ 1.11;
$\CL^\Phi_{\text{dR}} (X,T)^{\CO\text{-triv}}$, $\CL^\Phi_{?}
(X,T;K)_T$ 1.13; $\CM (\CC (U))$ 4.1; $O^\times_{D,c}$ 1.1;
$O^\times_{T_b}$ 1.1; $P_{D,c}$ 1.1; Rat 1.4; $r_{L,\nu} $ 2.5;
$\CV_{\text{crys}} $ 1.2; $\CV^{\text{sm}}$ 1.3; $\tau_\psi$ 2.9;
$\tau_\CN$ 4.3; $T$, $T^c_S$ 1.1; $\gamma_{\rho^\varepsilon}^* (L)$,
$\gamma_{\rho^\varepsilon}^* (\lambda )$ 6.2;
$\gamma_{\rho^\varepsilon}^! (L)$, $\gamma_{\rho^\varepsilon}^!
(\lambda )$ 6.1;  $\lambda_P$ 2.3; $\lambda (F)$ 4.2; $\Lambda (L)$
6.1; $\mu^\nabla_P$ 2.4;  $\mu^\nu_P$ 2.5; $\mu^\psi$ 2.9;
$\mu^{\nabla/\kappa}_P$, $\mu^{\kappa/\nu}_P$ 2.10; $\pi_0 (X)$,
$\pi_0 (\CL )$, $\pi_1 (\CL )$ 1.1;   $\phi_\kappa$ 2.10;
$\nabla^\varepsilon$ 2.11; $\omega (X,T)$ 1.12; $\Omega^1
(K^\times_T )^{\text{inv}}$ 1.10;  $\tau_\psi$ 2.9; $1^!_b$, $1^*_b$
4.7; $2^T$ 1.1.



\begin{thebibliography}{BDE}

\bibitem[A]{A}
G.~Anderson, {\em Local factorization of determinants of twisted DR
cohomology groups}, Compositio Math. \textbf{83} (1992), no.~1,
69--105.

\bibitem[B]{B}
A.~Beilinson, {\em Topological $\CE$-factors}, Pure Appl.~Math.~Q.
\textbf{3} (2007), no.~1, 357--391.

\bibitem[BBE]{BBE}
A.~Beilinson, S.~Bloch, H.~Esnault, {\em $\CE$-factors for
Gau\ss-Manin determinants}, Moscow Mathematical Journal \textbf{2}
(2002), no.~3, 477--532.


\bibitem[BD]{BD}
 A.~Beilinson, V.~Drinfeld, {\em Quantization of Hitchin's
integrable system and Hecke eigensheaves}, {\tt
http://www.math.uchicago.edu/$\sim$mitya/langlands.html}

\bibitem[BG]{BG}
A.~Beilinson, D.~Gaitsgory, {\em A corollary of the $b$-function
lemma}, {\tt math.AG 0810.1504} (2008).

\bibitem[Bj]{Bj}
J.-E.~Bj\"ork, {\em Analytic $\CD$-modules and applications},
Mathematics and its Applications, Kluwer Academic Publishers,
Dordrecht-Boston-London, 1993.

\bibitem[BDE]{BDE}
S.~Bloch, P.~Deligne, H.~Esnault, {\em Periods for irregular
connections on curves} (2005).

\bibitem[BE]{BE}
S.~Bloch, H.~Esnault, {\em Homology for irregular connections},
J.~Th\'eor.~Nombres~Bordeaux \textbf{16} (2004), no.~2, 357--371.

\bibitem[CC]{CC}
C.~Contou-Carr\`ere, {\em Jacobienne locale, groupe de bivecteurs de
Witt universel, et symbole mod\'er\'e},  C.~R.~Acad.~Sci.~Paris
S\'er.~I Math. \textbf{318} (1994), no.~8, 743--746.

\bibitem[Del]{Del}
P.~Deligne, {\em Seminar at IHES, Spring 1984, handwritten notes by
G.~Laumon}, \newline {\tt
http://www.math.uchicago.edu/$\sim$mitya/langlands.html}

\bibitem[Den]{Den}
C.~Deninger, {\em Local L-factors of motives and regularized
determinants}, Invent. Math. \textbf{107} (1992), 135--150.

\bibitem[Dr]{Dr}
V.~Drinfeld, {\em Infinite-dimensional vector bundles in algebraic
geometry}, The Unity of Mathematics. In Honor of the 90th
Birthday of I.M.~Gelfand. Progress in Mathematics, vol.~244,
Birkh\"auser, Boston-Basel-Berlin, 2006, pp.~263--304.

\bibitem[E]{E}
H.~Esnault, {\em Talk at the Tokyo conference  ``Ramification and
vanishing cycles", 2007}, \newline {\tt
http://www.ms.u-tokyo.ac.jp/$\sim$t-saito/conf/rv/rv.html}

\bibitem[G]{G}
W.~Goldman, {\em The complex-symplectic geometry of SL$(2,\Bbb
C)$-characters over surfaces}, Algebraic groups and arithmetic, Tata
Inst.~Fund.~Res., Mumbai, 2004, pp.~375--407.

\bibitem[Gr1]{Gr1}
A.~Grothendieck, {\em Sur certains espaces de fonctions holomorphes,
I,II},  J.~Reine Angew. Math. \textbf{192} (1953), 35--64, 77--95.

\bibitem[Gr2]{Gr2}
A.~Grothendieck, {\em La th\'eorie de Fredholm},  Bulletin de la
S.M.F. \textbf{84} (1956), 319--384.

\bibitem[I]{I}
R.~Ishimura, {\em Homomorphismes du faisceau des germes de fonctions
holomorphes dans lui-m\^eme et op\'erateurs diff\'erentiels},
Memoirs of the Faculty of Science, Kyushu University \textbf{32}
(1978), 301--312.

\bibitem[L]{L}
G.~Laumon, {\em Transformation de Fourier, constantes d'equations
fonctionelles et conjecture de Weil},  Publ.~Math.~IHES \textbf{65}
(1987), 131--210.


\bibitem[M]{M}
B.~Malgrange, {\em \'Equations diff\'erentielles \`a coefficients
polynomiaux}, Progress in Mathematics, vol.~96, Birkh\"auser,
Boston, MA, 1991.

\bibitem[LS]{LS}
F.~Loeser, C.~Sabbah, {\em \'Equations aux diff\'erences finies et
d\'eterminants d'int\'egrales de fonctions multiformes},
Comment.~Math.~Helv. \textbf{66} (1991), no.~3, 458--503.

\bibitem[Me]{Me}
Z.~Mebkhout, {\em Une \'equivalence de cat\'egories}, Compositio
Math. \textbf{51} (1984), no.~1, 51--62.

\bibitem[P]{P}
D.~Patel, {\em Thesis},  University of Chicago (2008).


\bibitem[PX1]{PX1}
D.~Pickrell, E.~Xia, {\em Ergodicity of mapping class group actions
on representation varieties, I. Closed surfaces},
Comment.~Math.~Helv. \textbf{77} (2002), 339--362


\bibitem[PX2]{PX2}
D.~Pickrell, E.~Xia, {\em Ergodicity of mapping class group actions
on representation varieties, II. Surfaces with boundary},
Transformation Groups \textbf{8} (2003), no.~4, 397--402.

\bibitem[PS]{PS}
 A.~Pressley, G.~Segal, {\em Loop groups}, Oxford Mathematical Monographs, The
Clarendon Press, Oxford University Press, New York, 1986.

\bibitem[PSch]{PSch}
 F.~Prosmans, J.-P.~Schneiders,
{\em A topological reconstruction theorem for $\CD^\infty$-modules},
Duke Math~J. \textbf{102} (2000), no.~1, 39--86.



\bibitem[ST]{ST}
T.~Saito, T.~Terasoma, {\em Determinant of period integrals},
J.~Amer.~Math.~Soc. \textbf{10}, (1997), no.~4, 865--937.

\bibitem[SW]{SW}
G.~Segal, G.~Wilson, {\em Loop groups and equations of KdV type},
Publ.~Math.~IHES \textbf{61} (1985), 5--65.

\bibitem[T]{T}
T.~Terasoma, {\em A product formula for period integrals},
Math.~Ann. \textbf{298} (1994), 577--589.


\end{thebibliography}
\end{document}